\documentclass[10pt,letterpaper,reqno]{article}

\usepackage{epsf}
\usepackage{amsmath}
\usepackage{amsfonts}
\usepackage{amssymb}
\usepackage{amstext}
\usepackage{amsbsy}
\usepackage{amsthm}
\usepackage{amscd}
\usepackage[all]{xy}

\newcommand{\eeight}{\left( E_8 \times E_8 \right ) \rtimes \mathbb{Z}_2}
\newcommand{\spintt}{{\rm Spin}(32)/\mathbb{Z}_2}
\newcommand{\mhet}{\mathcal{M}_{{\rm het}}^G}

\newcommand{\m}{\mathcal}
\newcommand{\w}{\widetilde}
\newcommand{\mc}{\theta_G}

\newcommand{\Zee}{\mathbb{Z}}
\newcommand{\Cee}{\mathbb{C}}

\newcommand{\psl}{{\rm PSL}(2, \mathbb{Z})}
\newcommand{\abcd}{\left (
\begin{array}{cc}
a & b \\
c & d \\
\end{array}
\right)}

\newtheorem{teo}{Theorem}
\newtheorem{prop}[teo]{Proposition}
\newtheorem{coro}[teo]{Corollary}
\newtheorem{rem}[teo]{Remark}
\newtheorem{dfn}[teo]{Definition}

\newtheorem{lem}[teo]{Lemma}
\newtheorem{cl}[teo]{Claim}

\title{Heterotic String Data and Theta Functions} 
\author{Adrian Clingher
\thanks{Department of Mathematics, Stanford University, Stanford CA 94304.
{\tt Email:{\bf clingher@math.stanford.edu}}
}
}

\begin{document}
\maketitle
\begin{center} 
\abstract{\noindent We use the language of differential cohomology to give an analytic description of the moduli space of classical vacua for heterotic string theory in eight dimensions. The complex structure of this moduli space is then related to the character of the appropriate Kac-Moody algebra.}
\end{center}
%
%
%
%
%
\section{Introduction}
\label{fir}
The general framework of heterotic string theory \cite{narain1} \cite{narain2} \cite{wgs} 
involves the following defining geometric data:
\begin{enumerate}
\item a space-time which is a smooth spin manifold $X$, and
\item a principal $G$-bundle $P\to X$.
\end{enumerate}
The background fields $ \m{F} $ consist of triplets $ \left ( g, A, B \right ) $ consisting of:
\begin{enumerate}
\item a riemannian metric $g$ on $X$,
\item a connection $A$ on $P$,
\item a $B$-field $B$, which is, at least locally, a two-form on $X$.
\end{enumerate}
The first conditions required by cancellation of anomalies are topological. The space-time $X$ has to be 
ten-dimensional. The Lie group $G$ must be
$ \eeight $ or $ \spintt  $. 
\par The three background fields are to be considered up to gauge 
equivalence. This means that the connection is considered up to bundle isomorphism, that the metric 
is considered up to diffeomorphism equivalence, and that the two-form $B$ is considered up to 
equivalence by adding exact two-forms.  One of the major problems at this point has been 
to understand exactly the mathematical nature of the B-field. In the original physics 
interpretation \cite{narain2}, the B-field $B$ is 
manipulated as a globally defined two-form on $X $ with its contribution to the 
world-sheet path-integral being given by the phase:
\begin{equation}
\label{bbbb111}
{\rm exp} \left ( i \int_{\Sigma} x^*B\right ). 
\end{equation}
However, this definition is not fully satisfactory. First, $B$ in this formulation lacks a gauge
invariance. Second, the B-field interpretation in heterotic string theory 
must obey the world-sheet anomaly cancellation mechanism, as explained for instance by Witten in 
\cite{MR1748791}. This requires, among other conditions, that the path-integral contribution 
$ (\ref{bbbb111}) $ of $B$ is not a complex number of unit length but rather an element in the total space
of a non-trivial circle bundle. In addition, the field strength $H_B $, which is a globally defined three-form, 
has to satisfy the equation.     
\begin{equation}
\label{annc}
d H_B \ = \ \frac{1}{4 \pi c_G} \ {\rm Tr} \left ( F_A \wedge F_A \right ) 
- \frac{1}{4\pi} \ {\rm Tr} \left ( F_g \wedge F_g \right ) .
\end{equation}    
Here $ F_A $ and $F_g $ represent the curvature forms of the connection $A$ and the Levi-Civita
connection associated to the metric $g$ in $ TX $, respectively. The number $c_G $ is the dual 
Coxeter number of $G$. In the naive 
interpretation of $B$ as a globally defined two-form, $H_B = dB $. In particular 
$ H_B $ is a closed 3-form and cannot satisfy $ (\ref{annc}) $. 
\par The above issues suggest that the B-field is not in fact 
a globally defined two-form but rather must be given by a more subtle object, which locally resembles 
a two-form but has a non-trivial global structure. In the recent years, it has been noted by a number of
authors (see for instance \cite{MR1748791} and \cite{f1}) that, in order to satisfy the above conditions, 
the B-field has to be understood within a gerbe-like formalism. 
In this paper we follow ideas developed by Freed in \cite{f1} and Hopkins and Singer in \cite{hopkins}, and 
interpret the B-field as a non-flat 2-cochain in differential 
cohomology\footnote{In the general heterotic string setting, this approach introduces the B-field as an 
object in differential KO-theory. However, when one compactifies over 
the two-torus, as is the case we shall be dealing with here, 
the topological type of the space-time is sufficiently low-dimensional that one can 
ignore the differences between KO-theory and standard cohomology.}. 
These cochains appear as a mixture of Cech cocycles and 
local differential forms and can be regarded, in an open covering $ U_a $ of $X$, as multiplets:
$$ B \ = \ \left ( \ H, \ \omega^2_a, \ \omega^1_{ab} , \ \omega^0_{abc} , \ \omega^{-1}_{abcd} \ 
\right )$$ 
with $H $ being a globally defined three-form (the field strength $H_B$) and $ \omega^i $ local
i-forms (constant functions to $2\pi \Zee $ if $i=-1$). The anomaly cancellation condition 
appears in this formulation as:
\begin{equation}
\label{cs}
\check{d} B \ = \ \m{CS}_A - \m{CS}_g 
\end{equation}
where $ \check{d} $ represents the differentiation operator in differential cohomology
and $ \m{CS}_A$ and $ \m{CS}_g$ are Chern-Simons terms regarded as non-flat 3-cocycles.
One can introduce then a concept of gauge invariance for B-fields. Two B-fields are said 
to be gauge equivalent if the two differential 2-cochains underlying them differ by a coboundary of a flat differential 
1-cochain. It follows that the globally defined 3-form $H_B$, the field strength of $B$,
is gauge invariant. However, it is not necessarily true that $H_B $ is closed. In fact, 
one recovers easily the expected equation $ (\ref{annc}) $ by just 
writing the anomaly cancellation condition $ (\ref{cs}) $ at the field strength level.      
\par Moreover, apart from their rather technical construction, the B-fields as 
differential 2-cochains carry a very nice geometrical interpretation. They
represent sections in a two-gerbe \cite{hit} and carry holonomy denoted:
\begin{equation}
\label{bhol}
{\rm exp} \left ( i  
\int_{\Sigma} x^* B 
\right ) 
\end{equation}
along any map $ x \colon \Sigma \rightarrow X $ in the same way circle bundle connections carry
holonomy along loops. This is the proper meaning of $ (\ref{bbbb111}) $, 
the B-field contribution to the world-sheet path-integral. However the quantity $ (\ref{bhol}) $ 
is not an unitary complex number but rather a point in the total space of a circle fibration, which
fits exactly the picture required anomaly cancellation (see \cite{MR1748791} for details). 
\par In this framework, the heterotic classical vacua are obtained by taking the triplets 
$ \left ( A, g, B \right ) $ satisfying the cancellation condition $ (\ref{cs}) $ and imposing
over them the classical Equations of Motion. These equations can be written as:     
\begin{equation}
\label{eqofmot}
F_A =0, \ \ \ {\rm Ric}(g) =0 , \ \ H_B = 0. 
\end{equation}
The moduli space of heterotic classical vacua $ \mhet $ is then constructed out of all 
solutions to $ (\ref{eqofmot}) $ modulo gauge equivalence. 
\par The goal of this paper is to work out a precise mathematical description of the moduli 
space of classical vacua for heterotic string theory when the space-time is compactified along a 
two-torus. Namely, $ X = \mathbb{R}^8 \times E $ with $E$ a two-torus and the fields are considered
up to Euclidean equivalence on the $ \mathbb{R}^8$ factor. The moduli space of quantum vacua associated
to this case has been described in \cite{narain1} and \cite{narain2} from a purely physical perspective.  
In that formulation, the B-field appears as a globally defined two-form $B$ and the metric and $B$ fit 
together to form the imaginary and, respectively, the real part of the so-called complexified Kahler 
class. The physical quantum features of the theory depend then on a certain lattice of momenta $ \m{L}_{(A,g,B)} $ 
described by Narain in \cite{narain1}. This is a rank $20$ lattice that lives inside a fixed ambient real space 
$ \mathbb{R}^{2,18} $, is well-defined up to $ O(2) \times O(18) $ rotations 
and is even and unimodular. The real group $ O(2,18) $ acts transitively on the set of all 
$ \m{L}_{(A,g,B)} $ and, in 
this regard, one can consider the physical momenta as being parameterized by the $36$-dimensional
real homogeneous space:
\begin{equation}
\label{coset} 
O(2,18) / O(2) \times O(18). 
\end{equation}         
One identifies then the points in $ (\ref{coset}) $ corresponding to equivalent quantum 
theories. This amounts to factoring out the left-action of the group $ \Gamma $ of integral 
isometries of the lattice. 
The quantum (Narain) moduli space of distinct heterotic string 
theories compactified on the two-torus appears as:
\begin{equation}
\label{narainmod}  
\m{M}^{\rm quantum}_{{\rm het}} \ = \ \Gamma \backslash O(2,18) / O(2) \times O(18).
\end{equation} 
\par From a geometric interest point of view, the above physics-inspired Narain construction
has a couple of disadvantages. First, since $(\ref{narainmod})$ provides a description of quantum states, not all
the identifications are accounted by classical geometry. As explained for example in \cite{morrison1}, part of the $\Gamma $-action models the so-called quantum corrections and results in
identifications of momenta for pairs of triplets $(A,g,B) $ which are not isomorphic. 
Second, the Narain construction does not provide a holomorphic description. Technically, one 
can endow the homogeneous quotient $ (\ref{narainmod}) $ with a complex structure, but holomorphic families 
of elliptic curves and flat connections do not embed as holomorphic subvarieties in  
$ \m{M}^{\rm quantum}_{{\rm het}} $ (see the Appendix of \cite{clingher} for an outline of this issue).     
\par In this paper, we take a completely geometrical approach and use the interpretation of the B-fields as differential cochains in order to give an analytic description for the moduli space of classical vacua $ \mhet $. 
In section $ \ref{sec4} $ we prove the following:
\begin{teo}  
\label{het1}
The moduli space of classical solutions $ \mhet $ can be given
the structure of a 18-dimensional complex variety with orbifold singularities. 
Moreover, $ \mhet $ represents the total space of a holomorphic
$ \mathbb{C}^* $-fibration
\begin{equation}
\label{fib}
\mhet \rightarrow  \m{M}_{E, G}
\end{equation} 
over the moduli space $ \m{M}_{E, G} $ of isomorphism classes of 
pairs of elliptic curves and flat $ G$-bundles, where $ G=\eeight $ in the
case of the $ E_8 \times E_8 $ heterotic theory and $ G = \spintt $ 
in the case of $ {\rm Spin}(32)/ \Zee_2 $ theory.  
\end{teo}
\noindent In this context, it is important then to analyze the holomorphic 
type of fibration $ (\ref{fib}) $. The structure of the base space $  \m{M}_{E, G} $ is 
well-known (see \cite{morgan1}). Given a fixed elliptic 
curve $E$, the family of equivalence classes of flat $G$-bundles (for G simply-connected) 
can be identified with the complex quotient:
$$ 
W  \backslash \left ( {\rm Pic}^o(E) \otimes_{\Zee} \Lambda_G \right )  
$$  
where $ \Lambda_G $ represents the coroot lattice of $G$ and $ W $ is the Weyl
group. Using a coordinate oriented model to track down the variation of the elliptic 
curve, one can regard  $ \m{M}_{E, G} $ as a complex orbifold:
\begin{equation}
\label{modelll}
\Pi_G \backslash V_G . 
\end{equation} 
where $ V_G = \m{H} \times \left ( \mathbb{C} \otimes_{\Zee} \Lambda_G \right ) $.
$ \m{H} $ represents the complex upper half-plane and $\Pi_G$ is a 
modular group acting on $ V_G $ by mixing as a semi-direct product the 
$ {\rm SL}(2, \Zee) $ action on
$ \m{H} $ and the affine Weyl group action on $ \mathbb{C} \otimes_{\Zee} \Lambda_G $.
\par The model $ (\ref{modelll}) $ allows us to analyze easily holomorphic 
$ \mathbb{C}^* $-fibrations over $ \m{M}_{E, G} $. Such fibrations over a complex 
orbifold are best described in terms of equivariant line bundles over the universal 
cover. Those are holomorphic line bundles  
$ L \rightarrow V_G  $
where the action of the modular group $ \Pi_G $ on the base is given a lift to the
fibers. All holomorphic line bundles over  
$ V_G   $
are trivializable and a lift of the action of $ \Pi_G $ to fibers can be obtained 
through a set of automorphy factors $ \left ( \varphi_a \right )_{a \in \Pi_G} $ 
with $ \varphi_a \in H^0( V_G , \ \m{O}^*_{V_G} ) $ satisfying:
$$ \varphi_{ab}(x) = \varphi_a(bx)  \cdot \varphi_b(x). $$ 
Such a set generates a class in the group cohomology 
$ H^1(\Pi_G, H^0( V_G , \ \m{O}^*_{V_G} )) $. Two automorphy factor sets provide 
isomorphic fibrations if and only if they determine the same group cohomology
class.  
\par There is one particular important holomorphic $\mathbb{C} $-fibration over
$ \m{M}_{E, G} $, the $\Lambda_G$-character fibration. This fibration can be  
defined using the model $ (\ref{modelll}) $ as the fibration supporting the 
$\Lambda_G$-character function:
\begin{equation}
\label{char1}
B_{\Lambda_G} \colon V_G \rightarrow \mathbb{C}, \ \ \ \ B_{\Lambda_G}
(\tau,  z) = 
\frac{\Theta_{\Lambda_G}(\tau,  z)}{\eta(\tau)^{16}}
\end{equation}
where $ \Theta_{\Lambda_G}(\tau,  z) $ represents the holomorphic theta-function
associated to the lattice $ \Lambda_G $ and $ \eta $ is Dedekind's eta function.
The holomorphic function $ B_{\Lambda_G} $ has an important interpretation. 
It represents (see \cite{kac1}) the zero-character of the level $ l=1 $ basic highest weight 
representation of the Kac-Moody algebra associated to $G$. Under the action of 
the modular group, $ B_{\Lambda_G} $ descends to a section in a non-trivial $ \mathbb{C} $-fibration:
\begin{equation}
\label{character}
\m{Z} \rightarrow  \m{M}_{E, G} 
\end{equation} 
We call this the $\Lambda_G$-character fibration. Then, as we prove in section $\ref{prooof}$, one has: 
\begin{teo}
\label{het2}
The heterotic $ \mathbb{C}^* $-fibration $ (\ref{fib}) $ can be holomorphically identified with 
the $\mathbb{C}^*$-fibration induced by the $\Lambda_G$-character fibration $ (\ref{character}) $.  
\end{teo}
\noindent One can then conclude that, in the light of Theorem $ \ref{het1} $, the moduli space 
$ \mhet $ of heterotic classical vacua compactified along the 2-torus can be seen naturally as 
the total space of the $\mathbb{C}^*$-fibration associated to $ (\ref{character}) $. 
\par As a final note, we comment on the relation between the space $\mhet$ that we construct here 
and the Narain moduli space $\m{M}^{\rm quantum}_{{\rm het}}$ of quantum vacua $(\ref{narainmod})$. 
It turns out that, there exists a diffeomorphism between an open neighborhood along the zero section 
of the zero-section in the $ \mathbb{C}^*$-fibration of $ \mhet$ and an open subset of 
$\m{M}^{\rm quantum}_{{\rm het}} $ on which the volume of the elliptic curve is large. That is precisely
the region where physics predicts the quantum corrections became insignificant and the moduli spaces
of classical and quantum vacua should coincide. This issue, as well as the place of $ \mhet$ in the
framework of the F-Theory/Heterotic String duality in eight dimensions, is the subject of a joint 
paper \cite{clingher} with John Morgan.  
\par This paper is a revised version of the author's Ph.D. thesis at Columbia University. The author 
would like to take this opportunity to thank his teacher there, John Morgan, for his assistance, 
encouragement and patience. 
%
%
%
%
\section{The B-field as a Differential Cochain}
\label{ch2}
This section introduces the framework needed for giving a precise 
definition of the B-field. We shall treat this object in the language of differential integral 
cohomology following ideas of Freed \cite{f1} and Hopkins-Singer \cite{hopkins}. 
This approach mixes together differential forms and integral cohomology and has its roots in earlier 
works of Deligne \cite{b1} and Cheeger-Simons \cite{cheeger}. However, as we shall see, the B-field is a 
cocycle rather than a class in this theory. 
\subsection{Differential Cohomology.}
\label{aa}
\par Let X be a smooth manifold. Choose $ \m{U} = \left ( U_i \right) _{i \in I} $, an open covering 
of $X$ indexed by an ordered set $I$. For $ r, s \in \Zee $, we set: 
$$ \check{C} ^{r,s} (\m{U}) =
\begin{cases}
0 & \text{if $ s<-1 $ or $ r<0 $ ; \ } \\
C^r(\m{U}, \Omega ^s _X ) & \text{if $ 0 \leq s $ \ ;} \\
C^r(\m{U}, 2 \pi \Zee ) & \text{if $ s=-1 $. } 
\end{cases} $$ 
where $ C^r(X, \Omega ^s _X ) $  represents the set of Cech r-cochains
with values on the sheaf of s-forms on X, $ \Omega ^s _X $:
$$ C^r(\mathcal{U}, \Omega ^s _X ) = \prod_{i_0, i_1, .... i_r} 
{ \Omega^s_X( U_{i_0} \cap U_{i_1} \cap ... \cap U_{i_r} ) }. $$   
Following \cite{f1}, a ${\bf flat}$ ${\bf differential}$ ${\bf  n}$-${\bf cochain}$ is defined 
as an element of the direct sum:  
$$ \check{C} ^n (\m{U}) := \bigoplus _{r+s=n} \check{C} ^{r,s} (\m{U}) $$  
In other words, one can see a flat differential n-cochain as a multi-plet of 
differential Cech cochains: 
$$ 
\omega = ( \omega^n _{a} , \ \omega^{n-1}_{a_1 a_2} , \ 
\omega^{n-2}_{a_1 a_2 a_3}, \ ...... \ \omega^{-1}_{a_1 a_2 ....a_{n+2}} ) 
$$ 
There are two derivation operators acting on the above differential 
cochains. A vertical 
differentiation operator 
$ d \colon \check{C} ^{r,s} (\m{U}) \rightarrow \check{C} ^{r,s+1} (\m{U}) $
represents the usual differentiation on each form component 
(inclusion if $s=-1$). A second horizontal operator 
$ \delta \colon  \check{C} ^{r,s} (\m{U}) \rightarrow 
\check{C} ^{r+1,s} (\m{U}) $
is the Cech coboundary operator: 
$$ 
\delta(\omega^n)_{a_1, a_2, .... a_{n+2}} = \ 
\sum_{j=1}^{n+2} {(-1)^{j+1} \omega^n_{a_1...a_{j-1} a_{j+1} ...a_{n+2}}
 \ \vert _{U_{a_1} \cap .... \cap U_{a_{n+2}}} }.
$$ 
The flat cochains together with the two derivations build a double complex 
$ ( \check{C} ^{r,s} (\m{U}), \ d, \ \delta ) $. As usual in such situations
a total differentiation can be introduced: 
$$ \check{d} \colon \check{C} ^n (\m{U}) \rightarrow
\check{C} ^{n+1} (\m{U}) \ , \ \ \check{d}\vert_{\check{C}^{r,s}}= 
\delta \vert_{\check{C}^{r,s}} + (-1)^{r+1} d \vert_{\check{C}^{r,s}}. $$
$ ( \check{C}^*(\m{U}), \ \check{d} ) $ is a cochain complex, defining 
flat n-cocycles and flat n-coboundaries which in turn generate 
$ \m{U} $-cohomology groups $ \check{H}^n(\m{U}) $. 
\par Let $ \m{V} = (V_j)_{j \in J} $ be a refinement of the open covering 
$ \m{U} = (U_i)_{i \in I} $ and $ \sigma \colon J \rightarrow I $ be 
a subordination map such that $ V_j \subset U_{\sigma(j)} $ for any 
$ j \in J $. One then has the restriction homomorphism: 
$$ \sigma^* \colon \check{C}^n(\m{U}) \rightarrow \check{C}^n(\m{V}). $$
commuting with the two differentiation operators $ d $ and $ \delta $. 
This defines a morphism of cochain complexes inducing a homomorphism 
at cohomology level:
$$ \sigma^{\m{U}}_{\m{V}} \colon 
\check{H}^n(\m{U}) \rightarrow \check{H}^n(\m{V}). $$
As in Cech theory:
\begin{lem} 
\label{lemacucech}
The homomorphism $ \sigma^{\m{U}}_{\m{V}} $ depends only
on the open covering $ \m{U} $ and refinement $ \m{V} $ and not on the choice
of subordination map $ \sigma \colon J \rightarrow I $. Furthermore, 
$ \sigma^{\m{U}}_{\m{U}} $ is identity, and if $ \m{W} $ is a refinement of
$ \m{V} $ then $ \sigma^{\m{U}}_{\m{W}} = \sigma^{\m{V}}_{\m{W}} \circ 
\sigma^{\m{U}}_{\m{V}} $.
\end{lem} 
\begin{proof}
Assume $ \sigma, \sigma' \colon J \rightarrow I $ are two subordination
maps with $ V_j \subset U_{\sigma(j)} \cap U_{\sigma'(j)} $. We define a 
family of homomorphisms:
$$ k^q \colon \check{C}^q(\m{U}) \rightarrow \check{C}^{q-1}(\m{V}) $$
by the following pattern. If
$$ \omega \ = \ ( \omega^q _{i_1} , \ \omega^{q-1}_{i_1 i_2} , \ 
\omega^{q-2}_{i_1 i_2 i_3}, \ ...... \ \omega^{-1}_{i_1 i_2 ....i_{q+2}} ) 
\in  \check{C}^q(\m{U}) $$
then $ k^q(\omega) = \eta \in  \check{C}^{q-1}(\m{V}) $ where
$$ \eta \ = \ ( \eta^{q-1}_{j_1} , \ \eta^{q-2}_{j_1 j_2} , \ 
\eta^{q-3}_{j_1 j_2 j_3}, \ ...... \ \eta^{-1}_{j_1 j_2 ....j_{n+1}} ) 
\in  \check{C}^q(\m{U}) $$  
with components:
\begin{equation}
\label{comphomot}
\eta^{q-r}_{j_1 j_2 \cdots j_r} \ = \ \sum_{t=1}^{r} \ 
(-1)^{t+1} \  
\omega^{q-r}_{\sigma(j_1) \cdots \sigma(j_t) \sigma'(j_t) \sigma'(j_{t+1}) 
\cdots \sigma'(j_r) } \ \vert _{  
V_{j_1} \cap V_{j_1} \cap \cdots \cap V_{j_r} }.
\end{equation} 
We claim that: 
\begin{equation}
\label{alacech}
\check{d} k^{q} + k^{q+1} \check{d}  \ = \ \sigma^* - (\sigma')^*.   
\end{equation}
Indeed, for a multi-index $ (j) = ( j_1 j_2 \cdots j_r) $ one can write:
$$
\left ( \check{d} k (\omega) \right )_{(j)} \ = \ (-1)^{r} d k(\omega)_{(j)} 
+ \left [ \delta k(\omega) \right ] _{(j)} = 
$$
$$ 
= \ (-1)^r \ \sum_{t=1}^r \ (-1)^{t+1} \ d \omega_{\sigma(j_1) 
\cdots \sigma(j_t) \sigma'(j_t) \cdots \sigma'(j_r)} \ + \ 
\sum_{l=1}^{r} (-1)^{l+1} (k \omega)_{j_1 \cdots j_{l-1} j_{l+1} \cdots 
j_r } \ = 
$$
$$ 
= \ \ (-1)^r \ \sum_{t=1}^r \ (-1)^{t+1} \ d \omega_{\sigma(j_1) 
\cdots \sigma(j_t) \sigma'(j_t) \cdots \sigma'(j_r)} \ + 
$$
$$   
+ \ \sum_{l=1}^{r} \sum_{t=1}^{l-1} \ (-1)^{t+l} 
\omega_{\sigma(j_1) \cdots \sigma(j_t) \sigma'(j_t) \cdots \sigma'(j_{l-1})
\sigma'(j_{l+1}) \cdots \sigma'(j_r)) } \ +  
$$
$$
+  \sum_{l=1}^{r} \ \sum_{t=l}^{r-1} \ (-1)^{t+l} \ 
\omega_{\sigma(j_1) \cdots \sigma(j_{l-1}) \sigma(j_{l+1}) \cdots 
\sigma(j_{t-1}) \sigma'(j_{t-1}) \cdots \sigma'(j_r)) }.   
$$
Similarly, 
$$
\left [ k(\check{d} \omega ) \right ]_{(j)} \ = \ \sum_{t=1}^{r} \ 
(-1)^{t+1} \ \left ( \check{d} \omega \right )_{\sigma(j_1) \cdots 
\sigma(j_t) \sigma'(j_t) \cdots \sigma'(j_r)} \ = \     
$$
$$
=  (-1)^{r+1} \sum_{t=1}^{r} \ 
(-1)^{t+1} \ d \omega_{\sigma(j_1) \cdots 
\sigma(j_t) \sigma'(i_t) \cdots \sigma'(i_r)} + 
$$
$$
+ \sum_{t=1}^{r} \ \sum_{l=1}^{t-1} (-1)^{t+l} \ 
\omega_{\sigma(j_1) \cdots \sigma(j_{l-1}) \sigma(j_{l+1}) \cdots 
\sigma(j_{t}) \sigma'(j_{t}) \cdots \sigma'(j_r) } + 
$$ 
$$
+ \ \sum_{t=1}^{r} \ \left ( \ - \   
\omega_{\sigma(j_1) \cdots \sigma(j_{t-1}) \sigma'(j_t) \cdots 
\sigma'(j_r) } \ + \ 
\omega_{\sigma(j_1) \cdots \sigma(j_t) \sigma'(j_{t+1}) \cdots 
\sigma'(j_r) }
\right ) \ +  
$$
$$  
+ \ \sum_{t=1}^r \ \sum_{l=t+2}^{r+1} \ (-1)^{t+l} \ 
\omega_{\sigma(j_1) \cdots \sigma(j_{t}) \sigma'(j_{t}) \cdots 
\sigma'(j_{l-1}) \sigma'(i_{l+1}) \cdots  \sigma'(j_r) } \ . 
$$
After canceling the terms with opposite signs, one obtains:
$$ 
\left ( \check{d} k (\omega) \right )_{(j)} \ + \ 
\left [ k(\check{d} \omega ) \right ]_{(j)} \ = 
\omega_{\sigma(i_1) \sigma(i_2) \cdots \sigma(i_r) } \ - \ 
\omega_{\sigma'(i_1) \sigma'(i_2) \cdots \sigma'(i_r) }. 
$$ 
Relation $ (\ref{alacech}) $ shows that $ k^q $ represents a homotopy operator 
between the two
cochain complex morphisms $ (\sigma')^* $ and $ \sigma^* $. Hence,
they must induce identical morphisms at cohomology level. This proves 
the first part of lemma. The second part follows immediately.  
\end{proof}
\vspace{.1in}
\noindent Lemma $ \ref{lemacucech}$ shows that $ \left ( \check{H}^r(\mathcal{U}), \ 
 \sigma^{\m{U}}_{\m{V}} \right ) $ forms a direct system. One then 
takes the direct limit over all possible coverings, defining:
$$ \check{H}^r(X) = \operatornamewithlimits{lim}_{ {\rightarrow} 
\atop{ \mathcal{U} }} \ \left ( \check{H}^r(\mathcal{U}), \ 
 \sigma^{\m{U}}_{\m{V}} \right ). $$
These ${\bf flat}$ differential cohomology groups are known to form the 
smooth Deligne cohomology of $X$ \cite{b1}. 
\par The differential objects we wish to study and make use of are non-flat extensions 
of the above cochains. By definition, a ${\bf non}$-${\bf flat}$ differential cochain is a pair 
$ (H, \ \omega ) $ consisting of a 
global ($n+1$)-form $H \in \Omega^{n+1} (X)$ and a flat n-cochain $ \omega $. 
One represents such an object, in an open covering, as a multi-plet of 
form-valued Cech cochains: 
$$ \omega = ( H, \ \omega^n _{a} , \ \omega^{n-1}_{a_1 a_2} , \ 
\omega^{n-2}_{a_1 a_2 a_3}, \ ...... \ \omega^{-1}_{a_1 a_2 ....a_{n+2}} ). 
$$ 
The top global form H represents the ${\bf field \ strength}$ of $ \omega $. Let us denote 
the set of non-flat differential 
n-cochains defined over an open covering $ \m{U} $ by: 
$$ N\check{C}^n (\m{U}) = \check{C}^n (\m{U}) \times \Omega^{n+1}(X). $$ 
The differentiation operator $ \check{d} $ extends naturally to 
$ \check{C}^n (\m{U}) $ . 
As in the flat case, $ ( N\check{C}^* (\m{U}) , \ \check{d}) $ 
determines a cochain complex defining cocycles $ N \check{Z}^n(\m{U}) $ 
and coboundaries $ N \check{Z}^n(\m{U}) $. 
\par The discussion in Lemma $ \ref{lemacucech} $ can be immediately 
reformulated to the new context. 
For a refinement $ \m{V} $ of $ \m{U} $ with two distinct 
subordination maps 
$ \sigma , \sigma'  \colon J \rightarrow I $, the two induced 
restriction homomorphisms:
$$ \sigma^* , (\sigma')^* \colon \check{C}^n(\m{U}) \rightarrow 
\check{C}^n(\m{V}). $$ differ by: 
\begin{equation}
\check{d} k^{q} + k^{q+1} \check{d}  \ = \ \sigma^* - (\sigma')^* .   
\end{equation}
In above expression, $ k^* $ is the homotopy operator
$$ k^q \colon N \check{C}^q(\m{U}) \rightarrow N \check{C}^{q-1}(\m{V}) $$ 
defined by:
$$
k^q \left ( H, \ \omega^q _{i_1} , \ \omega^{q-1}_{i_1 i_2} , \ 
\omega^{q-2}_{i_1 i_2 i_3}, \ ...... , \ \omega^{-1}_{i_1 i_2 ....i_{q+2}} 
\right ) \ 
= \  \ \left ( 0, \ \eta^{q-1}_{j_1} , \ \eta^{q-2}_{j_1 j_2} , \ 
\eta^{q-3}_{j_1 j_2 j_3}, \ ...... , \ \eta^{-1}_{j_1 j_2 ....j_{n+1}} 
\right ) 
$$
with components $ \eta^{q-r}_{j_1 j_2 ....j_{r}}, \ 1 \leq r \leq q $ given
by formulas $ (\ref{comphomot}) $.   
In a similar manner to the flat case, projection
$$ \sigma^{\m{U}}_{\m{V}} \colon N \check{H}^q(\m{U}) \rightarrow
N \check{H}^q(\m{V}) $$
does not depend on the choice of subordination assignment.  
Differential non-flat cohomology groups can then be defined: 
\begin{equation}
\label{nonflatch}
N\check{H}^*(X) =\operatornamewithlimits{lim}_{ {\rightarrow} 
\atop{ \mathcal{U} }} \ \left ( N \check{H}^r(\mathcal{U}), \ 
\sigma^{\m{U}}_{\m{V}} \right )  . 
\end{equation}
However, we are more interested here in non-flat cochains as explicit 
differential objects rather than as a framework for the above non-flat 
cohomology groups. In fact, one can show that the non-flat cohomology 
$ (\ref{nonflatch}) $ recovers the Cech cohomology of $X$ with coefficients
in the sheaf of smooth functions to the circle, 
$ N \check{H}^n(X) \simeq H_{{\rm Cech}}^n(X, \ S^1) $. 
\par Let us denote the set of non-flat differential cochains on $X$ by 
$$ N \check{C}^n(X) \ = \ \bigcup_{\m{U}} \ N  \check{C}^n(\m{U}).  $$ 
Among the cochains  in $ N \check{C}^n(X) $ we give special 
consideration to those associated to global differential 
forms. Let $ T \in \Omega^n(X) $ be such a $n$-form. 
In a given open covering, $ \m{U} = (U_i)_{i \in I} $, 
one can consider a differential cocycle with 
only two non-vanishing components as follows:
\begin{equation}
\label{formm} 
\left ( d T , \ T \vert _{U_i}, \ 0, \ 0, \cdots  \right ). 
\end{equation}   
This makes a non-flat differential n-cocycle. 
In future considerations we shall use differential forms in cocycle-like
computations, always interpreting the form as in $ (\ref{formm}) $. 
\par There exists a natural equivalence relation on the set of non-flat differential
n-cocycles.  
\begin{dfn}
Two cocycles $ \omega_1 \in N \check{Z}^n(\m{U}_1)$ and $ 
\omega_2 \in N \check{Z}^n(\m{U}_2) $ are said to be equivalent 
(denoted $ \omega_1 \sim \omega_2 $) if there exists a common refinement $ \m{U} $ and
corresponding subordination maps such that the difference of the two 
restrictions of $ \omega_1 $ and $ \omega_2 $ on $ N \check{Z}^n(\m{U}) $ is a 
coboundary of a differential flat (n-1)-cochain on $ \m{U} $. 
\end{dfn}
\noindent Clearly, two flat n-cocycles are equivalent if and only if they determine the same 
class in $ \check{H}^n(X) $. The relation $ \sim $ does not extend as equivalence relation
for general non-flat n-cochains. That is because there exist n-cochains  
$ \omega \in N \check{C}^n(\m{U}) $ which restricted on a refinement 
$ \m{U}' $ through two distinct subordination maps give 
restrictions $ \omega' , \omega'' \in N \check{C}^n(\m{U}') $ such that
the difference
$$ \omega' - \omega'' = \check{d} k^n \omega + k^{n+1} \check{d} \omega $$
is not necessarily a coboundary of a flat $(n-1)$-cochain.
However, relation $ \sim $ is still an equivalence relation on 
the subset $ \m{B}^n(X) \subset N \check{C}^n(X) $ defined as:    
\begin{equation}
\label{impdef}
\m{B}^n(X) \ = \ \left \{ \omega \in N \check{C}^n(\m{U}) \ \vert \ \m{U} \ 
{\rm open \ covering}, \ 
\check{d}\omega \in \Omega^{n+1}(X) \right \}. 
\end{equation}
\subsection{Geometrical Interpretation}
\label{sec222}
The differential cochains introduced above have a nice geometrical interpretation. Roughly 
speaking, they represent connections on higher dimensional 
analogues of circle bundles, n-gerbes. 
Let us start by looking at the low-dimensional models.  
\par In a given open covering $ \m{U}= \left ( U_a \right )_{ a \in \m{A}}  $ of $X$, a 
non-flat differential zero-cochain appears as a triplet $ \sigma = ( T, f_a, t_{ab} ) $ with $T$ 
being a global 1-form, $f_a$ local functions and $ t_{ab} $ $ 2 \pi \Zee $-valued assignments. 
If $ \sigma $ is 
a cocycle, relations $ dT =0 $, $ T_a = df_a $ and $ f_a - f_b = t_{ab} $
assure us that the local functions $ q_a = {\rm exp}(i f_a) $ glue together to a global
function $ q \colon X \rightarrow S^1 $. Moreover, the global 1-form $T$ is obtained by just
pulling back through $q$ the Maurer-Cartan form of $ S^1 $. Two zero-cocycles
are equivalent if and only if they determine the same global map $ q \colon X \rightarrow S^1 $. 
\par In a similar open neighborhood, a non-flat one-cochain appears as:
$$ \omega \ = \ \left ( H, \theta_a, h_{ab} , r_{abc} \right ) $$        
where $ H $ is a global 2-form, $ \theta_a $ are local 1-forms, $ h_{ab} $ are local 
functions and $ r_{abc} $ are $ 2 \pi \Zee $-valued integral assignments. Again, if $ \omega $ 
is a cocycle then the local ingredients above are related by $ dH=0 $, $  H_a = d \theta_a $, 
$ \theta_a - \theta_b = h_{ab} $, $ h_{bc} - h_{ac} + h_{ab} = r_{abc} $ and 
$ r_{bcd}- r_{acd} + r_{abd} - r_{abc} =0 $. Let $ g_{ab} \colon U_a \cap U_b \rightarrow S^1 $
be the local functions $ g_{ab} = {\rm exp}(i h_{ab} ) $. They satisfy the relation 
$ g_{ab} g_{ac}^{-1} g_{bc} =1 $ and therefore make a Cech 1-cocycle with 
a cohomology class in $ H^1(X, S^1) $. Such an object is known to determine a circle
bundle $ L \rightarrow X $. The rest of the cocycle conditions tell that local 1-forms 
$ \theta_a $ glue together to form a circle connection on $L$. The global two-form $ H $ 
is naturally the curvature of the connection. One can read the connection holonomy
along 1-loops directly from the cocycle $ \omega $.  Let us assume 
$ \gamma \colon S^1 \rightarrow X $ is a smooth loop
in X. The open covering $ (U_a)_{a \in \mathcal{A}} $ will, by restriction, 
cover $ \gamma $.
We choose a triangulation of the loop $ \gamma $ subordinated
to the covering. Such a feature is determined by a union: 
$$ \gamma = \Delta_1 \cup \Delta_2 \cup .... \cup \Delta_n $$ 
with $ \Delta_i \colon [0,1] \rightarrow \gamma $ , 
$ \Delta_i(1)=\Delta_{i+1}(0) $ for $ 1 \leq i \leq n-1 $ and 
$ \Delta_n(1)=\Delta_1(0) $. 
The triangulation $ \gamma $ is subordinated to the covering 
$ (U_a)_{a\in \m{A}} $ in the sense that a subordination map 
$ \rho \colon \{1, \cdots n \} \rightarrow \mathcal{A} $ is chosen such that 
$ \Delta_i \subset U_{\rho(i)} $.  
The holonomy of the connection around the loop $ \gamma $ can then be
computed as the $ {\rm exp}( i \cdot {\rm hol}_{\gamma}) $ where:  
\begin{equation}
\text{hol}_ {\gamma}= \ \sum_{i} \int_{\Delta_i}\theta_{\rho (i)} \ 
- \sum_{i} \int_{\Delta_i (0)}{h_{\rho (i) \rho (i-1) }}  
\end{equation} 
Here, $ \rho (-1) = \rho (n) $. So far, it seems that the above formula
for $ {\rm hol}_{\gamma} $ depends on the choice of triangulation and 
subordination map. However, one can show that under variations of 
triangulation of the loop or subordination map $ \rho $, the expression 
$ \text{hol}_{\gamma} $ gets modified by an element in $ 2 \pi \Zee $. 
The exponential $ \text{exp}(i \cdot \text{hol}_{\gamma}) $ remains
therefore unchanged under such modifications. 
\par Geometrically, a non-flat differential
one-cocycle is therefore just a circle bundle connection. One verifies immediately
that two one-cocycles are equivalent if and only if the corresponding connections
are gauge equivalent. 
\par Assume now that under the above conditions one introduces a non-flat one-cochain $ \sigma $ 
satisfying $ \check{d} \sigma = \omega $. Using the earlier notation, one obtains: 
$ H=dT$, $ \theta_a = T_a - df_a $, $ h_{ab} = -f_a + f_b + t_{ab} $ and 
$ r_{abc} = t_{bc} - t_{ac} + t{ab} $. From here, one deduces that $ g_{ab} = q_a^{-1} q_{b} $ 
and therefore the local functions $ q_a $ glue together to form a global section 
trivializing the circle bundle $ L $. Moreover, the holonomy of the connection is
described in this trivialization by just integrals of the global 1-form $ T $. One can
say therefore that relation $ \check{d} \sigma = \omega $ realizes $ \sigma $ as 
a ``geometrical'' section in $ L $, in the sense that the information provides a section 
together with its behavior under holonomy of the existing $S^1$-connection. 
\par The above arguments can be generalized to fit higher dimensional differential cochains. 
Let  $ \omega \in N \check{Z}^n(X) $ be a differential non-flat n-cocycle described
in an open covering $ (U_a)_{a \in \m{A}} $ as a multi-plet:
\begin{equation}
\omega = ( H, \ \omega^n_{a}, \ \omega^{n-1}_{a_1 a_2}, \ 
\omega^{n-2}_{a_1 a_2 a_3}, \ .... \ \omega^{0}_{a_1 a_2 ....a_{n+1}}, \ 
\omega^{-1}_{a_1 a_2 .... a_{n+2}} ).
\end{equation} 
The upper index represents the degree of the corresponding local form. The $-1$ index corresponds
to locally constant functions taking values in $ 2 \pi \Zee $. One defines local $S^1$-valued 
functions:
$$ g_{a_1 a_2 ....a_{n+1}} \colon U_{a_1} \cap U_{a_2} \cap \cdots \cap U_{a_{n+1}} \rightarrow 
S^1 , \ \ g_{a_1 a_2 ....a_{n+1}} = {\rm exp} ( i \cdot \omega^{0}_{a_1 a_2 ....a_{n+1}} ). $$
It is a straightforward computation that $ \left ( \delta g \right )_{a_1 a_2 .... a_{n+2}} = 1 $
and therefore $ g $ defines a Cech n-cocycle. This data defines (see \cite{hit} for details) a 
(n-1)-gerbe. These are objects which, due to their geometrical features, can be considered
higher dimensional analogues of circle bundles. The cocycle $ \omega $ defines then a connection
on such a (n-1)-gerbe. Such a connection carries, as we shall explain shortly, a holonomy
along any embedded closed n-manifold. Two connections are said to be gauge equivalent if
they determine similar holonomies. Moreover an equality of type $ \check{d} \sigma = \omega $ 
with $ \sigma \in \in N \check{C}^{n-1}(X) $ can also be explained geometrically. Similarly 
to the circle bundle case, the (n-1)-cochain $ \sigma $ determines a section trivializing
the n-gerbe underlying $ \omega $.    
\par Let us explain how the holonomy associated to a gerbe connection
of type $ \omega $ is defined. Suppose $ Y $ is a n-dimensional closed 
submanifold embedded in X. In order to define the holonomy 
of $ \omega $ along Y one needs two additional ingredients: 
\begin{enumerate}
\item $ \bf{A \ dual \ cell \ decomposition \ for \ Y} $. That is a 
decomposition dual
to a triangulation. (We need each vertex to be adjacent to $n$ edges,
each edge to be adjacent to $n-1$ 2-cells, and so on.) Let us describe 
the top cells as $ ( \Delta_i )_{i \in I} $. They inherit orientation
from the orientation of Y. We denote by $ \Delta_{i_1 i_2 ...i_k} $ 
the n+1-k dimensional cell obtained by intersecting $ \Delta_{i_1}, 
\ \Delta_{i_2}, \ ....\ \Delta_{i_k} $ (if such a cell exist). We
assume the following orientation convention. A cell 
$ \Delta_{i_1 i_2 ...i_k} $ receives orientation as boundary component
in $ \Delta_{i_1 i_2 ...i_{k-1}} $. That means $ \Delta_{i_1 i_2} $ 
is oriented as boundary component in  $ \Delta_{i_1} $ which is 
oriented a priori. $ \Delta_{i_1 i_2 i_3} $ gets orientation as part
of the boundary of $ \Delta_{i_1 i_2} $, and so on. We will refer to 
a multi-index in the form $ (i) = (i_1 i_2 \cdots i_k) $.
$ \Delta_{(i)} $ will then be a ($n+1-k$)-cell with a certain orientation.
Permuting the elements in the index does not change the cell but
the orientation changes according to the signature of the permutation. 
\item $ \bf{A \ subordination \ map} $, $ \ \rho \colon I \rightarrow 
\mathcal{A} $ such that $ \Delta_i \subset U_{ \rho (i) } $. 
\end{enumerate}
These being settled, one defines the holonomy of $ \omega $ along Y
as $ {\rm exp} \left (i \cdot {\rm hol}_Y (\omega) \right  ) $ where: 
\begin{equation}
\label{2000}
\text{hol}_{Y} ( \omega ) = \ \sum_{k=1}^{n+1} \ \
(-1)^{k+1} 
\sum_{ (i) = ( i_1 > i_2 > \cdots >i_k) }   
\ \ \int_{\Delta_{(i)}} 
\omega^{n+1-k}_{\rho(i_1) \rho(i_2) ... \rho(i_k) }
\end{equation}
\noindent Here we just pick an order relation on I to make sure we do not use
the same cell twice during summation. It can be seen that above formula 
generalizes the holonomy description for 1-cocycles presented earlier. 
The holonomy does not depend on the choice of dual cell decomposition 
or subordination map. That is justified by:
\begin{cl}
\label{eq2011}
The expression $ \text{hol}_{Y} ( \omega ) $ above varies by an element
in $ 2 \pi \Zee $ under modifications of the dual cell decomposition 
$ \Delta_i $ or of the subordination map $ \rho $.   
\end{cl} 
\begin{proof}
 We follow two steps. Firstly, we show that varying the 
subordination map for the same cell decomposition changes expression
$ (\ref{2000}) $  by an element in $ 2 \pi \Zee $. Secondly, 
we'll see that refining a 
cell decomposition
under the same subordination map does not change $ (\ref{2000}) $. 
These will be 
enough to prove Claim $ \ref{eq2011} $ . 
\par To start with the first step, let us assume that the initial subordination
map of $ ( \Delta_{i} )_{i \in I } $, : 
$$ \rho \colon I  \rightarrow \mathcal{A} $$ is modified over a unique
cell $ \Delta_{i_o} $ such that: 
$$ 
\widetilde{\rho}(i) \ = \ 
\begin{cases}
\rho(i) \  & \text{if} \ i \neq i_o  \ ; \\
\widetilde{\rho}(i_o) \neq \rho(i_o) & \text{if} \ i=i_o \ ; \\
\end{cases} $$ 
\par Let $ \text{hol}_{Y} ( \omega ) $ be the holonomy defined using the
subordination map $ \rho $ and 
$ \widetilde{\text{hol}}_{Y} ( \omega ) $ be the holonomy defined
using $ \widetilde{\rho} $. We claim the following happens: 
$$ 
\text{hol}_{Y} ( \omega ) \ - \ \widetilde{\text{hol}}_{Y} ( \omega )
\ = \ 
\sum_{(i)=(i_o=i_1>i_2> \cdots i_{n+1})} \ \int_{\Delta_{(i)}} \ 
\omega^{-1}_{\rho(i_1) \widetilde{\rho}(i_1) \rho(i_2) \cdots \rho(i_{n+1})} 
$$ 
The quantity on the right is then a sum of numbers in $ 2 \pi \Zee $ 
(integrals are just
functions evaluated in points) and the difference of the two holonomies 
is therefore $ 2 \pi $ times an  integer. That completes the first step.
Let us prove the relation above. Assume the order relation for indices 
$ i \in I $ is chosen such that $ i_o $ is the largest index. We make 
the following notations: 
$$ 
P_k \ = \ \sum_{ (i) = ( i_o=i_1 > i_2 > \cdots >i_k) }   
\ \ \int_{\Delta_{(i)}} 
\left ( 
\omega^{n+1-k}_{\rho(i_1) \rho(i_2) ... \rho(i_k) } \ - \ 
\omega^{n+1-k}_{\widetilde{\rho}(i_1) \rho(i_2) ... \rho(i_k) } 
\ \right ) 
$$ 
and 
$$ 
S_k \ = \ 
\sum_{ (i) = ( i_o=i_1 > i_2 > \cdots >i_{k}) } \int_{\Delta_{(i)}}   
d \omega^{n-k}_
{\rho(i_1) \widetilde{\rho}(i_1) \rho(i_2) \cdots \rho(i_{k})}    
$$ 
Therefore 
$$ 
\text{hol}_{Y} ( \omega ) \ - \ \widetilde{\text{hol}}_{Y} ( \omega )
\ = \ 
\sum_{k=1}^{n+1} \ \
(-1)^{k+1} \ P_k  
$$
But $ P_1 \ = \ S_1 $ and in general for $ 2 < k $ : 
$$
S_{k-1} \ + \ (-1)^{k+1} \ P_{k} \ = \ S_{k} 
$$ 
This results by just applying Stokes' Theorem. Hence: 
$$ 
\text{hol}_{Y} ( \omega ) \ - \ \widetilde{\text{hol}}_{Y} ( \omega )
\ = \ 
\sum_{k=1}^{n+1} \ \
(-1)^{k+1} \ P_k \ = \ S_1 \ +  \sum_{k=2}^{n+1} \ ( S_k - S_{k-1} ) \ 
 = \ S_{n+1}  
$$ 
And since: 
$$ 
S_{n+1} \ = \ 
\sum_{ (i) = ( i_o=i_1 > i_2 > \cdots >i_{n+1}) } \int_{\Delta_{(i)}}   
\omega^{-1}_
{\rho(i_1) \widetilde{\rho}(i_1) \rho(i_2) \cdots \rho(i_{n+1})}    
$$ 
the needed identity is proved. 
\par We cover now the second step. Let us assume we have two different dual
cell decompositions $ ( \Delta_i )_{i \in I} $ and 
$ ( \Delta_j )_{j \in J} $ with the latter one being a refinement of the 
former. Say, there's a map:  
$$ 
\varphi \colon J \rightarrow I \ \text{such that} \ 
\Delta_j \subset \Delta_{\varphi(j)} \ .
$$ 
Moreover, both decompositions are subordinated to the covering $ U_a $ 
through maps: 
$$ 
\rho \colon I \rightarrow \mathcal{A} 
$$ 
$$
\rho \circ \varphi \colon J \rightarrow \mathcal{A} \ . 
$$ 
Assume the indices in both I and J are ordered such that the refinement
map $ \varphi $ is increasing. We have two expressions for holonomy,
depending on which decomposition we are using: 
$$ 
\ \sum_{k=1}^{n+1} \ \
(-1)^{k+1} 
\sum_{ (i) = ( i_1 > i_2 > \cdots >i_k) }   
\ \ \int_{\Delta_{(i)}} 
\omega^{n+1-k}_{\rho(i_1) \rho(i_2) ... \rho(i_k) }
$$  and
$$ 
\ \sum_{k=1}^{n+1} \ \
(-1)^{k+1} 
\sum_{ (j) = ( j_1 > j_2 > \cdots >j_k) }   
\ \ \int_{\Delta_{(j)}} 
\omega^{n+1-k}_{\rho(j_1) \rho(j_2) ... \rho(j_k) }
$$ 
They are the same though. That is (roughly speaking) because 
$ \omega^{n+1-k}_{\rho(j_1) \rho(j_2) ... \rho(j_k) } $ vanishes
as soon as two sub-indices are the same. This completes the
second step. 
\end{proof} 

\subsection{Definition of the B-Field.}
We are in position to give a definition for the B-field \cite{f1}. Recall the notation 
in $ (\ref{impdef})$. 
\begin{dfn}
Let space-time $X$  be a smooth manifold.
A B-field on X is a non-flat differential 2-cochain 
$ B \in \m{B}^2(X) $ defined over an open covering $ \m{U} $.  
\end{dfn} 
\noindent Let us recall that for non-flat differential n-cochains in: 
$$  
\m{B}^n(X) \ = \ \left \{ \omega \in N \check{C}^n(\m{U}) \ \vert \ \m{U} \ 
{\rm covering}, \ 
\check{d}\omega \in \Omega^{n+1}(X) \right \} $$
one has the equivalence relation $\sim $ defined in $ \ref{aa} $ which extends
the standard cocycle equivalence. 
\begin{dfn}
Two B-fields $ B_1, B_2 \in \m{B}^2(X) $ are said to be (gauge) equivalent
if $ B_1 \sim B_2 $.
\end{dfn}
\noindent Two B-fields cannot be normally added up unless they are defined on the
same open covering. However, there is a well-defined summation rule 
on the set of equivalence classes:
\begin{equation}
\label{gr}
\m{B}^2 (X) / \sim.
\end{equation}  
For $ B_1, B_2 \in \m{B}^2(X) $, one defines $ [B_1 ] + [B_2] = [B'_1 + 
 B'_2] $, where $ B'_1, \ B'_2 $ are restrictions of $ B_1, \ B_2 $ on an open
covering refining both coverings underlying $ B_1 $ and $ B_2 $. This
definition does not depend on the choice of refinement or subordination maps.
In this respect, $ (\ref{gr}) $ becomes a group.
\par Let $ \omega \in \Omega^3(X) $ be a fixed 3-form and: 
$$ \m{T}_{\omega} = 
\{ [B] \ | \ B \in \m{B}^2(X), \ \check{d}B=\omega \}.$$
The equivalence relation does not modify the field strength $ H_B $ 
of a B-field $ B \in \m{B}^2(X) $. There is then an exact sequence:
$$ 0 \rightarrow \check{H}^2(X) \hookrightarrow \ \m{B}^2(X) / \sim \  
\ \stackrel{\Phi}{\rightarrow} \ \Omega^3(X) \times \Omega^3(X) $$ 
with $ \Phi(B) = (H_B, \ \check{d}B) $. Clearly, 
$  \m{T}_{\omega} = \Phi^{-1} ( \Omega^3(X) \times \{ \omega \} ) $.
One can therefore conclude that the field strength projection map:
$$ \m{T}_{\omega} \rightarrow \Omega^3(X) , \ \ [B] \to H_B $$
realizes a fibration over the image, the fibers being principal homogeneous
spaces for the second smooth Deligne cohomology, $ \check{H}^2(X) $. 
\par Let us review the geometrical meaning attached to these objects. 
Any B-field contains a 3-form  $ \omega = \check{d}B $ which (as a differential 
3-cocycle in the sense of $ (\ref{formm}) $) can be seen as a 
connection $ \mathcal{A}_{\omega} $ on the trivial 2-gerbe over X. 
$ \mathcal{A}_{\omega} $ has exact field strength  
(given by $ d H_B = d \omega $) and it may very well carry 
holonomy along closed 3-manifolds embedded in X. The B-field can be seen 
as a (geometrical) section
trivializing the 2-gerbe. Since B is not required to be flat, the section
is not necessarily parallel with respect to $ \mathcal{A}_{\omega} $. In
fact, its covariant derivative with respect to $ \m{A}_{\omega} $ is $ H_B $. 
\par In general, according to above definition, B-fields are just non-flat 
differential 2-cochains. They do not carry holonomy in the standard sense.
(Non-flat 2-cocycles can be viewed as connections on gerbes and do carry
holonomy along closed surfaces.) However, as mentioned above, 
$ B $ can be regarded then as a section trivializing a 2-gerbe.
But a non-flat 2-cocycle can be seen (integrated down, see the appendix for details) 
to give a circle bundle connection over the space of embedded closed surfaces 
in X. Its underlying circle bundle is trivial but has no natural 
trivialization.
A B-field which is a section of a 2-gerbe integrates down to produce a 
section of  this circle bundle. Thus, a B-field assigns to any closed 
surface $ \Sigma $ mapping to $X$, a point in the circle fiber over 
the point of the space of mappings given by $\Sigma \rightarrow X$. 
We denote this by:
$$ {\rm exp} \left ( i \int_{\Sigma} B \right ) $$ 
and interpret it as the holonomy of the B-field along $ \Sigma $. If B is
a cocycle ($ \check{d}B =0 $) the circle bundle has a canonical 
trivialization and with respect to that, the quantity above can be seen as 
a unitary complex number, well-defined up to an overall phase factor
independent of the mapping. It realizes the standard gerbe connection holonomy.
\subsection{The Chern-Simons Cocycle.}
\label{csss}
Let $G$ be a compact, simply-connected, simple Lie group and 
$P \rightarrow X$ a principal $G$-bundle. To any choice of connection 
$A$ on $P \rightarrow X$, one can associate a non-flat differential 
Chern-Simons 3-cocycle $ \m{CS}_A \in N \check{Z}^3(X) $. As shown 
by Freed in \cite{f1}, a general definition for $ \m{CS}_A $ can be achieved 
by pulling back a universal choice on $ BG $ through classifying maps. 
However, in this way $ \m{CS}_A $ is well-defined just up to a flat 
differential 3-coboundary. For computational reasons we would like to 
introduce $ \m{CS}_A $ directly as a 3-cocycle. Although not canonically,
this can be achieved by fixing a local trivialization of the bundle. Our
task is further eased by the fact that all $G$-bundles we shall be dealing with 
throughout this paper are trivializable. 
\par Let $ \mathbf{g} $ be the Lie algebra of $G$ and $ {\rm ad} \colon
G \rightarrow {\rm End}({\bf g}) $ the adjoint representation. 
Chern-Weil theory provides an isomorphism:
$$ H^4(BG, \mathbb{R}) \ \simeq \ I^2(G) $$ 
where $ I^2(G) $ represents the family of ad-invariant quadratic 
forms on $ {\bf g} $. Such a quadratic form is called integral if it 
represents an integral cohomology class in $ H^4(BG, \mathbb{R}) $. For
compact, simply connected, simple, Lie groups $ G $, 
$ H^4(BG, \Zee) \simeq \Zee$. The
integral ad-invariant forms make a free group a rank one. A generator is 
given by:
\begin{equation}
\label{innerprod}
q(a) \ = \ \frac{1}{16 \pi^2 c_G} \ <a, a>_k.
\end{equation}
Here, the right-hand side bracket denotes the Killing form on $ {\bf g} $. 
The number $ c_G $ is the Dynkin index of $ G $ which is always
integer or half-integer . $c_G = 1/2 $ for $ G = {\rm SU}(n) $.  For
$ G=E_8 $ the number $ c_G $ equals the dual Coxeter number $30$.
\par The Lie bracket on $ \bf{g} $ can be naturally extended to a graded 
Lie algebra bracket on  $ \Omega ^*(X, \ \bf{g} ) $. Precisely,  
$$ 
\left [\omega^p, \ \omega^q \right ] (v_1, v_2 \cdots v_{p+q}) = \ 
$$
$$  
= \ \frac{1}{(p+q)!} 
\sum_{\sigma \in S_{p+q}} (-1)^{{\rm sgn} (\sigma)} \left [ \ \omega^p
(v_{\sigma(1)} , 
v_{\sigma(2)} \cdots v_{\sigma(p)} ) , \ \omega^q (v_{\sigma(p+1)} , \cdots 
v_{\sigma(p+q)} ) \ \right ] \ .
$$
It satisfies  
$$
\left [ \omega^p, \ \omega^q \right ] = (-1)^{pq+1} 
\left [\omega^q, \ \omega^p \right ]
$$
as well as the Jacobi identity
$$
\left [ \omega^p , \ \left [ \omega^q, \ \omega^r  \right ] \right ] = 
\left [ \left [ \omega^q, \ \omega^r  \right ] , \ \omega^p \right ]
+ (-1)^{pq} 
\left [ \omega^q , \ \left [ \omega^p, \ \omega^r  \right ] \right ] \ .
$$
The inner product $ (\ref{innerprod}) $ extends as well to a pairing:
$$ \left < \cdot , \cdot \right > \colon 
\Omega^p(X, {\bf g} ) \otimes \Omega^q(X, {\bf g} ) \rightarrow 
\Omega^{p+q}(X) $$ 
$$
\left < \omega^p, \ \omega^q \right > (v_1, v_2 \cdots v_{p+q}) = 
$$
$$ 
= \ \frac{1}{(p+q)!} 
\sum_{\sigma \in S_{p+q}} (-1)^{{\rm sgn} (\sigma)} \left < \ \omega^p
(v_{\sigma(1)} , 
v_{\sigma(2)} \cdots v_{\sigma(p)} ) , \ \omega^q (v_{\sigma(p+1)} , \cdots 
v_{\sigma(p+q)} ) \ \right >_k .
$$ 
One can check that:
$$ \left < \omega^p, \ \omega^q \right > = (-1)^{pq} 
\left < \omega^q, \ \omega^p \right > \ \ \text{and} \  
\left < \omega^p, \ \left [ \omega^q, \ \omega^r  \right ] \right > = 
\left < \left [ \omega^p, \ \omega^q  \right ], \ \omega^r \right > .
$$  
On $G$, one has the $ {\bf g} $-valued Maurer-Cartan 1-form $ \theta_{G} $, 
which assigns to each vector its left-invariant extension. This satisfies:
$$ L_{g}^* \mc = \mc , \ \ {\rm and} \ \ R_g^* \mc = ad_{g^{-1}} \mc  $$
where $ L_{g},R_g \colon G \rightarrow G $ represent left and right 
multiplication with $ g \in G $. Its differential verifies the Maurer-
Cartan equation:
$$ d \mc + \frac{1}{2} \left [ \mc , \ \mc \right ] = 0. $$
Combining the two pairing operators we obtain a new 3-form on G:
$$ \m{W}_{G} = - \frac{1}{6} \ 
\left < \mc, \ \left [ \mc, \ \mc \right ] \right > \ \in 
\Omega^3(G, \mathbb{R} ). $$
It is a closed bi-invariant 3-form with integral periods. 
\par We introduce then the connection data. 
A connection on $ P \rightarrow X $ can be seen as a 1-form 
$ A \in \Omega^1(P, \ {\bf g}) $ satisfying:
$$ L_{p}^* A = \mc , \ \ {\rm and} \ \ R_g^* A = ad_{g^{-1}} A . $$
The curvature, 
$$ F = dA + \frac{1}{2} \ \left [ A, \ A \right ] $$ 
verifies:
$$ L_{p}^* F = 0 , \ \ {\rm and} \ \ R_g^* F = ad_{g^{-1}} F $$
as well as the Bianchi identity;
$$ dF + \frac{1}{2} \ \left [ A, \ F \right ] = 0 . $$
The Chern-Simons 3-form is by definition a global form on the total space:
$$ CS_A = \left < A, \ dA + \frac{1}{3} \left [ A, \ A \right ] \right > = 
\left < A, \ F - \frac{1}{6} \left [ A, \ A \right ] \right > \ 
\in \Omega^3(P, \mathbb{R} ) . $$
Its basic properties are:
\begin{itemize}
\item $ L_p ^* CS_A = \m{W}_G $
\item $ R_g ^* CS_A = CS_A $ 
\item $ d CS_A = \left < F, \ F \right > $  
\end{itemize}
the last relation being the standard computation for the Chern-Simons 
action.
\par So far, all objects were defined on the total space of the bundle. 
However, by employing a choice of local trivializations, all data can be 
transfered down to X. Here we simplify our discussion. We assume that
the $G$-bundle $ P \rightarrow X $ is trivial. We fix, once for all, 
a global section $ s \colon X \rightarrow P $ trivializing the bundle.
In this setting, everything can be pulled-back on $X$. We shall keep
the same notation for the curvature 1-form and its curvature and 
Chern-Simons three-form although what we really mean is their pull-back
through the section $s$. The Chern-Simons non-flat differential 3-cocycle 
can then be defined as the 3-cocycle induced by the global 3-form:
$$ 2 \pi \ 
\left < A, \ dA + \frac{1}{3} \left [ A, \ A \right ] \right > = 
2 \pi \ \left < A, \ F - \frac{1}{6} \left [ A, \ A \right ] \right > \ 
\in \Omega^3(X, \mathbb{R} ) . $$
In differential cocycle language the Chern-Simons 3-cocycle
can be represented in a random open covering as:
\begin{equation}
\label{cstr}
\m{CS}_A \ = \ 2 \pi \ \left ( \left < F, \ F \right > , \ 
\left < A, dA + \frac{1}{3} \left [ A, A \right ] \right >
 \vert _{U_a} , \ 0, \ 0, \ 0, \ 0 \right ).  
\end{equation}   
This makes a well-defined cocycle since the global 4-form 
$ 2 \pi  H = 2 \pi  \left < F, \ F \right >\in \Omega^4(X) $ has 
periods in $ 2 \pi \Zee $. If $ G={\rm SU}(n) $, $ H = {\rm ch}_2(F) $ and 
its associated 4-cohomology class recovers the second Chern class 
$ c_2(P) \in H^4(X, \Zee) $. Since the bundle is trivial this cohomology
class vanishes.
\par The formulation $ (\ref{cstr}) $ of $ \m{CS}_A $ depends on the choice 
of trivialization $s$. However, a variation of $s$ changes $ \m{CS}_A $ 
by a flat 3-coboundary. The gauge class of the 3-cocycle does not
change. The above definition of the Chern-Simons 3-cocycle can be
extended to non-trivial $G$-bundles. The construction process involves
fixing a family of local trivializations. However, we shall not deal with the
non-trivial bundle case here. 
\par We continue with: 
\begin{lem}
The holonomy of $ \m{CS}_A $ along a closed 3-manifold $W$ embedded in $X$ 
recovers the Chern-Simons invariant $ {\rm cs}(W, A) $.  
\end{lem} 
\begin{proof}
This follows from standard arguments in Chern-Simons theory. Let 
$ W \hookrightarrow X $ be an embedded closed 3-manifold. There always
exists a 4-manifold $ M $ such that $ \partial M = W $. The bounding 
4-manifold $ M $ is not necessarily embedded in $X$. Due to its triviality,
the bundle $ P $ extends to a $G$-bundle $ \w{P} $ over $M$. The global
section $S$ and connection
$A$ extend also to $ \w{s} $ and $ \w{A} $ on $ \w{P} $. Therefore 
the restriction on $W$ of the Chern-Simons cocycle $ \m{CS} $ extends
to a non-flat 3-cocycle $ \m{CS}_{\w{A}} $ on $ M $. In this setting:
\begin{equation}
\label{form}
{\rm hol}_{W}(\m{CS}_A) \ = \ {\rm exp} \left ( 2 \pi \cdot \int_M 
\left < \w{F} , \ \w{F} \right > \right ).
\end{equation}   
But $ (\ref{form}) $ is exactly the Chern-Simons invariant $ {\rm cs}
(W, A) $.   
\end{proof}
\vspace{.1in}
\noindent In what follows we analyze the way Chern-Simons cocycles
(seen here as global 3-forms due to triviality of the bundle) change
under the action of symmetry group $ \m{G} $ of bundle automorphisms 
$ \varphi $ for $ P $ covering orientation preserving diffeomorphisms 
$ \bar{\varphi} $ on X. 
\vspace{.1in}
\begin{teo} \
\label{thetath}
\begin{enumerate}
\item For each $ \varphi \in \m{G} $ and connection $A$ there exist 
a differential flat 2-chain $ \theta_{(A, \varphi)} \in \m{B}^2(X) $ such 
that
\begin{equation}
\label{eqimp}
\m{CS}_{\varphi^*A} \ = \ \bar{\varphi}^* \m{CS}_A + 
\check{d} \theta_{(A, \varphi)} .
\end{equation}
\item For any two $ \varphi_1 , \varphi_2 \in \m{G} $, the quantity:
\begin{equation}
\label{comb} 
\left [ \theta_{(A, \varphi_1 \circ \varphi_2)} \right ] - 
\left [ \theta_{(\varphi_1 ^*A, \varphi_2)} \right ]  - 
\left [ \bar{\varphi}_2^*\theta_{(A, \varphi_1)} \right ]  
\end{equation}
vanishes in $ \m{B}^2(X) / \sim $. 
\end{enumerate}  
\end{teo} 
\begin{proof}
We look at the first part of the theorem. The symmetry group 
$ \m{G} $ consist of all automorphisms $ \varphi $ 
of the bundle $P $ covering orientation preserving diffeomorphisms 
$ \bar{\varphi} $ on $X$. It includes as a normal subgroup the group
of standard gauge transformations $ \m{G}(P) $. Technically, we have 
the following short exact sequence of groups:
\begin{equation}
\label{eqseq}
\{ 1 \} \rightarrow \m{G}(P) \rightarrow \m{G} \rightarrow 
{\rm Diff}^+(X) \rightarrow \{ 1\}. 
\end{equation}
Due to triviality of the bundle the above sequence splits. Indeed there
is a map $ {\rm Diff}^+(X) \rightarrow \m{G} $, $ \bar{\varphi} 
\leadsto \bar{\varphi}_o $ sending a diffeomorphism $ \bar{\varphi} $
of $X$ to the unique automorphism of $P$ leaving the section $s$ 
invariant. This map builds a section for the second projection in the
exact sequence $ (\ref{eqseq}) $. Therefore any automorphism 
$ \varphi \in \m{G} $ can be decomposed uniquely as:
\begin{equation}
\label{deco}
\varphi \ = \ \psi \circ \bar{\varphi}_o  
\end{equation}
with $ \psi \in \m{G}(P) $. The symmetry group can then be understood as
a semi-direct product:
$$ \m{G}(P) \rtimes {\rm Diff}^+(X). $$ 
\par Let $ \psi $ be a standard gauge transformation. In the chosen 
trivialization $s$, $ \psi $ can be seen through a smooth function
$ t \colon X \rightarrow G $. For any connection $ A $ we get:
$$ \psi ^* A = {\rm ad} _{t^{-1}} A + t^* \theta_G. $$
The Chern-Simons three-form gets modified as follows:
\begin{equation}
\label{rell}
\left < \psi^*A, d\psi^*A + \frac{1}{3} 
\left [ \psi^*A, \psi^*A \right ] \right > 
\ - \ \left < A, dA + \frac{1}{3} \left [ A, A \right ] \right >
\ = \ d \left < {\rm ad}_{t^{-1}} A , \ t^* \theta_G \right > +
t^* \m{W}_G.   
\end{equation}
The 3-form $ \m{W}_G $ is closed and has integral periods. 
Therefore, one could interpret $ 2 \pi \cdot \m{W}_G $ as a flat 3-cocycle 
with trivial holonomy. Hence, we can pick a universal choice of a flat 
2-cochain $ \eta \in \check{C}^2(G) $ 
satisfying $ \check{d} \eta = 2 \pi \cdot \m{W}_G $. Based
on this assumption we define:
\begin{equation}
\label{th}
\theta_{(A, \psi)} \ = \ 
2 \pi \cdot \left < {\rm ad}_{t^{-1}} A , \ t^* \theta_G \right > \ + \ 
t^* \eta 
\end{equation}
The first term on the right-hand side of $ (\ref{th}) $ is to be interpreted 
as a non-flat differential 2-cochain representable in a random covering
$ U_a $ as:
$$ 2 \pi \cdot \left ( 0 , \ 
\left < {\rm ad}_{t^{-1}} A , \ t^* \theta_G \right > \vert _{U_a} , \ 
0, \ 0, \ 0 \ \right ). $$ 
In this setting $ \theta_{(A, \psi)} $ as defined in $ (\ref{th}) $ 
becomes a differential flat 2-cochain. A straight-forward computation
based on $ (\ref{rell}) $ shows that:
\begin{equation}
\label{step1}
\m{CS}_{\psi^*A} \ = \ \m{CS}_A + 
\check{d} \theta_{(A, \psi)} .
\end{equation}
Here we make use of the fact that $ \m{CS}_A $ and $ \m{CS}_{\psi^*A} $, 
as 3-cocycles associated to global 3-forms, can be represented in any choice 
of open covering, in particular on the covering underlying 
$ \theta_{(A, \psi)} $. The above expression is then $ ( \ref{eqimp}) $ in 
the case of a standard gauge transformation.
\par Let us consider $ \varphi \in \m{G} $. According to $ (\ref{deco}) $,
we can uniquely decompose $ \varphi $ as:
$$ \varphi \ = \psi \circ \bar{\varphi}_o \ $$ 
where $ \bar{\varphi}_o $ is the unique automorphism covering 
$ \bar{\varphi} $ on base space, 
preserving the trivializing section $s$ and  $ \psi \in \m{G}(P) $. 
Clearly:
\begin{equation}
\label{relll}
\m{CS}_{\varphi^*A} \ = \ \m{CS}_{(\psi \circ \bar{\varphi}_o)^* A} \ 
= \ \m{CS}_{\bar{\varphi}_o ^*(\psi ^* A)} \ = 
\ \bar{\varphi}^*\m{CS}_{\psi^*A} \ = \  
\bar{\varphi}^* \left ( \ \m{CS}_A + \check{d} \theta_{(A , \psi)} \ 
\right ).  
\end{equation}  
Defining:
$$ \theta_{(A, \varphi)} \ = \ \bar{\varphi}^*\theta_{( A , \psi)} $$
and using this in $ (\ref{relll}) $ we obtain:
$$ \m{CS}_{\psi^*A} \ = \ \bar{\varphi}^* \m{CS}_A + 
\check{d} \theta_{(A, \psi)} $$
which is exactly equation $ (\ref{eqimp}) $. 
\par We prove now the second part of the theorem. As before, we start by 
analyzing the case when the two automorphisms involved are just gauge 
transformations. Let $ \psi_1, \psi_2 \in \m{G}(P) $. We have:
$$
\begin{array}{ccccc}
P & \stackrel{\psi_1}{\rightarrow} & P & \stackrel{\psi_2}{\rightarrow} &
P \\
\Big\downarrow & & \Big\downarrow & & \Big\downarrow \\
X & = & X & = & X .\\
\end{array} 
$$
The three flat 2-cochains involved, $ \theta_{(A, \psi_1 \circ \psi_2)} $, 
$ \theta_{(\psi_1 ^*A, \psi_2)} $ and $ \theta_{(A, \psi_1)} $ live on
different open coverings. However, our goal is to prove an equality involving
their equivalence classes. We take a common refinement and make choices of 
subordination maps restricting therefore the three flat 2-cochains on a unique 
common open covering where they can be added up. The equivalence classes
won't be affected by the choice of refinement. In this setting:   
\begin{equation}
\label{bbb}
\theta_{(A, \psi_1 \circ \psi_2)} - \theta_{(\psi_1 ^*A, 
\psi_2)} - \theta_{(A, \psi_1)} \ = 
\end{equation}
$$ =  
2 \pi \cdot 
\left < {\rm ad}_{(t_2t_1)t^{-1}} A , (t_2t_1)^* \theta_G \right > - 
2 \pi \cdot 
\left < {\rm ad}_{t_2^{-1}} ({\rm ad}_{t_1^{-1}}A + t_1^* \theta_G ) , 
t_2^* \theta_G \right > - 
$$
$$
- \ 2 \pi \cdot \left < {\rm ad}_{t_1^{-1}} A, t_1^* \theta_G \right > + 
(t_2t_1)^* \eta - t_2^*\eta - t_1^*\eta = $$
\begin{equation}
\label{bbbb}
= 2 \pi \cdot 
\left < {\rm ad}_{t_2^{-1}} t_1^*\theta_G, \ t_2^* \theta_G \right > + 
(t_2t_1)^* \eta - t_2^*\eta - t_1^*\eta.  
\end{equation}
Therefore, expression $ (\ref{bbb}) $ does not depend on connection A. 
Now, clearly, $ (\ref{bbbb}) $ makes a differential 2-cocycle. That is
because:
$$ \check{d} \left ( 
2 \pi \cdot 
\left < {\rm ad}_{t_2^{-1}} t_1^*\theta_G, \ t_2^* \theta_G \right > + 
(t_2t_1)^* \eta - t_2^*\eta - t_1^*\eta \right ) \ = \ 
$$   
$$
= 2 \pi \cdot 
d \left < {\rm ad}_{t_2^{-1}} t_1^*\theta_G, \ t_2^* \theta_G \right > +
(t_2t_1)^* \m{W}_G - t_2^* \m{W}_G - t_1^* \m{W}_G \ = 0. 
$$
It is also a flat 2-cocycle since all 2-cochains $ \theta_{(A, \psi)} $ 
are flat. We show $ (\ref{bbbb}) $ is actually a coboundary of a flat
differential 1-cochain. That follows from the fact that its holonomy
along any embedded compact oriented 2-manifold vanishes. Let $ \Sigma \subset 
X $ be an embedded smooth surface. There always exist a compact, 
oriented 3-manifold $W$ such that $ \partial W = \Sigma $. It is 
not necessarily that
$W$ is embedded in $ X$. The bundle $ P $ extends to $ \w{P} $ 
over $W$ and so does the trivialization $s$. Obstruction theory 
(based on $ \pi_1(G)=0 $) shows that the two gauge transformations 
$ \psi_i , \ i \in \{1, 2 \} $ extend to gauge transformations $ 
\w{\psi}_i , \ i \in \{ 1 , 2 \} $ in $ \w{P} $. They correspond
to functions: $ \w{t}_i \colon W \rightarrow G $. Accordingly,
the flat 2-cocycle $ (\ref{bbbb}) $ extends to a flat 2-cocycle
\begin{equation}
\label{bbbbb}
2 \pi \cdot \left < {\rm ad}_{\w{t}_2^{-1}} \w{t}_1^*\theta_G, \ 
\w{t}_2^* \theta_G \right > + 
(\w{t}_2\w{t}_1)^* \eta - \w{t}_2^*\eta - \w{t}_1^*\eta.
\end{equation}
But, in such situations the holonomy of $ (\ref{bbbb}) $ along
$ \Sigma $ is just the exponential of the strength field of 
$ (\ref{bbbbb}) $ integrated over $ W $. The latter vanishes since
the strength field of $ (\ref{bbbbb}) $ is zero. 
\par We look now at the general case. Assume that $ \varphi_1, 
\varphi_2 \in \m{G} $ covering $ \bar{\varphi}_1, \bar{\varphi}_1 \in 
{\rm Diff}^+(X) $  
$$
\begin{array}{ccccc}
P & \stackrel{\varphi_1}{\rightarrow} & P & \stackrel{\varphi_2}{\rightarrow} &
P \\
\Big\downarrow & & \Big\downarrow & & \Big\downarrow \\
X & \stackrel{\bar{\varphi}_1}{\rightarrow} & X & 
\stackrel{\bar{\varphi}_1}{\rightarrow} & X \\
\end{array}  
$$
and they decompose as:
$$ \varphi_i \ = \ \psi_i \circ \bar{\varphi}_{io} , \ \ i \in \{1,2\}. $$
By definition 
$$ \theta_{(A, \varphi_i)} \ =  \bar{\varphi}_i^*\theta_{( A , \psi_i)}. $$
Let 
$ \psi'_1 = \bar{\varphi}_{2o} \circ \psi_1 \circ 
\bar{\varphi}_{2o}^{-1} \in \m{G}(P) $. For any connection
$A$ we then have:
\begin{equation}
\label{eqbb}
\theta_{(\bar{\varphi}_{2o}^*A, \psi_1)} \ = \ \bar{\varphi}_2 ^* 
\theta_{(A, \psi'_1)}.
\end{equation}  
Therefore, expression $ (\ref{comb}) $ becomes:
\begin{equation}
\label{ee}
\left [ \theta_{(A, \varphi_2 \circ \varphi_1)} \right ] - 
\left [ \theta_{(\varphi_2 ^*A, \varphi_1)} \right ] - 
\left [ \bar{\varphi}_1^*\theta_{(A, \varphi_2)} \right ]  \ = \ 
(\bar{\varphi}_{2o} \circ \bar{\varphi}_{1o})^* 
\left [ \theta_{(A, \psi_2 \circ \psi'_1)}\right ]  - 
\end{equation} 
$$  
- \ \bar{\varphi}_1^* \left [ \theta_{(\varphi_2^*A, \psi_1)} \right ] - 
\bar{\varphi}_1^* \bar{\varphi}_2^* 
\left [ \theta_{(A, \psi_2)} \right ]  . 
$$
Now, from the case of pure gauge transformations we know that, on some 
common refining sub-covering:
\begin{equation}
\label{eee}
\theta_{(A, \psi_2 \circ \psi'_1)} - \theta_{(\psi_2 ^*A, 
\psi'_1)} - \theta_{(A, \psi_2)}
\end{equation}
is a coboundary of a differential flat 1-cochain. Putting together
relations $ (\ref{ee}) $ and $ (\ref{eee}) $ we obtain that:
$$ 
\left [ \theta_{(A, \varphi_2 \circ \varphi_1)} \right ] - 
\left [ \theta_{(\varphi_2 ^*A, \varphi_1)} \right ] - 
\bar{\varphi}_1^* \left [ \theta_{(A, \varphi_2)} \right ] \ = \ 
\bar{\varphi}_1^* \bar{\varphi}_2^* 
\left [ \theta_{(\psi_2 ^*A, \psi'_1)} \right ] \ - 
$$
$$ - \ 
\bar{\varphi}_1^* \left [ \theta_{(\varphi_2^*A, \psi_1)} \right ] + 
\left [ \check{d} \{ {\rm flat \ 1-cochain} \} \right ] .
$$ 
Using equality  $ (\ref{eqbb}) $ with connection $ \psi_2^*A $ we notice
that the first term in the right-hand side of above expression vanishes.
We can therefore conclude that:
$$ \left [ \theta_{(A, \varphi_2 \circ \varphi_1)} \right ] - 
\left [ \theta_{(\varphi_2 ^*A, \varphi_1)} \right ] - 
\bar{\varphi}_1^* \left [ \theta_{(A, \varphi_2)} \right ] \ = 0 .$$  
\end{proof}
\noindent The above features are the basic facts about Chern-Simons cocycles 
when $G$ is a simply-connected, simple, Lie group. However, the two compact groups
we shall use throughout the paper, $ \eeight $ and $ \spintt $, do not 
quite fall in this category. In what follows we adapt the previous discussion to those cases.
\newline  
\newline 
$ {\bf G = \eeight. } $ A trivializable $G$-bundle $ P \rightarrow X$ 
can be understood as a sum $ P_1 \times P_2 $ with $ P_i \rightarrow X, \ i \in 
\{ 1, 2 \} $ trivializable $ E_8 $-bundles. A section $ s $ is $ P $ 
can be fixed in the form of a pair $ (s_1, s_2) $ with $ s_i $ 
section in $ P_i $. Following the pattern, one can view a $G$-connection on 
$ P $ as a pair $ A = ( A_1, A_2) $ with $ A_i $, $i=1,2$ as $E_8 $-connections on
$ P_i $'s. We define then the Chern-Simons cocycle as:
\begin{equation}
\m{CS}_A^{G} \ = \ \m{CS}_{A_1}^{E_8} + \m{CS}_{A_2}^{E_8}.
\end{equation}        
The differential 2-cochains $ \theta_{(A, \varphi)} $ can be reconstructed 
as well. The group of automorphisms $ \varphi $ of $ P $ has as index-two subgroup 
the group of automorphisms $ \varphi = 
( \varphi_1 , \varphi_2) $ with $ \varphi_i $ automorphisms of $ P_i $. 
The $ \mathbb{Z}_2$ quotient is generated by the automorphism 
that exchanges the two $E_8$ factors. We go ahead and define for $ \varphi = 
( \varphi_1 , \varphi_2) $:
\begin{equation}
\theta^G_{(A, \varphi)} \ = \ \theta^{E_8}_{(A_1, \varphi_1)} +
 \theta^{E_8}_{(A_2, \varphi_2)}.  
\end{equation}
A straightforward computation shows that Theorem $ \ref{thetath} $ is 
still true. 
\newline
\newline
$ {\bf G = Spin(32) / \Zee_2. } $ Let $ P \rightarrow X $ be a trivializable
$ G = {\rm Spin}(32) / \Zee_2 $-connection. We set a section s. The Lie group 
projection  
$ {\rm Spin}(32) \rightarrow {\rm Spin}(32) / \Zee_2 $ induces an
isomorphism at Lie algebra level. We define the Chern-Simons cocycle
corresponding to a connection $ A $ by using the trace pairing
of the standard 32-dimensional representation of $ {\rm Spin}(32) $. 
To view it in a different way, say $ \w{P} $ is the trivial lifting
of $ P $ to a $ {\rm Spin}(32) $-bundle
\begin{equation}
\begin{array}{ccc}
\w{P} & \stackrel{p}{\rightarrow} & P \\
\Big\downarrow & & \Big\downarrow \\
X & = & X \\
\end{array}
\end{equation}
together with a section $ \w{s} = p^* s$. 
There is then a $1$-to-$1$ correspondence (given by pull-back through p) 
between $ G $-connections
on $ P $ and $ {\rm Spin}(32) $-connections on $ \w{P} $. The Chern-Simons
cocycle can then be defined as:
\begin{equation}
\m{CS}^G_A \ = \ \frac{1}{2} \cdot \m{CS}^{{\rm Spin}(32)}_{p^*A}.
\end{equation}   
The $\theta_{(A, \varphi)} $ 2-cochains can be extended similarly. Here
we make use of the fact that, as we mentioned before, the symmetry group 
$ \m{G} $ in the $ {\rm Spin}(32) / \Zee_2$ case consists of only 
liftable automorphisms. So, any $ \varphi \in \m{G} $ covering 
$ \bar{\varphi} \in {\rm Diff}^+(X) $ can be seen as coming from 
an automorphism $ \w{\varphi} $ of $ \w{P} $. We define then:
\begin{equation}
\label{33}
\theta^G_{(A, \varphi)} \ = \ \frac{1}{2} \cdot \theta_{(p^*A, \w{\varphi})}.
\end{equation}   
It can be seen that, up to a coboundary of a differential flat 1-cochain, 
 the right side of $ (\ref{33}) $ does not depend
on the choice of automorphism lifting $ \w{\varphi} $. Indeed, say
$ \w{\w{\varphi}} $ is a different lifting. Then $ \w{\w{\varphi}} =
 \xi \circ \w{\varphi} $ where $ \xi $ is the gauge transformation in
$ \w{P} $ which consist in each fiber of right multiplication with
$ q $, the exotic element in $ \Zee_2 $. Let $ \w{\varphi} = 
\psi \circ \bar{\varphi}_o $ with $ \psi $ gauge transformation. Then,
for any connection $ \w{A} $ on $ \w{P} $, 
$$ \theta_{(\w{A}, \w{\varphi}) } \ = \ \bar{\varphi}^* 
\theta_{(\bar{\varphi}^*\w{A}, \psi) } \ \ {\rm and} \ \ 
\theta_{(\w{A}, \w{\w{\varphi}}) } \ = \ \bar{\varphi}^* 
\theta_{(\bar{\varphi}^*\w{A}, \xi \circ \psi)}. 
$$   
But, up to a coboundary of a flat 1-chain, 
$$ \theta_{(\bar{\varphi}^*\w{A}, \xi \circ \psi)} =  
\theta_{(\bar{\varphi}^*\w{A}, \psi) } + 
\theta_{(\psi^*\bar{\varphi}^*\w{A}, \xi) }. $$
The last term on the right side above vanishes. That happens because,
denoting $ B = \psi^*\bar{\varphi}^*\w{A} $ and applying definition   
$ (\ref{th}) $ one can write:
\begin{equation}
\label{mmm}
\theta_{(B, \xi) } = 2 \pi \cdot \left < {\rm ad}_{t^{-1}} B, t^* \theta_G 
\right > + t^* \eta 
\end{equation}
where $ t \colon X \rightarrow {\rm Spin}(32) $ is the map associated
to the gauge transformation in the given trivialization. In the case of 
$ \xi $, the map $ t $ is constant. Therefore, both terms in the right
side of $ (\ref{mmm}) $ vanish. Going back to the two equations above, 
that is enough to conclude:     
$$ \theta_{(\w{A}, \w{\varphi}) } \ = \ \theta_{(\w{A}, \w{\w{\varphi}}) }
$$
up to a coboundary of a flat 1-cochain. 
\par Expression $ (\ref{33}) $ gives then a well-defined formulation for
$ \theta^G_{(A, \varphi)} $. Theorem $ \ref{thetath} $ is verified. \\ \\
We finish this section with a short remark about the gravitational 
$ SO(10)  $ Chern-Simons.
\newline
\newline
$ {\bf G = SO(10). } $ The gravitational Chern-Simons term 
$ \m{CS}_g $ appears in the anomaly cancellation condition $ (\ref{anccond}) $ 
described in next chapter. It represents 
the Chern-Simons cocycle associated to the Levi-Civita connection induced
by riemannian metric $g$ on $ TX $. That can be seen as a connection on
the $ {\rm SO}(10) $-principal bundle of orthonormal oriented frames 
associated to $TX$. So $ \m{CS}_g $ is basically the Chern-Simons 
cocycle of a ${\rm SO}(10) $-connection. However, the space-time $X$ 
is endowed from the beginning with a spin structure. The Levi-Civita
connection lifts uniquely to a $ {\rm Spin}(10) $ connection on the
fixed spin bundle. We define $ \m{CS}_g $ to be half of 
the Chern-Simons associated to the $ {\rm Spin}(10) $-lifted Levi-Civita. 
The Lie group $ {\rm Spin}(10) $ is simply-connected (and in our particular
case of space-time the spin structure is trivializable) so the entire
previous discussion applies.  
%
%
%
%
\section{The String Data} 
\label{sec3}
In this section, we analyze the parameter data for heterotic string theory. In general, a 
$10$-dimensional spin manifold $X$ is chosen as a space-time and 
a gauge bundle is fixed to be a $G$-principal bundle $ P \rightarrow X $. 
For anomaly cancellation reasons \cite{wgs}, the Lie 
group $G$ must be one of the two choices $ \eeight $ or 
$ \spintt $. The string data specifying a theory 
consist of a triplet
$(g, \ A, \ B ) $ involving a metric on $X$, a $G$-connection on $P$ 
and a B-field in the form of a differential non-flat 2-cochain on $X$. 
Equations of Motion \cite{wgs} impose the following conditions:
\begin{enumerate}
\item The metric $g$ on X is Ricci flat. 
\item The $G$-connection $A$ on $P$ has vanishing curvature.
\item The B-field B has vanishing field strength.   
\end{enumerate}  
These classical solutions to the Equations of Motion are to be 
considered up to gauge equivalence. However, before factoring to equivalence 
classes, there is one more constraint to be taken into account. This is the 
anomaly cancellation condition interpreted by Freed \cite{f1} as:
\begin{equation}
\label{anccond}
\check{d}B = \m{CS}_A -  \m{CS}_g .
\end{equation}
$ \m{CS}_A $ and $ \m{CS}_g  $ are the Chern-Simons 3-cocycles 
corresponding to the connection $A$, respectively the lifted 
Levi-Civita connection
associated to the metric $g$ on the fixed spin bundle of $X$. 
Their construction was explained in section $ \ref{csss} $. 
The possibility of imposing such a 
constraint requires an a priori topological condition on space-time $X$. 
Namely,
\begin{equation}
\label{topocond}
\lambda(P) = {\rm p}_1(TX) \ {\rm in} \ H^4(X, \ \mathbb{Z} ). 
\end{equation} 
Here $ \lambda \in H^4(BG, \Zee) $ represents the level of the theory. 
\par In case $ G = \eeight $ one obtains 
$ H^4(BG, \Zee) = H^4(BE_8 , \Zee) \times H^4(BE_8, \Zee) $. It is a known
fact \cite{mt} that for a simply-connected, simple Lie group $H$, 
$ H^4(BH, \Zee) \simeq \Zee $. A generator $ \xi $ for 
$ H^4(BE_8, \Zee) $ is obtained via Chern-Weil theory from the 
ad-invariant quadratic form:
$$ I_o \colon e_8 \rightarrow \mathbb{R} , \ \ 
\ \ I_o ( a) = \frac{1}{16 \pi^2 c_{E_8}} \  
\left < a, a \right >_k $$    
where $e_8 $ is the Lie algebra of $E_8 $, 
$ \left < \cdot , \cdot \right > $ represents the Killing form on 
$ e_8 $ and $ c_{E_8} = 30 $ is the dual Coxeter number. The level
of the $ G = \eeight $ theory is then chosen as 
$ \lambda = \left ( \xi , \xi \right ) \in H^4(BG, \Zee) $.     
\par The case $ G = {\rm Spin}(32) / \Zee_2 $ can be treated similarly.
One obtains $ H^4(BG, \Zee) \simeq \Zee $ and the level of the theory
is chosen to be the generator $ \lambda \in H^4(BG, \Zee) $ representing
in Chern-Weil theory the quadratic form: 
$$ \frac{1}{2} I_o \colon {\rm spin}(32) \rightarrow \mathbb{R} , \ \ 
\ \ \frac{1}{2} I_o ( a) = \frac{1}{2} \cdot 
\frac{1}{16 \pi^2 c_{{\rm Spin}(32)}} \  
\left < a, a \right >_k. $$    
\par The particular case of eight-dimensional heterotic string is simpler from this 
point of view. One compactifies from ten to eight dimensions. The space time is 
$ X = E \times \mathbb{R}^8 $, with $E$ a 2-dimensional torus. Condition
$ (\ref{topocond}) $ is satisfied by default since there is no $4$-cohomology.
In this particular framework, the classical solutions to the 
Equations of Motion are invariant under orientation preserving 
isometries of $ \mathbb{R}^8 $ and we consider only 
$ \mathbb{R}^8 $-invariant solutions. Therefore we can equivalently regard the three 
objects involved as:
\begin{enumerate}
\item g, flat metric on $E$.
\item $ A $, G-connection on a fixed principal bundle $ P \rightarrow E $. 
\item $B$, flat differential cochain on $E$ with 
$ \check{d}B = \m{CS}_A - \m{CS}_g $. 
\end{enumerate}  
If $G= \eeight $ all $G$-bundles over a 2-torus are topologically 
trivial. However, for $ G= {\rm Spin}(32)/ \mathbb{Z}_2 $ there exist a 
non-trivial topological type. The topological type of the bundle is 
characterized by the generalized 
Stiefel-Whitney class $ \w{w}_2 \in H^2(X, \mathbb{Z}_2) $ which measures 
the obstruction to lifting the structure group of the bundle to 
${\rm Spin}(32) $. Nevertheless, the ${\rm Spin}(32)/\mathbb{Z}_2 $ heterotic string 
assumes that the gauge bundle allows such a lifting (P carries a vector structure in 
physics terminology). Therefore we consider from now on that 
the $G$-bundle $ P \rightarrow E$ is topologically trivial.   
\par Let us denote:  
$$
\w{\m{M}}_{{\rm het}}= 
\left \{ (g, \ A, \ [B]) \ \vert \ 
g \ {\rm flat \ metric} , \ F_A = 0, \ H_B = 0, \ 
[B] \in \m{T}_{\omega} \ {\rm where} \ \omega = 
\m{CS}_A -  \m{CS}_g \right \}
$$ 
and 
$$ \m{P}_G = \left \{ (g, \ A) \ \vert \ 
g \ {\rm flat \ metric} , \ F_A = 0 \right \}. 
$$ 
One cannot consider the equivalence
class of a general differential 2-cochain. 
However, all bundles involved here are 
trivializable. Therefore, as explained in section 
$ \ref{csss} $, the two Chern-Simons
cocycles $ \m{CS}_A $ and $ \m{CS}_g $ are global 3-forms. The B-fields
appearing in the definition of $ \w{\m{M}}_{{\rm het}} $ are then included 
in $ \m{B}^2(X) $ and their
equivalence class is well-defined. 
\par The above $\w{\m{M}}_{{\rm het}}$ and $\m{P}_G$ spaces 
carry natural smooth structures. There is a projection map:
\begin{equation}
\label{proj1}
\pi \colon \w{\m{M}}_{{\rm het}} \rightarrow \m{P}_G, 
\end{equation}
realizing a smooth fibration. 
The space of gauge classes of flat 2-cocycles, $ \check{H}^2(E) \simeq 
S^1 $ acts freely and transitively on each fiber.
In this sense, $ (\ref{proj1}) $ becomes a principal circle bundle. This
raw picture must be divided out by symmetries. The symmetry group $ \m{G} $ 
consists of bundle automorphisms covering any orientation preserving 
diffeomorphism of the base:
$$
\begin{displaystyle}
\begin{array}{ccc}
P & \stackrel{\varphi}{\rightarrow} & P \\
\downarrow & & \downarrow \\ 
E & \stackrel{\bar{\varphi}}{\rightarrow} & E \\
\end{array}
\end{displaystyle} 
$$
\begin{equation}
\m{G} \stackrel{{\rm def}}{=} \left \{ \varphi \in {\rm Aut}_G(P) \ \vert 
\ \bar{\varphi} \in {\rm Diff}^+(E) \right \}.  
\end{equation}     
(Strictly speaking this is the symmetry group for $ \eeight $ case.
If $ G= {\rm Spin}(32)/ \Zee_2 $, $ \m{G} $ is only made out of those 
automorphisms that can be lifted to $ {\rm Spin}(32) $ automorphisms.)
The group $ \m{G} $ is a non-trivial extension of the 
diffeomorphism group of $E$ by 
the group of changes of gauge of the bundle $P$: 
$$ \{ 1 \} \rightarrow \m{G}(P) \rightarrow \m{G} \rightarrow  {\rm Diff}^+(E)
\rightarrow \{1\}. $$ 
$\m{G} $ acts naturally on both, the total and base spaces of 
$ (\ref{proj1}) $. 
The action on $ \m{P}_G $ is the standard one: 
$$ \m{P}_G \times \m{G} \rightarrow \m{P}_G , \ \ (g, \ A). \varphi = 
( \bar{\varphi} ^* g, \ \varphi ^* A). $$
In order to establish the action on the total space, we use the following
facts: 
$$ 
\m{CS}_{\bar{\varphi} ^* g} = \bar{\varphi}^* \m{CS}_g
$$ 
and 
$$
\m{CS}_{\varphi ^* A } = \bar{\varphi}^* \m{CS}_A + \check{d} 
\theta_{(A, \ \varphi)}. $$
Here $ \theta_{(A, \ \varphi)} $ is a flat 2-cochain in $ \m{B}^2(X)  $
as described in section $ \ref{csss} $. The action $ \w{\m{M}}_{{\rm het}} \times \m{G} 
\rightarrow \w{\m{M}}_{{\rm het}} $ can then be described as:
$$ (g, \ A, \ [B]). \varphi = 
( \bar{\varphi} ^* g, \ \varphi ^* A, \ 
[\bar{\varphi} ^* B ] + [ \theta_{(A, \ \varphi)}] ). $$
It is a well defined action since, as is established in $ \ref{csss} $, 
$$ 
[ \theta_{(A, \ \varphi_1 \circ \varphi_2)} ] -  [ 
\theta_{(\varphi_1 ^* A, \ \varphi_2)} ] - 
[ \bar{\varphi}_2^*\theta_{(A, \ \varphi_1)} ]  = 0.    
$$ 
The symmetry group action commutes with the projection $ \pi $ factoring
the circle bundle $ (\ref{proj1}) $. However, as we shall see, there are
points in the base for which the stabilizer group does not act trivially
on the fiber above. Nevertheless, all stabilizer groups are finite. 
The circle bundle $ (\ref{proj1}) $ descends to a circle fibration:
\begin{equation}
\label{proj2}
\w{\m{M}}_{{\rm het}}  \ / \m{G}  \rightarrow  \m{P}_G \ / \m{G}. 
\end{equation}
The total space above represents exactly the moduli space of heterotic 
string parameters $ \m{M}_{{\rm het}} $. Our goal is to describe this space 
in detail and characterize the circle fibration $ (\ref{proj2}) $.
\vspace{.1in}
%
%
%
\section{Anomaly Cancellation: Chern-Simons vs. Pfaffians}
\label{sec4}
The discussion in section $ \ref{sec3} $ was based on a smooth family of metrics and 
$G$-connections on a 2-torus. The associated space of B-field equivalence 
classes forms 
the total space of a principal circle bundle over the parameter space. 
Moreover, this bundle carries a natural connection. Now, let $ \rho $ be a 
complex unitary representation of the Lie group G. 
Using $ \rho $ one can construct
for each pair $ (A, g) $ a complex elliptic operator, the coupled Dirac 
operator. A determinant line
bundle can be associated to this family of operators. It carries canonical
Quillen connection and metric. Under certain conditions, a tensor combination 
of such determinant bundles recovers completely the B-field circle bundle
with its connection.
The determinant tensor combination forms the so-called fermionic anomaly 
of the family.  
The identification existing between the two bundles is known as the 
Green-Schwartz anomaly cancellation. This section provides the details of the
construction.  
\subsection{The Chern-Simons Bundle}
To begin with, we establish the following general setting. Let 
$ p \colon Z \rightarrow Y $ be a smooth fibration of manifolds with all 
fibers isomorphic to a 2-torus $E$. We assume the following geometrical data:
\begin{itemize}
\item A spin structure and a metric $g^{(Z/Y)} $ on the tangent bundle along 
the fibers $ T(Z/Y) \rightarrow Z $.
\item A projection $TZ \rightarrow T(Z/Y) $.
\item A (trivializable) $ SO(n) $-bundle $ U \rightarrow Z $ endowed with
a connection $ A^{{\rm grav}} $.
\item A (trivializable) $G$-principal bundle $ Q \rightarrow Z $ endowed with 
a connection $ A^{{\rm gauge}} $.
\item A complex unitary representation 
 $ \rho \colon  G \rightarrow U(r) $ carrying
a real structure (meaning that representation $ \rho $ 
can be obtain by complexification from a real representation). 
\end{itemize}
The tangent bundle over the fibers of $p$ can be endowed with a canonical 
connection. Indeed, let $ g^{Y} $ be an arbitrary riemannian metric on $Y$.
Since $ TZ = p^* TY \oplus T(Z/Y) $, the pull-back $ p^* g^{Y} $ 
together with $g^{(Z/Y)} $ generates a riemannian metric on Z. Let
$ \nabla^{Z} $ be its Levi-Civita connection. Projecting $ \nabla^{Z} $ 
on the tangent bundle along fibers we obtain a connection 
$ \nabla^{(Z/Y)} $. It is independent of the choice of metric $ g^{Y} $ 
on $Y$ \cite{MR88h:58110a}. 
\par The Chern-Simons cocycle associated to the family:
\begin{equation}
\label{cs99}
\bold{CS} = \m{CS}_{A^{\rm gauge}} - \m{CS}_{A^{\rm grav}}. 
\end{equation} 
is a non-flat 3-cocycle\footnote{Due to the fact that
the two bundles involved are trivializable, the Chern-Simons cocycle 
$ \bold{CS} $ actually represents, as explained in section $ \ref{csss} $, 
a global 3-form.} in $N \check{Z}^3(Z)$. 
It can be seen as a
connection on a 2-gerbe over Z. Fixing an $ y \in Y $ we obtain a 2-torus
$E_y = p^{-1}(y) $ together with a G-connection 
$ A^{{\rm gauge}}_y = A^{\rm gauge} \vert _{E_y} $ 
and a metric $ g_y $ on $ TE_y $. The restriction
$ A^{{\rm grav}}_y = A^{\rm grav} \vert _{E_y} $ 
recovers the Levi-Civita connection associated to $g_y$. Therefore
$$ \bold{CS} \vert _{E_y} = \m{CS}_{A^{\rm gauge}_y} - 
\m{CS}_{A^{\rm grav}_y} \in N \check{Z}^3(E_y). $$
The Chern-Simons 3-cocycle $ (\ref{cs99}) $ can be pushed-forward 
(see section $\ref{appendix} $) to $Y$, defining a 1-cocycle:
\begin{equation}
\label{cs100}
\omega = \int_{Z/Y} \bold{CS} \in N \check{Z}^1(Y) 
\end{equation}     
which can be geometrically interpreted as a circle connection on a 
certain principal circle bundle
\begin{equation}
\label{a11}
\m{R}_Y \rightarrow Y. 
\end{equation}
We call this the Chern-Simons bundle associated to the family.
\par Now, we incorporate the B-fields. Let:
\begin{equation}
\label{bundle}
\w{M}_Y =  \left \{ (y, \ [B]) \ \vert \ 
[B]  \in \m{T}_{\omega} , \ {\rm where} \ \omega = 
\m{CS}_{A^{\rm gauge}_y} - 
\m{CS}_{A^{\rm grav}_y} \right \} .
\end{equation}
$ \w{M}_Y  $ can be endowed with a canonical topology and smooth structure.
The projection $ \pi \colon \w{M}_Y \rightarrow Y $ makes then a circle
bundle under the action on the total space of 
$ \check{H}^2(E) \simeq S^1 $.  
\vspace{.1in}
\begin{teo}
\label{thbfield}
There exist a canonical circle bundle isomorphism $ \Phi_Y $ :
$$
\begin{displaystyle}
\begin{array}{ccc}
\w{\m{M}}_Y & \stackrel{\Phi_Y}{\rightarrow} & \m{R}_Y \\
\Big\downarrow & & \Big\downarrow \\
Y & = & Y \\ 
\end{array} 
\end{displaystyle}
$$ 
identifying the B-field parameter space $ \w{\m{M}}_Y $ to the total space
of the circle bundle $ (\ref{a11}) $.   
\end{teo} 
\begin{proof} The correspondence $ \Phi_Y $ goes as follows. Let 
$ (y, \ [B] ) $ be an element in $ \w{\m{M}}_Y $. B is
a non-flat 2-cochain on $E_y$ satisfying:
\begin{equation}
\label{eq11}
\check{d}B =  \m{CS}_{A^{\rm gauge}_y} - 
\m{CS}_{A^{\rm grav}_y}.
\end{equation}
We take: 
\begin{equation}
\label{inttt}
\Phi_Y ([B]) \ = \ \left [ \int_{E_y} B \ \right ]
\end{equation}
The integral inside $ (\ref{inttt}) $ is explained in section $\ref{appendix} $. Although
the construction of the integral itself on 2-cochains in not canonical, its
(gauge) equivalence class is well-defined.   
$ \Phi_Y([B]) $ is simply an equivalence class of a 0-cochain over the point 
$ y \in Y $. (Roughly speaking such a 0-cochain over a point is just a real
number, two 0-cochains being equivalent if the difference of the two
numbers is in $ 2 \pi \Zee $.) But $  \m{CS}_{A^{\rm gauge}_y}  -   
\m{CS}_{A^{\rm grav}_y} = \bold{CS} \vert _{E_y} $ and equation 
$ (\ref{eq11}) $ shows
B can be interpreted as a section in the 2-gerbe underlying 
$ \m{CS}_{A^{\rm gauge}_y} - \m{CS}_{A^{\rm grav}_y}   $.
The integration mechanism commutes with differentiation 
(see section $ \ref{appendix} $). 
Hence, formally:
$$ \check{d} \left ( \int_{E_y} B \ \right ) \ = 
\ \int_{E_y} \left ( \m{CS}_{A^{\rm gauge}_y} - 
\m{CS}_{A^{\rm grav}_y} \right ).  $$ 
In this light,
$ \Phi_Y([B]) $ can be regarded as a geometrical section in the circle
bundle underlying the 1-cocycle
\begin{equation}
\int_{E_y} \left (  \m{CS}_{A^{\rm gauge}_y} - \m{CS}_{A^{\rm grav}_y} 
\right )  \ = \ 
\left ( \int_{Z/Y}  \bold{CS} \right ) \Big{\vert} _{y} 
\ = \ \omega  \vert _{y}  
\end{equation}
Yet, this circle bundle is just the circle $ \m{R}_{y} $ over 
$ y$.  Therefore
a section represents just a a point in $ \m{R}_{y} $. Hence,  
$ \Phi_Y([B]) $ can be seen as being an element inside the total space of the 
Chern-Simons bundle $ (\ref{a11}) $. 
\par The map $ \Phi_Y $ is well-defined since two gauge-equivalent B-fields
integrate to gauge-equivalent 0-cochains generating the same section-point
in $ \m{R}_Y $.  One checks immediately that it realizes an equivariant isomorphism.
\end{proof}  
\vspace{.1in}
\begin{rem}
\label{rem1}
As seen above, the total space of the Chern-Simons bundle $ (\ref{a11}) $ 
describes the space of equivalence classes of B-fields associated to a family 
of metrics and connections $Y$. The connection $ \omega $ naturally induced
on this circle bundle has curvature:
$$ \Omega_{\omega} = \int_{Z/Y} 
\left ( H_{A^{\rm gauge}} - H_{A^{\rm gauge}} \right ) $$
where  $ H_{A^{\rm gauge}} $ and $ H_{A^{\rm grav}} $ are the strength
fields of the Chern-Simons cocycles $ \m{CS}_{A^{\rm gauge}} $ and 
$ \m{CS}_{A^{\rm grav}} $. The holonomy along a smooth loop $ \gamma 
\colon S^1 \rightarrow Y $ can be obtained as:
$$ {\rm hol}_{\omega}(\gamma) = {\rm cs}(W, A_{{\rm gauge}}) \cdot    
{\rm cs}(W, A_{{\rm grav}})^{-1}, $$ 
the two factors on the right side representing Chern-Simons
invariants along the closed 3-manifold $W$ swept by $ \gamma $ inside
$Z$.
\end{rem} 
\vspace{.1in}
\subsection{Chern-Simons vs. Pfaffians.}  
Next, we establish a relation between the above Chern-Simons circle
bundle and determinant bundles for analytic Dirac operators. We follow
arguments described in \cite{f5}. 
There are fixed metrics and spin structures on 
each 2-torus $E_y$ in the family. Therefore, one can define Dirac 
operators: 
\begin{equation}
\label{nondirac}
\m{D}^+_y \colon S^+_y \rightarrow S^-_y. 
\end{equation}
There is a holomorphic interpretation for these operators. A metric on a 
surface induces a complex structure through its conformal class. 
The spin structure is then a holomorphic root $ K^{1/2} $ for 
the canonical line $ K $. Dirac
operators become just $ \overline{\partial} $-operators:
\begin{equation}
\label{part}
\overline{\partial}_{K^{1/2}} \colon \m{E}^0(K^{1/2}) \rightarrow
\m{E}^{0,1}(K^{1/2}).  
\end{equation} 
They build a determinant line bundle $ {\rm Det}(\m{D}^+_y) \rightarrow Y $
which comes endowed with a Quillen metric and connection \cite{qq}.
Moreover, the natural pairing:
\begin{equation}
\m{E}^0(K^{1/2}) \otimes \m{E}^{0,1}(K^{1/2}) \rightarrow \mathbb{C} , \ \ \ 
a \otimes b  \ \leadsto \int_E a \otimes b 
\end{equation}
gives an isomorphism:
$$ \m{E}^{0,1}(K^{1/2}) \simeq \m{E}^0(K^{1/2})^*. $$
The Dirac operators $ (\ref{part}) $ can then be regarded 
as skew-adjoint operators:
\begin{equation}
\label{part2}
\overline{\partial}_{K^{1/2}} \colon \m{E}^0(K^{1/2}) \rightarrow
\m{E}^{0}(K^{1/2})^*.  
\end{equation}   
In these conditions the corresponding determinant line bundle admits a square root. 
This square root is the pfaffian complex line bundle (see \cite{f5} for details):
\begin{equation}
\label{pf}
{\rm Pfaff}(\m{D}^+_y) \rightarrow Y. 
\end{equation} 
It comes equipped with a metric and connection which are half the ones
on the determinant line. 
\par The framework can be further developed. One couples the Dirac 
operators $ (\ref{nondirac}) $ to the following families of connections: 
\begin{itemize}
\item connection $ \nabla^{{\rm grav}}_y $ induced by $A^{{\rm grav}}_y$ 
in the rank n complex bundle associated to $U_y$ through the (complexified) 
standard representation of $SO(n)$.
\item connection $ \nabla^{{\rm gauge}}_y $ induced by $A^{{\rm gauge}}_y$ 
in the rank r complex hermitian bundle associated to $Q_y$ through the 
representation $ \rho $.  
\item connection $ \nabla^{Z/Y}_y $ on 
$ TE_y \otimes_{\mathbb{Z}} \mathbb{C} $.   
\end{itemize}
We obtain three families of elliptic operators: $ \m{D}^{{\rm grav}}_y$, 
$ \m{D}^{{\rm gauge}}_y$ and $\m{D}^{(Z/Y)}_y$. Each family builds 
a pfaffian line bundle. The pfaffians can be constructed since
all three complex bundles supporting the coupling connections have
natural real structures.  
\par We consider the following fermionic combination depending on two integers
$ \alpha $ and $ \beta $.:  
\begin{equation}
\label{ferman}
\m{L}_{{\rm ferm}} \ = \ 
{\rm Pfaff}( \m{D}^{{\rm gauge}}_y )  \otimes
\left ( {\rm Pfaff}( \m{D}_y) \ \right )^{\otimes \beta} 
\otimes \ 
\left ( {\rm Pfaff} ( \m{D}^{(Z/Y)}_y) \ \right )^{\otimes \alpha} 
\otimes \ \left ( {\rm Pfaff} ( \m{D}^{{\rm grav}}_y) \ \right )
^{\otimes (-\alpha)} .     
\end{equation}
The tensor product $ (\ref{ferman}) $ produces a complex line bundle 
$ \m{L}_{{\rm ferm}} \rightarrow Y $ and is
endowed with a canonical Quillen metric $ g_q$ and compatible unitary 
connection$ \nabla_q$ \cite{f5} \cite{MR88h:58110a}. Moreover, in the 
case of a holomorphic family $ Y $, $ (\ref{ferman}) $ inherits a holomorphic
structure compatible with the connection. Restricting $ (\ref{ferman}) $ to
unit circles in each fiber we obtain a circle bundle 
$ \w{\m{L}}_{{\rm ferm}} \rightarrow Y $, the phase pfaffian. 
The Quillen connection descends to $ \w{\m{L}}_{{\rm ferm}} $. 
\par Under certain conditions the phase pfaffian can be identified (up to
an overall phase) 
to the Chern-Simons circle bundle $ (\ref{a11}) $. Explicitly, let us 
assume that $G$ is a simply-connected, compact Lie group and 
$ \lambda \in H^4(BG, \ \Zee) $ is the integral class giving the level
of the Chern-Simons cocycle. The strength field of the Chern-Simons 
cocycle is then a de Rham representative for the image of $ \lambda $ 
in real cohomology. The unitary representation 
$ \rho \colon G \rightarrow U(r) $ induces a cohomology map at the level
of classifying spaces:
$$ B \rho^4 \colon H^4(BU(r), \ \mathbb{Z}) \rightarrow 
H^4(BG, \ \mathbb{Z}). $$ 
\vspace{.01in}
\begin{teo}
\label{ancancel}
If $ B \rho ^4 (c_2) = 2 \alpha \cdot \lambda $ and 
$ \alpha (n+ 22) = r + \beta $ 
then the Quillen connection $ \nabla_q $ on the phase pfaffian 
$ \w{\m{L}}_{{\rm ferm}} $ 
and the $\alpha$-multiple of the Chern-Simons connection $ \omega $ existing
on the $\alpha$-th tensor power of $ (\ref{a11}) $ have the same curvature 
and holonomy.
\end{teo}
\begin{proof} 
To begin with, we compute the curvature of the fermionic anomaly 
$ (\ref{ferman}) $. The basic ingredient here is Bismut-Freed formula 
\cite{f5} which computes the curvature of the Quillen connection on a general 
pfaffian line bundle. In our particular case we obtain:
\begin{equation}
\label{integg}
\Omega^{\m{L}_{{\rm ferm}}} \ = \ \pi  \ \left ( \int_{Z/Y} \ 
\hat{A}(\Omega^{(Z/Y)}) \left ( {\rm ch}(\Omega^{{\rm gauge}}) + \beta 
+ \alpha \ 
 {\rm ch}(\Omega^{(Z/Y)}) - \alpha \ {\rm ch}(\Omega^{{\rm grav}}) \right )
\right )_{(2)}.
\end{equation} 
Let: 
$$ p = p_1 \left ( \nabla^{(Z/Y)} \right ) $$
be the Chern-Weil representative for the first Pontriagin class written
in terms of connection $ \nabla^{(Z/Y)} $ on the vertical tangent
space $ T (Z/Y) \rightarrow Y $.
Then:
$$ \hat{A}(\Omega^{(Z/Y)}) \ = \ 1 - \frac{p}{24} + \cdots $$ 
$$ {\rm ch}(\Omega^{{\rm gauge}}) \ = \ r + 
{\rm ch}_1(\Omega^{{\rm gauge}}) + \ 
{\rm ch}_2(\Omega^{{\rm gauge}}) + \cdots $$
$$ {\rm ch}(\Omega^{{\rm grav}}) \ = \ n + 
{\rm ch}_1(\Omega^{{\rm grav}}) +  \ 
{\rm ch}_2(\Omega^{{\rm grav}}) + \cdots $$
$$ {\rm ch}(\Omega^{(Z/Y)}) \ = \ 2 + p + \cdots $$
Rewriting the relevant terms inside the integral $ (\ref{integg}) $, one
obtains:
$$
\Omega^{\m{L}_{{\rm ferm}}} \ = \ \pi  \ \int_{Z/Y} \ 
\left ( \frac{\alpha(22+n)-\beta -r}{24} \ p + 
{\rm ch}_2(\Omega^{{\rm gauge}}) 
- \alpha \ {\rm ch}_2(\Omega^{{\rm grav}}) \right ). $$
Assuming $ r + \beta =\alpha(22+n)  $, the above expression becomes:
\begin{equation}
\label{vvv}
\Omega^{\m{L}_{{\rm ferm}}} \ = \ \pi  \ \int_{Z/Y} \ 
\left ( {\rm ch}_2(\Omega^{{\rm gauge}}) 
- \alpha \ {\rm ch}_2(\Omega^{{\rm grav}}) \right ). 
\end{equation}          
The strength field of the gauge Chern-Simons is:
$$ H_{CS_{A^{{\rm gauge}}}} \ = \ 
2 \pi \cdot \frac{1}{2 \alpha}  \cdot {\rm ch}_2(\Omega^{{\rm grav}}) $$
whereas for the gravitational Chern-Simons:
$$ H_{CS_{A^{{\rm grav}}}} \ = \ 
2 \pi \cdot \frac{1}{2}  \cdot {\rm ch}_2(\Omega^{{\rm grav}}). $$
Combining this with $ (\ref{vvv}) $ we obtain 
$$ 
\Omega^{\m{L}_{{\rm ferm}}} =  \alpha \ \int_{Z/Y} \ \left ( 
H_{CS_{A^{{\rm gauge}}}} - H_{CS_{A^{{\rm grav}}}} \right ). $$
\par The equality of holonomies is more delicate. It relies on a key 
relation between $ \xi $-invariants and Chern-Simons invariants
\cite{aaa}. The facts we need here are best summarized in the following 
lemma (\cite{f2}, Proposition 3.20):
\begin{lem}
\label{atiyah}
Let $M$ be a spin three-manifold, $g$ a riemannian metric on $M$ and 
$ E \rightarrow M $ a complex hermitian vector bundle with compatible
connection $A$. We denote by $ \m{D}^+_{(g,A)} $ the standard Dirac 
operator coupled to connection $A$. Its $ \xi$-invariant is defined
as:
$$ \xi(A) \ = \ \frac{\eta(A)+h(A)}{2}$$
where $ \eta(A) $ is the spectral eta-invariant \cite{MR88h:58110b}
 of $\m{D}^+_{(g,A)} $ 
and $h(A) $ is the dimension of the kernel of this operator. 
\par Assume that for a certain positive integer $N$, the formal combination:
$$ \left [ N \cdot \hat{A}(\Omega^g) \ {\rm ch}(\Omega^A) \right ]
_{(4)} $$
represents the Chern-Weil representative of an integral characteristic
class $ c \in H^4(BU, \mathbb{Z} ) $. We denote by $ {\rm cs}(A) \in U(1) $
the level $c$ Chern-Simons invariant associated to the connection $A$. 
Then, the unitary 
complex number:
\begin{equation}
\label{xics}
e^{2 N \pi i  \cdot \xi(A) } \cdot {\rm cs}(A)^{-1}  
\end{equation}
is a spin bordism invariant depending only on the class of 
$ E \rightarrow M $ in $ \Omega_3^{{\rm spin}} (BU) $. In particular,
if $E$ extends over a spin 4-manifold bounding $M$, then $ (\ref{xics}) $ 
vanishes. 
\end{lem}
We use this technical lemma for our purposes. Let $ \gamma \colon 
S^1 \rightarrow Y $ be a smooth loop. We denote by $W $ 
the spin 3-manifold swept by $ \gamma $ inside $Z$. We pick 
a metric on $S^1$. It induces a metric on $g_o $ on $W$. For each 
$ \epsilon $ we construct the Dirac 
operator associated to the scaled metric $ g_o / \epsilon ^2 $. Coupling
this Dirac operator with the restrictions of each of the three connections
existing on $Z$, $ A^{{\rm gauge}} $, $ A^{{\rm grav}} $ and 
$ \nabla^{(Z/Y)} $ we obtain three elliptic operators. The
three $ \xi$-invariants associated to them are obtain on the following 
pattern \cite{MR88h:58110b}:
$$ \xi_{\epsilon} \ = \ \frac{\eta_{\epsilon}+ h_{\epsilon} }{2} $$
where $ \eta_{\epsilon} $ represents the eta-invariant associated to the
corresponding Dirac operator and $ h_{\epsilon} $ is the dimension of the 
kernel. In fact, as shown in \cite{f2},  in dimension 3 the 
$ \xi $-invariant is independent of the metric and therefore the 
index $ \epsilon $ can be ignored. The fermionic holonomy along $ \gamma $ 
can be obtained through Witten's adiabatic limit formula \cite{MR88h:58110b}
which for this particular case stands as:
\begin{equation}
{\rm hol}_{\m{L}_{{\rm ferm}}}(\gamma) \ = \ 
e^{- \pi i  \xi } 
\ \ {\rm where} \ \ \xi = 
\xi^{{\rm gauge}} + \alpha \ \xi^{(Z/Y)} - 
\alpha \ \xi^{{\rm grav}} .
\end{equation} 
Now, the curvature considerations explained earlier together with Lemma 
$ \ref{atiyah} $ led us to conclude that:
\begin{equation}
\label{combb}
e^{- \pi i  \xi } \cdot {\rm cs}(W, A_{{\rm gauge}})^{-\alpha} \cdot 
{\rm cs}(W, A_{{\rm grav}})^{\alpha} 
\end{equation}     
is a spin bordism invariant. Since 
any given closed spin 3-manifold can be realized as the boundary of a 
spin 4-manifold
and all connections considered here extend, $ (\ref{combb}) $ must
vanish. 
\end{proof}
\vspace{.1in} 
\noindent The theorem $\ref{ancancel}$ allows one to identify 
(up to an overall phase) the two circle 
bundles:
\begin{equation}
\begin{array}{ccc}
\m{R}_Y^{\otimes \alpha} & \simeq & 
\w{\m{L}}_{{\rm ferm}} \\
\Big\downarrow &  & \Big\downarrow \\
Y & = & Y .\\ \\
\end{array}
\end{equation}
Next, we are going to apply this pattern in two particular cases
relevant to our discussion. 
\subsection{The adjoint pfaffian.}
Let us return to heterotic variables on a 2-torus. 
For simplicity we analyze the case  $ G = \eeight $. Let us consider 
$ Y =\m{P}_G $ be the family 
of pairs of flat metrics on $E$ and flat G-connections on the trivial
G-bundle $ P \rightarrow E $. The group of symmetries $ \m{G} $ consists
of bundle automorphisms covering orientation preserving diffeomorphisms
on the base space $X$. As explained
earlier the moduli space of string variables $ \m{M}_{{\rm het}} $ makes
the total space of a fibration:
\begin{equation}
\label{proj5}
\w{\m{M}}_{{\rm het}} \ / \m{G}  \rightarrow  
\m{P}_G  \ / \m{G} 
\end{equation}
which is obtained by factoring out the action of the symmetry group $ \m{G} $ 
from a circle line bundle:
\begin{equation}
\label{proj6}
\pi \colon \w{\m{M}}_{{\rm het}} \rightarrow \m{P}_G. 
\end{equation}
Moreover, the action of $ \m{G} $ on the base space has finite stabilizer
groups. $ (\ref{proj5}) $ makes a circle fibration.
\par We apply the earlier discussion to a particular framework involving 
a family $Y$ given by the following geometrical data:
\begin{itemize}
\item $ Y = \m{P}_G $ , $ Z = E \times \m{P}_G $, fibration $ p \colon Z 
\rightarrow Y  $ is just projection on the first factor.  
\item the tangent bundle along the fibers, $ T(Z/Y) $ is just 
$ TE \times \m{P}_G \rightarrow E \times \m{P}_G $. It is endowed
with a metric $ g^{(Z/Y)} $ and a connection $ \nabla^{(Z/Y)} $.
\item $U \rightarrow E \times \m{P}_G $ is the rank 10 real vector 
bundle obtained by direct summing $ T(Z/Y) $ with the rank 8 trivial
vector bundle over $ E \times \m{P}_G $. It is endowed with metric
$ g^{{\rm grav}} = g^{(Z/Y)} \oplus g_{{\rm prod}} $ and connection 
$ A^{{\rm grav}} = \nabla^{(Z/Y)} \oplus \nabla_{{\rm prod}} $ ( 
$ g_{{\rm prod}} $ and $ \nabla_{{\rm prod}} $ are the product metric and 
connection on the rank 8 trivial real bundle over $ E \times \m{P}_G $).     
\item $ Q = P \times \m{P}_G $. We endow it with a gauge connection
$ A^{{\rm gauge}} $ such that $ A^{{\rm gauge}} \vert _{P \times (A, \ g) } 
= A $. 
 \item $ \rho $ is the adjoint representation of G. It has rank 496. \\
\end{itemize}
By theorem $ \ref{thbfield} $, the line 
bundle $ (\ref{proj6}) $ can
be canonically identified to the Chern-Simons circle bundle associated
to the family $ Y =\m{P}_G $:
\begin{equation}
\begin{array}{ccc}
\w{\m{M}}_{{\rm het}} & \stackrel{\simeq}{\rightarrow} & \m{R}_{\m{P}_G} \\
\Big\downarrow &  & \Big\downarrow \\
\m{P}_G & = & \m{P}_G .\\ \\
\end{array}  
\end{equation}
The right hand side bundle comes equipped with a Chern-Simons circle 
connection. We consider
the fermionic pfaffian phase combination $ \w{\m{L}}_{{\rm ad}} $ associated 
to this particular geometric data with $ \alpha = c_{E_8} =30 $ and 
$ \beta = 464 $. Interpreting theorem $ \ref{ancancel} $ 
accordingly we obtain:
\begin{coro}
There exist a circle bundle isomorphism (unique up to multiplication by  
a unitary complex number) between the order $ c_{E_8} $ tensor power of  
circle bundle  $ (\ref{proj6}) $ and the adjoint phase pfaffian 
$ \w{\m{L}}_{{\rm ad}} $,
\begin{equation}
\label{izo}
\begin{array}{ccc}
\w{\m{M}}_{{\rm het}}^{\otimes c_{E_8}} & \simeq  & 
\w{\m{L}}_{{\rm ad}}  \\
\Big\downarrow &  & \Big\downarrow \\
\m{P}_G & = & \m{P}_G .\\ \\ 
\end{array}  
\end{equation} 
identifying the Chern-Simons connection to a $c_G$ order fraction
of the Quillen connection.   
\end{coro}
\begin{proof}   
For a simply-connected, simple, compact Lie group G, $ H^4(BG, \Zee) 
\simeq \Zee $. An explicit generator $ \theta $ corresponds to the integral
ad-invariant quadratic form:
$$ q \colon {\bf g} \rightarrow \mathbb{R} , 
\ \ \ \ q(a) \ = \ \frac{1}{16 \pi ^2 c} \ <a, a>_k $$ 
where the above right-hand side pairing $ < \cdot, \cdot > _k $ represents
the Killing form. The number $ c$ is always integer or half-integer. 
For $ G = {SU}(n) $, $ c= 1/2$. For $G = E_8 $, $ c $ equals the dual 
Coxeter number $ c_{E_8} = 30 $. 
\par Here we analyze the case $ G = \eeight $. The Chern-Simons
level $ \lambda \in H^4(BG, \Zee) \simeq \Zee \oplus \Zee $ corresponds
to $ (\theta, \theta) $ where $ \theta $ is the generator in $ H^4(
BE_8 , \Zee ) $. The adjoint representation $ {\rm ad} \colon 
E_8 \rightarrow SU(248) $ induces a cohomology map:
$$ B{\rm ad}^4 \colon H^4(BSU(248), \Zee) \rightarrow H^4(BE_8, \Zee) 
, \ \ \ {\rm with} \ \  B{\rm ad}^4(c_2) = 2c_{E_8} \cdot  \theta = 
60 \cdot \theta.$$ 
The adjoint representation of $ G $ is then just $ \rho = ({\rm ad},
{\rm ad}) \colon E_8 \times E_8 \rightarrow SU(248) \times SU(248) 
\hookrightarrow SU(496) $ and induces:
$$ B\rho^4 \colon H^4(BSU(496), \Zee) \rightarrow H^4(BSU(248), \Zee)
\times H^4(BSU(248), \Zee) 
\stackrel{B{\rm ad}^4 \times B{\rm ad}^4}{\rightarrow} H^4(BG, \Zee) 
$$
with $ B\rho^4(c_2) \ = \ (2 c_{E_8} \cdot \theta, 
 2 c_{E_8} \cdot \theta) \ = \ 2 c_{E_8} \cdot \lambda $.
Taking then $ \alpha = c_{E_8} = 30 $ and $ \beta = 464 $, both conditions
in Theorem $ \ref{ancancel} $ are satisfied.   
\end{proof}
\vspace{.1in}
In the light of above statement, the heterotic fibration 
$ \w{\m{M}}_{{\rm het}} \rightarrow \m{P}_G $ can be regarded naturally 
as a root of order $c_{E_8} $ for the adjoint pfaffian phase 
$ \w{\m{L}}_{{\rm ad}} $. The latter can be further simplified. 
One can write:
$$ \m{D}^{{\rm grav}}_{(A, \ g)} =  \m{D}^{(Z/Y)}_{(A, \ g)} \oplus 
\m{D}_g(\mathbb{C}^8). $$
The second operator on the right hand side represents just the Dirac
operator on $E$ associated to the metric g coupled with the product
connection on the trivial rank 8 vector bundle over $ E $. This does
not depend on connection $ A$ and is just eight times the standard
Dirac operator $ \m{D}_g $ associated to metric g. Therefore:
$$ {\rm det}(\m{D}^{{\rm grav}}_{(A, \ g)}) =  
{\rm det}(\m{D}^{(Z/Y)}_{(A, \ g)}) \otimes 
{\rm det}(\m{D}_g)^{8}. $$
Rewriting the fermionic pfaffian combination $ (\ref{ferman}) $,  
one obtains:
\begin{equation}
\label{newpfaff}
\m{L}_{{\rm ad}} = \  
{\rm Pfaff}( \m{D}^{{\rm gauge}}_{(A, \ g)} ) \
\otimes \ {\rm Pfaff} 
( \m{D}_g) ^{ \otimes (\beta - 8 \alpha)} \ = \ 
{\rm Pfaff}( \m{D}^{{\rm gauge}}_{(A, \ g)} ) \
\otimes \ {\rm Pfaff} 
( \m{D}_g) ^{ \otimes 224}      
\end{equation}  
The $ \w{\m{L}}_{{\rm ad}} $ represents just the phase of
$ (\ref{newpfaff}) $. 
\par The action of symmetry group $ \m{G} $ on $ \m{P}_G $ preserves
the geometric data. Therefore,  the group $ \m{G} $ acts \cite{f5} \cite{f2} 
on each of the pfaffians involved in $ (\ref{newpfaff}) $ preserving
the Quillen metrics and connections. Hence, there is an action of $ \m{G} $ 
on the
phase pfaffian $ \w{\m{L}}_{{\rm ad}} $. It makes $ (\ref{izo}) $ 
an equivariant isomorphism. The identification can be pushed down 
to quotients. One obtains an isomorphism of fibrations:
\begin{equation}
\label{izoseif}
\begin{array}{ccc}
\w{\m{M}}_{{\rm het}}^{\otimes c_{E_8}} \ / \m{G}& 
\stackrel{\simeq}{\rightarrow} & 
\w{\m{L}}_{{\rm ad}} \ / \m{G}  \\
\Big\downarrow &  & \Big\downarrow \\
\m{P}_G \ / \m{G} & = & \m{P}_G \ / \m{G} .\\ \\ 
\end{array}  
\end{equation}
\par This correspondence provides a first link relating the heterotic
parameters to fermionic pfaffian phases. However, it does not
offer quite enough input about the string parameter moduli space since 
one obtains information about, roughly speaking, the order $c_{E_8}$ tensor power of the moduli space and
not the moduli space itself. One way to use the same construction and
decrease the order of the power would be to choose a lower rank 
representation. $ E_8 $ does not have irreducible representations in
rank lower than 248. However, reducing the group structure one can 
get lower rank representations.
\subsection{The spin pfaffian}
\par Let $ H = {\rm Spin}(16) \times {\rm Spin}(16) $. There is a copy of 
$ {\rm Spin}(16) / \mathbb{Z}_2 $ sitting as a subgroup inside $ E_8 $.
Taking its double cover one obtains a Lie group morphism $ j \colon 
{\rm Spin}(16) \rightarrow E_8 $. It is known (see \cite{mt}) that the induced 
cohomology map:
\begin{equation}
\label{cohomap}
B j^4 \colon H^4(BE_8, \ \mathbb{Z}) \rightarrow 
H^4(B{\rm Spin(16)}, \ \mathbb{Z}) 
\end{equation}
is an isomorphism.
\par We denote by $ \m{P}_H $ the family of pairs $ (A, \ g) $ of flat
connections on the trivial $H$-bundle over $ E $ and  flat metrics on $E$. 
There is a projection map:
\begin{equation}
\label{sigmmamap}
\sigma \colon \m{P}_H \rightarrow \m{P}_G
\end{equation}
assigning to each $H$-connection its associated $G$-connection through:
$$ j \times j \colon {\rm Spin}(16) \times {\rm Spin}(16)
\rightarrow \eeight $$ 
\par We use the earlier geometric data on $ \m{P}_H $ taking as 
representation $ \rho $ the standard 32-dimensional representation available
for H. This will provide a pfaffian since it is actually a real representation.
The space of gauge classes of B-fields associated to the $ \m{P}_H $ family
makes the top space of the Chern-Simons circle bundle 
\begin{equation}
\label{csforh}
\m{R}_H \rightarrow \m{P}_H . 
\end{equation} 
Moreover, due to $ (\ref{cohomap}) $ the circle bundle $ (\ref{csforh}) $ 
is exactly the pull-back of the Chern-Simons bundle from $ \m{P}_G $ 
through projection $ \sigma $:
\begin{equation}
\begin{array}{ccccccc}
\m{R}_H & \simeq & \sigma ^* \m{R}_G & \rightarrow & \m{R}_G & \simeq & 
\w{\m{M}}_{{\rm het}} \\
\Big\downarrow & & \Big\downarrow & & \Big\downarrow & & \Big\downarrow \\
\m{P}_H & = & \m{P}_H & \stackrel{\sigma}{\rightarrow} & \m{P}_G & = & 
\m{P}_G.  \\
\end{array}
\end{equation}      
The Chern-Simons connection on $ \m{R}_H $ gets identified with the
pull-back through $ \sigma $ of the Chern-Simons connection on 
$ \m{R}_G $.
\par Next, we build the fermionic combination of pfaffians corresponding
to the 32-dimensional representation $ \rho $ and $ \alpha=1, 
\beta = 0  $:
$$ \m{L}_{\rho} = \ 
{\rm Pfaff}( \m{D}^{{\rm gauge}}_{(A, \ g)} )  \otimes \ 
{\rm Pfaff} ( \m{D}^{(Z/Y)}_g) \  
\otimes \ {\rm Pfaff} ( \m{D}^{{\rm grav}}_g)
^{\otimes -1} .
$$
It can be simplified to:
\begin{equation}
\label{phopfaf}
\m{L}_{\rho} = \   
{\rm Pfaff}( \m{D}^{{\rm gauge}}_{(A, \ g)} ) \
\otimes \ {\rm Pfaff} ( \m{D}_g) ^{ \otimes -8}  
\end{equation}
Employing Theorem $ \ref{ancancel} $, one obtains:
\begin{coro}
There is a circle bundle isomorphism (unique up to multiplication by a 
unitary complex number) 
between the phase pfaffian combination 
$ (\ref{phopfaf}) $ and the pull-back to $ \m{P}_H $ of the raw 
$ \eeight $ heterotic bundle $ (\ref{proj1}) $:
\begin{equation}
\label{bigdia}
\begin{array}{ccccc}
\w{\m{L}}_{\rho} & \simeq & \sigma ^* \w{\m{M}}_{{\rm het}} & \rightarrow & 
\w{\m{M}}_{{\rm het}} \\
\Big\downarrow & & \Big\downarrow & & \Big\downarrow \\
\m{P}_H & = & \m{P}_H & \stackrel{\sigma}{\rightarrow} & \m{P}_G \\
\end{array}
\end{equation}
transporting the Chern-Simons connection on  $ \w{\m{M}}_{{\rm het}} $ back to 
the Quillen connection existing on the pfaffian phase bundle. 
\end{coro}
\begin{proof}
We verify the conditions in Theorem $ \ref{ancancel} $. 
Let $ \xi $ be the generator in 
$ H^4(B{\rm Spin}(16), \Zee) \simeq \Zee $. If $ \theta $ is the
generator in $ H^4(BE_8, \Zee) $ then, as mentioned before,  under the map
$ j \colon {\rm Spin}(16) \rightarrow E_8 $, we get an isomorphism:
$$ Bj^4  H^4(BE_8, \Zee) \rightarrow  H^4(B{\rm Spin}(16), \Zee) $$ 
with $ Bj^4(\theta) = \xi $. The Chern-Simons cocycle of level
$ \lambda = (\theta, \theta) \in H^4(BG, \Zee) $ pulls then back to 
the Chern-Simons theory corresponding to level $ \w{\lambda} =
( \xi, \xi) \in H^4(BH, \Zee) $.    
\par We just have to check then the conditions needed in Theorem 
$ \ref{ancancel} $. Let $ s \colon {\rm Spin}(16) \rightarrow {\rm SU}(16) $ 
be the 
complexification of the standard 16-dimensional real representation
of $ {\rm Spin}(16) $. It induces a cohomology map:
$$ 
Bs^4 \colon H^4(BSU(16), \Zee) \rightarrow  H^4(B{\rm Spin}(16), \Zee) 
$$
with $ Bs^4(c_2) = 2 \cdot \xi $. The 32-dimensional representation 
$ \rho = ( s, s )  \colon H \rightarrow SU(16) \times SU(16) \hookrightarrow
SU(32) $ induces then:
$$ 
B\rho^4 \colon H^4(BSU(32), \Zee) \rightarrow H^4(BSU(16), \Zee) \times 
H^4(BSU(16), \Zee) \stackrel{Bs^4 \times Bs^4}{\rightarrow} H^4(BH, \Zee) 
$$
with $ B\rho^4(c_2) = (2 \cdot \xi, 2 \cdot \xi) = 2 \w{\lambda} $.  
\end{proof}
\vspace{.1in}
The symmetry group $ \m{G}_H $ acts on $ \m{P}_H $. It preserves the geometric
data necessary for building the pfaffians, and determines \cite{f5} \cite{f2}
an unitary action on  
$$  
{\rm Pfaff}( \m{D}^{{\rm gauge}}_{(A, \ g)} ) \
\otimes \ {\rm Pfaff} ( \m{D}_g) ^{ \otimes -8}. 
$$ 
Moreover, any H-gauge transformation induces a G-gauge transformation and
therefore we have a group morphism:
$$
i \colon \m{G}_H \rightarrow \m{G}.
$$
The two actions commute, making the isomorphisms in 
diagram $ (\ref{bigdia}) $ equivariant. There exist then
a circle fibration identification at quotient level:
\begin{equation}
\begin{array}{ccccc}
\w{\m{L}}_{\rho} \ / \m{G}_H& \simeq & \rho ^* 
\left ( \w{\m{M}}_{{\rm het}} \ / \m{G} \right ) & \rightarrow & 
\w{\m{M}}_{{\rm het}} \ / \m{G} \\
\Big\downarrow & & \Big\downarrow & & \Big\downarrow \\
\m{P}_H \ / \m{G}_H & = & \m{P}_H \ / \m{G}_H  & 
\stackrel{\sigma}{\rightarrow} & \m{P}_G \ / \m{G} \\ \\
\end{array} 
\end{equation}
Therefore, the  heterotic fibration pulls-back to recover the spin pfaffian. 
Summarizing the information obtained throughout this section we can say: 
\begin{prop} \
\label{e8case} 
\begin{enumerate}
\item The moduli space of $ G=E_8 \times E_8 $ heterotic parameters, 
$ \m{M}_{{\rm het}} $ makes the total space of a circle
fibration 
\begin{equation}
\label{seifert}
\m{M}_{{\rm het}} \rightarrow \m{P}_G \ / \m{G} .
\end{equation}
This is obtained by factoring out the action of the symmetry group $ \m{G} $ 
from circle bundle $ \w{\m{M}}_{{\rm het}} \rightarrow \m{P}_G $. 
\item As an equivariant model, $ (\ref{seifert}) $ represents a root of 
order $ c_{E_8} =30  $ for the adjoint pfaffian phase fibration:
$$ \w{\m{L}}_{{\rm ad}} \ / \m{G} \rightarrow \m{P}_G \ / \m{G} . $$ 
\item The pull-back of $ (\ref{seifert}) $ through the $ \sigma  $ projection of $(\ref{sigmmamap})$ recovers
the pfaffian phase corresponding to the $32$-dimensional representation
of ${\rm Spin}(16) \times {\rm Spin}(16)$.
\end{enumerate}
\end{prop}
\vspace{.1in}
We finish this section with a short note about the 
$ G = {\rm Spin}(32) / \Zee_2 $ case. This can be handled similarly 
to the $E_8 \times E_8 $ case. However, there are a couple differences 
we enumerate here. First of all, as mentioned before, the bundle $P$ is
chosen such that it carries a vector structure. Since we are working over
the 2-torus, $P$ must be topologically trivializable. We choose the 
trivial $ {\rm Spin}(32) $ bundle $ \w{P} $ as a lift for $P$. 
\par The Lie group projection $ {\rm Spin}(32) \rightarrow  {\rm Spin}(32)
 / \Zee_2 $ induces an isomorphism at Lie algebra level. Chern-Simons
cocycles for $ {\rm Spin}(32) / \Zee_2 $-connections on $P$ are then defined
(see section $ \ref{csss} $) using 
the Killing form of  $ {\rm spin}(32) $. The 
symmetry group $ \m{G} $ is made out of liftable symmetries, namely
those automorphisms of $ P $ which can be lifted to automorphisms on
$ \w{P} $. 
\par Based on this framework the entire $ E_8 \times E_8 $ discussion
in this section can be adapted to work for $ G ={\rm Spin}(32)
 / \Zee_2 $. The spin pfaffian arrives via the standard 32-dimensional 
spin representation $ \rho $ 
of $ H = {\rm Spin}(32) $. The analog of Proposition
$ \ref{e8case} $ can be formulated:
\begin{prop} \
\label{spincase}
\begin{enumerate}
\item The moduli space of $ G= {\rm Spin}(32) / \Zee_2  $ 
heterotic parameters, 
$ \m{M}_{{\rm het}} $ makes the total space of a circle
fibration 
\begin{equation}
\label{spinseifert}
\m{M}_{{\rm het}} \rightarrow \m{P}_G \ / \m{G} .
\end{equation}
This can be regarded as an equivariant model. $ (\ref{spinseifert}) $ 
is obtained by factoring out the action of the symmetry group $ \m{G} $ 
from a circle bundle $ \w{\m{M}}_{{\rm het}} \rightarrow \m{P}_G $. 
\item As an equivariant model, $ (\ref{spinseifert}) $ represents a root of 
order $ 30  $ for the adjoint pfaffian phase fibration:
$$ \w{\m{L}}_{{\rm ad}} \ / \m{G} \rightarrow \m{P}_G \ / \m{G} . $$ 
\item The pull-back of $ (\ref{spinseifert}) $ through 
$ \m{P}_H \rightarrow \m{P}_G  $ recovers
the pfaffian phase corresponding to the $32$-dimensional representation
of ${\rm Spin}(32)$. Moreover, the pull-back commutes with the actions
of the symmetry group $ \m{G} $, creating an identification of
equivariant models.
\end{enumerate}
\end{prop}  
\vspace{.1in}
One can interpret this description of $ {\rm Spin}(32) / \Zee_2 $ heterotic
parameters along the lines set by Witten \cite{MR1748791} in his analysis 
of world-sheet anomaly cancellation. Our definition of B-fields 
as differential 2-cochains recovers the properties Witten enlists in his
discussion on the nature of these fields. The 
$ {\rm Spin}(32) / \Zee_2 $ analog 
of pfaffian combination $ (\ref{phopfaf}) $ represents the supergravity
fermionic anomaly circle bundle and its cancellation mechanism can be 
seen through ideas similar to world-sheet anomaly cancellation. 
\par Evaluating the supergravity 
path integral for a fixed Riemann surface (namely the choice 
$ E \subset E \times \mathbb{R}^8 $) one obtains a product combination:
\begin{equation}
\label{pathint}
\m{L}^{-1}_{\rho}(A, g) \ {\rm exp} \left ( i \int_E B \right ).
\end{equation}   
Here, $ \m{L}^{-1}_{\rho}(A, g) $ represents the inverse of the pfaffian 
combination 
$ (\ref{phopfaf}) $ while the second factor is the ``holonomy''
of the B-field B along the 2-torus $E$. Both factors represent
sections of non-trivial line bundles. The first one makes a section
in the inverse of the line bundle $ \m{L}_{\rho} $. The second 
quantity, as discussed in section $ \ref{fir} $, takes values in the total space of 
a circle bundle, the 
Chern-Simons bundle $ \m{R} $. Both bundles carry unitary connections. 
In order for the path integral combination
$ (\ref{pathint}) $ to be a complex number, not just a point in the total
space of a tensor bundle, one needs to covariantly 
trivialize the tensor combination
$$ \m{L}_{\rho} ^{-1} \otimes \m{R}. $$
A trivialization can be achieved only if the two phase bundles can be made
to cancel each other. In other words, if a smooth isomorphism 
$ \m{\w{L}}_{\rho} \simeq \m{R} $ identifying the two circle holonomies
does exist. The third statement in
Proposition $ \ref{spincase} $ provides exactly such an isomorphism.  
%
%
%
%
\section{Moduli of Heterotic Data and the Character Fibration}
Let $G$ be one of the two Lie group choices $ E_8 \times E_8 $ 
or $ {\rm Spin}(32)/ \mathbb{Z}_2 $. 
We analyze the family $ \m{P}_G $ up to gauge equivalence. 
The family is made out of pairs consisting of flat $G$-connections $A$ and 
flat metrics $g$ over the 2-torus $E$. By convention, if
$ G = {\rm Spin}(32)/ \mathbb{Z}_2 $ we consider only those connections
that can be lifted to $ {\rm Spin}(32) $. The gauge group $ \m{G} $ is
represented by automorphisms $ \varphi $ of the bundle covering 
orientation preserving diffeomorphisms $ \bar{\varphi} $ on $ E $ (as 
mentioned earlier if $ G = {\rm Spin}(32)/ \mathbb{Z}_2 $ one considers
only liftable automorphisms). $ \m{G} $  enters a short exact sequence:
$$ \{ 1 \} \rightarrow \m{G}(P) \hookrightarrow \m{G} \rightarrow 
{\rm Diff}^+(E) \rightarrow \{ 1 \} $$
and acts on 
$\m{P}_G$ as $ \varphi \cdot (A, \ g) = ( \varphi^*A, \ \bar{\varphi}^*g)$.   
\par The space of flat metrics splits as:
 $$ {\rm Met}_o(E) = {\rm Conf}(E) \times \mathbb{R}_+^* $$
any given metric producing a conformal class and a volume. It is
known that, on a 2-torus, a conformal class produces a complex
structure.  The group of orientation preserving diffeomorphisms
$ {\rm Diff}^+(E) $ acts on $ {\rm Met}_o(E) $ leaving the volume component
invariant. Therefore:
$$  {\rm Met}_o(E)  /  {\rm Diff}^+(E) = 
\left ( {\rm Conf}(E) /  {\rm Diff}^+(E)  \right ) 
\times \mathbb{R}^*_+ . $$
Two conformal classes determine isomorphic complex structures
on $E$ if and only if they can be transformed one into the other through
an orientation-preserving diffeomorphism. The orbit space:
$$ \m{M}_E = {\rm Conf}(E)  /  {\rm Diff}^+(E) $$ 
represents then the family of isomorphism classes of elliptic curves. For 
practical reasons we shall identify a complex structure on a 2-torus by
a complex number $ \tau $ inside the upper half-plane $ \mathcal{H} $. The 
corresponding elliptic curve
is obtained by factoring the lattice 
$ \mathbb{Z} \oplus \tau \mathbb{Z} $ out of the complex plane $ \mathbb{C} $.
Two points in the upper half-plane $ \tau_1 $ and $ \tau_2 $ determine 
isomorphic complex structures if and only if:
$$ \tau_2 = \frac{a \tau_1+b}{c \tau_1 + d } \ \ {\rm with} \ \
\left ( 
\begin{array}{cc}
a & b \\
c & d \\
\end{array}
\right ) \in \ {\rm SL}(2, \ \mathbb{Z}). $$
One can therefore describe the moduli space of isomorphic classes 
of elliptic curves as a complex orbifold:
\begin{equation}
\m{M}_E = {\rm PSL}(2, \mathbb{Z}) \backslash \mathcal{H}.  
\end{equation}
\par Turning to connections, the family of gauge equivalence classes of
$G$-connections over a 2-torus is identified to $ {\rm Hom}(\pi_1(E) , \ 
G) / G $, the group G acting on the space of homomorphisms by conjugation.
A gauge class can therefore be seen as a pair of commuting elements
in $G$, up to simultaneous conjugation. For simply connected groups 
this picture can be refined (see for example \cite{morgan1}). Any two commuting elements in 
$G$ can be 
simultaneously conjugated inside a maximal torus. This also holds for
liftable $G$-connections if $G$ is connected but not necessarily 
simply-connected. We make a choice
of maximal torus $ T \hookrightarrow G $. Let $ W = W(T, G) $ be the
associated Weyl group  and $ \bold{t}_{\mathbb{R}} \subset \bold{g} $ be 
the Cartan sub-algebra. There is then a 1-to-1 correspondence:
\begin{equation}
\label{rel1}
{\rm Hom}(\pi_1(E) , \ G) \ / G \simeq {\rm Hom}(\pi_1(E) , \ T) \ / W .
\end{equation}
We consider the coroot lattice $ \Lambda_G \subset \bold{t}_{\mathbb{R}} $. 
The Weyl group acts on $ \Lambda_G $ preserving the set of coroots. There
is an identification $ T \simeq U(1) \otimes _{\mathbb{Z}} \Lambda_G $ 
commuting with the Weyl action. Under this identification:
\begin{equation}
\label{rel2}
{\rm Hom}(\pi_1(E), \ T) \simeq
{\rm Hom}(\pi_1(E), \ U(1) \otimes _{\mathbb{Z}} \Lambda_G) \simeq
{\rm Hom}(\pi_1(E), \ U(1)) \otimes _{\mathbb{Z}} \Lambda_G. 
\end{equation}  
The set $ {\rm Hom}(\pi_1(E), \ U(1)) $ represents the family of gauge 
equivalence classes of flat hermitian line bundles. In the presence of a 
complex structure $E_{\tau} $ this can be regarded as  
$ {\rm Pic}^o(E_{\tau}) $. Hence, using $ (\ref{rel1}) $ and $ (\ref{rel2}) $, 
on a complex torus $E_{\tau} $ the space of gauge classes of $G$-connections 
can be identified to the complex orbifold:
\begin{equation}
\label{orbi}
\left ( {\rm Pic}^o(E) \otimes _{\mathbb{Z}} \Lambda_G \right ) / W. 
\end{equation}  
For $ G = E_8 \times E_8 $, this is a product of two identical 8-dimensional 
complex weighted projective spaces \cite{loo}. If 
$ G = {\rm Spin} (32) / \Zee_2 $ then $ (\ref{orbi}) $ makes a quotient of
a 16-dimensional complex weighted projective space by $ \Zee_2 \times \Zee_2 $.
Bringing things together one concludes:
\begin{teo}
\label{t}
(\cite{morgan1} \cite{loo})
The moduli space $ \m{M}_{E, \ G} $ of equivalence classes of flat 
$G$-connections over elliptic curves can be given the structure of a 
17-dimensional complex variety. It fibers over the moduli space of 
elliptic curves:
\begin{equation}
\label{eee1}
\m{M}_{E, \ G} \rightarrow \m{M}_E 
\end{equation}
all fibers being 16-dimensional compact complex orbifolds. 
Moreover, there is a diffeomorphism of stratified spaces:
\begin{equation}
\label{split33}
\m{P}_G /  \m{G} \simeq \m{M}_{E, \ G} \times \mathbb{R}^*_+ ,
\end{equation}  
the second factor on the right-hand side representing the volume
of the flat metric on $E$.
\end{teo} 
Let us recapitulate now the facts of previous sections. Following section $\ref{sec3} $,
in both relevant cases $ G= \eeight  $ and $ G = \spintt  $, the moduli space of heterotic 
string data  $ \mhet $ is obtained by factoring out the action of the symmetry group 
$ \m{G} $ in the principal circle bundle of $(\ref{proj1})$:
\begin{equation}
\label{again}
\w{\m{M}}^G_{{\rm het}} \rightarrow \m{P}_G.
\end{equation} 
As discussed in section $ \ref{sec4} $, this bundle is endowed 
with a Chern-Simons $S^1$-connection. Furthermore, the quotient space $\m{P}_G/\m{G} $ splits as
in $(\ref{split33})$. Since the gravitational Chern-Simons invariant does not vary under 
flat conformal transformations, the bundle $(\ref{again})$ can be trivialized along the 
volume direction and therefore the space $\mhet$ can be seen as the total space of 
a $C^{\infty}$ fibration with fibers isomorphic to $\mathbb{R}^*_+ \times S^1=\mathbb{C}^*$:
\begin{equation}
\label{fibr}
\mhet \rightarrow \m{M}_{E,G}.  
\end{equation}
One can see at this point that $\mhet$ inherits a natural structure of analytic space.
Indeed, looking at the $\mathbb{C}^*$-fibration that covers $ (\ref{fibr})$, one sees that 
this fibration is associated to a $C^{\infty}$ fibration with complex lines over a smooth 
complex manifold which is endowed with a connection and compatible hermitian metric.
These features induce a natural complex structure on the total space of the line bundle
which descends to an analytic structure on $\mhet$. Moreover, in this setting, fibration 
$ (\ref{fibr}) $ becomes  a holomorphic (Seifert) $\mathbb{C}^*$-fibration. The analysis of 
this holomorphic structure of $ (\ref{fibr}) $ is the central goal of this paper. 
\par In future considerations we shall use an abstract, coordinate oriented, 
working model for 
$ \m{M}_{E, \ G} $. Namely, let:
$$ V_G \ \stackrel{{\rm def}}{=} \ \mathcal{H} \times \left ( 
\mathbb{C} \otimes_{\mathbb{Z}} \Lambda_G \right ). $$
Every pair $(\tau,z) \in V_G$ can be seen to determine an elliptic curve $E_{\tau} $ together 
with a flat $G$-connection connection on $E_{\tau}$. Indeed, one can take:
$$E_{\tau} \ = \ \mathbb{C} / \mathbb{Z} \oplus \tau \mathbb{Z} $$ and 
then the factor $z$ can be seen to determine an element in:   
\begin{equation}
\label{orbi2}
\left ( {\rm Pic}^o(E_{\tau}) \otimes _{\mathbb{Z}} \Lambda_G \right ) / W \ \simeq \ 
 \left ( E_{\tau} \otimes _{\mathbb{Z}} \Lambda_G \right ) / W  
\end{equation} 
which according to the previous arguments define a corresponding $G$-flat connection
on $E_\tau$. Let $ {\rm pr}_1 \colon V_G \rightarrow \mathbb{H}$ be the projection on the 
first factor. For each $ \tau \in \mathcal{H} $ take:
$$ \mathcal{L}_{\tau} \ \colon = \ \{\tau\} \times \left ( 
\left ( \mathbb{Z} \oplus \tau \mathbb{Z} \right ) \otimes \Lambda_G \right ) 
\subset {\rm pr}_1^{-1} (\tau). $$
Then $ \mathcal{L}_{\tau} $ represents a family of $32$-dimensional lattices 
sitting fiber-wise inside the fibration ${\rm pr}_1$. 
\begin{dfn}
Let $ \Pi_G $ be the group of holomorphic automorphisms of the fibration ${\rm pr}_1$ which preserve 
the lattice family $\m{L} $ and cover $ {\rm PSL}(2, \Zee) $ transformations on $ \mathcal{H} $. 
\end{dfn} 
\noindent It turns out that two elements $ (\tau, z) $ and $ (\tau', z') $ of 
$ \mathcal{H} \times \Lambda_{\Cee} $ determine isomorphic pairs of elliptic curves 
and flat $G$-connections if and only if they can be transformed one into another through an isomorphism 
in $ \Pi_G $. In this sense, the analytic space:
\begin{equation}
\label{model}
\Pi_G \ \backslash \left ( \mathcal{H} \times \Lambda_{\Cee} \right )   
\end{equation}     
can be seen as the moduli space of pairs of elliptic curves and flat 
$G$-bundles\footnote{Again, in the case $G= {\rm Spin}(32)/ \mathbb{Z}_2 $, one considers only 
${\rm Spin}(32)$-liftable connections} $\m{M}_{E, G}$. Also, there is a canonical projection
$ p \colon \Pi_G \rightarrow {\rm PSL}(2, \mathbb{Z}) $ which commutes with ${\rm pr}_1$ 
inducing a fibration:
$$ \Pi_G \backslash V_G \rightarrow {\rm PSL}(2, \mathbb{Z}) \backslash \mathcal{H}. $$
This is the model for projection  $(\ref{eee1})$ in theorem $\ref{t}$.
\par Let us single out the following three particular subgroups of $\Pi_G$:
\begin{itemize}
\item $ \m{S}_G \ = \  \left \{ \psi \in \Pi_G \ \Big\vert \ 
\psi(\tau,  z ) \ = \ \left ( 
\frac{a \tau + b}{c \tau + d} , \ \frac{z}{c\tau+d}  
\right ) \   {\rm where} \ 
\left (
\begin{array}{cc}
a & b \\
c & d \\
\end{array}
\right )
\in {\rm SL}(2,\mathbb{Z}) \ 
\right \} $
\item $ \m{T}_G  \ = \ 
\left \{ \psi \in \Pi_G  \ \Big\vert \ 
\psi ( \tau, z) = ( \tau, z + q_1 + \tau q_2 ) \ \ {\rm where} \ \ (q_1, q_2) \in 
\Lambda_G \oplus \Lambda_G   
\right \} $
\item $ \m{W}_G = \ \left \{ \psi = {\rm id} \oplus f  \in \Pi_G
\ \Big\vert \ f \in O(\Lambda)  \right \} $.  
\end{itemize}      
The three subgroups $ \m{S}_G$, $\m{T}_G$ and $\m{W}_G $ generate 
the entire $ \Pi_G $. In addition, note that $ \m{S}_G  \cap \m{W}_G  =  \{ \pm {\rm id} \} $
and that $ {\rm Ker}(p) $ is 
generated by $ \m{W}_G $ and $ \m{T}_G  $. 
One concludes from these facts that $ \Pi_G $ is a semi-direct product: 
\begin{equation}
\Pi_G \ = \ \m{T}_G   \rtimes 
\left ( \m{W}_G  \times_{\{ \pm {\rm id} \}} \m{S}_G  \right ).  
\end{equation} 
The model $ (\ref{model}) $ allows us now to describe easily
holomorphic line bundles over the moduli space of of elliptic curves
and flat $G$-bundles, $\m{M}_{E, \ G} $. Such bundles over a complex
orbifold are best described in terms of equivariant line bundles over
the cover. These equivariant bundles are line bundles $ L \rightarrow V_G $ 
where the action of the group $ \Pi_G $ on the base is given a lift to 
the fibers. All holomorphic line bundles over $ V_G $ are trivial and
a lift of the action to fibers can be obtained through a set of 
automorphy factors
$ (\varphi_a)_{a \in \Pi_G} $ with 
$ \varphi_a \in H^0(V_G, \m{O}^*_{V_G}) $ satisfying:
$$ \varphi_{ab}(x) = \varphi_a(b \cdot x) \varphi_b(x). $$
Such a set makes a 1-cocycle $ \varphi $ 
in $ Z^1(\Pi_G, H^0(V_G, \m{O}^*_{V_G})) $.
Two automorphy factors provide isomorphic line bundles on 
$ \m{M}_{E, \ G} $ if and only if they determine the same cohomology class
in $ H^1(\Pi_G, H^0(V_G, \m{O}^*_{V_G})) $. To state this rigorously, 
there is a canonical map $ \phi $ entering the following exact sequence:
\begin{equation}
\{1\} \rightarrow H^1(\Pi_G, H^0(V_G, \m{O}^*_{V_G})) 
\stackrel{\phi}{\rightarrow} H^1(\m{M}_{E, \ G}, \ \m{O}^*_{\m{M}_{E, \ G}}) 
\stackrel{}{\rightarrow} H^1(V_G, \ \m{O}^*_{V_G}) \simeq \{1\} .
\end{equation}      
\par There is an important holomorphic line bundle living on $ \m{M}_{E, \ G} $.
Leaving aside more advanced interpretation, one can minimally define this bundle 
as the fibration with lines supporting the character function of $ \Lambda_G$. 
Indeed, in both cases $ G=E_8 \times E_8 $ and 
$ G={\rm Spin}(32) / \mathbb{Z}_2 $, the coroot lattice $  \Lambda_G $ is 
positive definite, unimodular and even and therefore, one can define then 
an associated holomorphic theta function 
(see \cite{kac1} \cite{mum} for details) :
\begin{equation}
\Theta_{\Lambda_G} \colon V_G \rightarrow \mathbb{C}, \ \ \ \ 
\Theta_{\Lambda_G}(\tau, \ z) = \ \sum_{\gamma \in \Lambda} \ 
e^{ \pi i ( 2 (z,\gamma) + \tau (\gamma, \gamma))}.  
\end{equation}
The pairing appearing above represents the bilinear complexification of the 
integral pairing on $ \Lambda_G $. The ${\bf \Lambda_G}$-${\bf character}$ 
function can be written then as a quotient of $ \Theta_{\Lambda_G} $:
\begin{equation}
\label{char}
B_{\Lambda_G} \colon V_G \rightarrow \mathbb{C}, \ \ \ \ B_{\Lambda}
(\tau, \ z) = 
\frac{\Theta_{\Lambda_G}(\tau, \ z)}{\eta(\tau)^{16}}. 
\end{equation}
Here, $ \eta(\tau) $ is Dedekind's eta function:
$$ \eta(\tau) = e^{\pi i \tau /12} \prod_{m=1}^{\infty} \ \left ( 
1- e^{2 \pi i m \tau} \right ),$$
which is an automorphic form of weight $ 1/2$ and multiplier system 
given by a group homomorphism 
$ \chi \colon {\rm SL}_2(\Zee) \rightarrow \Zee / 24 \Zee $, in the
sense that \cite{sig}:
$$ \eta ( \gamma \cdot \tau ) \ = \ 
\chi(\gamma) \  \sqrt{c \tau+d} \ \eta(\tau) \ \ {\rm for} \ \ 
\gamma = \left (
\begin{array}{cc}
a & b \\
c & d \\
\end{array}
\right ) \in {\rm SL}_2(\Zee). 
$$
The character terminology for $ ( \ref{char} ) $ is justified by its role in the
representation theory of infinite-dimensional Lie algebras. The function 
$ B_{\Lambda_G} $ is the
zero-character of the level $l=1$ basic highest weight representation of the Kac-Moody
algebra associated to $G$ (see \cite{kac1} for details).
\par The function $B_{\Lambda_G}$ obeys the following transformation properties:
\begin{prop}   
\label{kac} 
(\cite{kac1}) Under the action of the group $ \Pi_G $, the $\Lambda_G$-character 
function $ (\ref{char}) $ transforms as :
$$  B_{\Lambda_G} \left ( g \cdot ( \tau, z) \right )  \ = \ 
\varphi^{{\rm ch}}_g (\tau, z) \cdot B_{\Lambda_G}( \tau, z), \ \ g \in \Pi_G. $$
The factors $ \varphi^{{\rm ch}}_a$ can be described on the generators of the group 
$\Pi_G$ as:
\begin{itemize}
\item $ \varphi ^{{\rm ch}}_g (\tau, z) = 
e^{\pi i(-2<q_2, z>-\tau<q_2,q_2> )}$  for $g\in \m{T}_G$ induced by $ (q_1, q_2) \in \Lambda_G
\oplus \Lambda_G $.  
\item $ \varphi ^{{\rm ch}}_g = e^{\frac{\pi i c <z, z>}{c \tau + d}} $ 
 for $g\in \m{S}_G$ induced by $ \left (
\begin{array}{cc}
a & b \\
c & d \\
\end{array} \right ) 
\in {\rm SL}(2, \mathbb{Z}) $. 
\item $ \varphi^{{\rm ch}}_g = 1 $ for $ g \in \m{W}_G $.
\end{itemize}
\end{prop}
\noindent The function $ B_{\Lambda_G} $ descends therefore to a holomorphic section of a holomorphic 
$ \mathbb{C} $-fibration:
\begin{equation}
\label{thetaline}
\m{Z} \rightarrow \Pi_G \backslash \left ( \m{H} \times \Lambda_{\mathbb{C}} \right ) = \m{M}_{E,  G}.
\end{equation}
for which a set of automorphy factors are given by the $ \varphi^{{\rm ch}}_g $ in proposition $\ref{kac}$.    
\par The main result of this paper asserts:
\begin{teo}
\label{main}
There exist a holomorphic isomorphism of $\mathbb{C}^*$-fibrations 
between $ (\ref{fibr}) $ and the $\mathbb{C}^*$ the $ \Lambda_G$-character 
fibration $ (\ref{thetaline}) $. The isomorphism is unique up to twisting with a 
non-zero complex number.  
\end{teo}  
The rest of the paper is dedicated to proving this theorem.
%
%
%
%
\section{Proof of Theorem $ \ref{main} $.}
\label{prooof}
Our strategy is as follows. We shall exploit the link described in section $ \ref{sec4} $ between 
the heterotic fibration:
\begin{equation}
\mhet \rightarrow \m{P}_G / \m{G} 
\end{equation}
and the pfaffian combinations $ (\ref{newpfaff}) $ and $ (\ref{phopfaf}) $ in order to explicitly 
compute a set of automorphy factors $ \varphi^{{\rm het}}_g, \ g \in \Pi_G $ for $ (\ref{fibr}) $. 
Then we shall compare $ \varphi^{{\rm het}}_g $ such computed with the automorphy 
factors $ \varphi^{{\rm ch}}_g $ of the $\Lambda_G$-character fibration $ (\ref{thetaline}) $ 
provided by proposition $ \ref{kac}$.
\subsection{Determinants of Flat Line Bundles over Elliptic Curves}
One can regard a pair consisting of an elliptic curve and a flat hermitian
line bundle as an element $ (\tau, \ u ) \in \m{H} \times \mathbb{C} $.
The imaginary complex number $ \tau $ determines a complex structure
$ E_{\tau} $ on the 2-torus whereas a complex number $ u = \alpha \tau + 
\beta \in \mathbb{C} $ determines a flat bundle of holonomy map:
$$ h \colon \pi_1(E) \simeq \Zee \oplus \Zee \rightarrow U(1) , \ \ 
h(m, n) = e^{2 \pi i ( \alpha m + \beta n)}. $$ 
The projection 
$ {\rm pr}_1 \colon \mathbb{H} \times \mathbb{C} \rightarrow \mathbb{H}  $
carries fiber wise a continuous family of lattices 
$U_{\tau} = \Zee \oplus \tau \Zee \subset \m{H} \times \mathbb{C}$. 
We denote by $ \Pi_{U(1)} $ the group of
automorphisms of projection $ {\rm pr}_1 $, preserving the lattice family $U_{\tau}$ 
and covering $ {\rm PSL}(2, \Zee) $ transformations on the base. Two
pairs $ ( \tau_1, u_1) $ and $ ( \tau_2, u_2) $ determine isomorphic
pairs of elliptic curves and flat line bundles if and only if they 
belong to the same orbit under the action of $ \Pi_{U(1)} $. The moduli 
space of elliptic curves and flat hermitian line bundles can then be seen as:
\begin{equation}
\label{modelnou}
\m{M}_{E,  U(1)} = \Pi_{U(1)} \backslash \left ( \m{H} \times \mathbb{C} \right ) .
\end{equation}
The induced projection:
$$ \m{M}_{E, U(1)} \rightarrow \m{M}_E =   
{\rm PSL}(2, \Zee) \backslash \m{H} $$ 
makes a holomorphic torus fibration, the fiber over $ [ \tau ] \in 
\m{M}_E $ being $ {\rm Pic}^o(E_{\tau}) $.  
The symmetry group $ \Pi_{U(1)} $ is generated by two particular 
subgroups:
\begin{itemize}
\item $ \m{S}_{U(1}) \ = \ \left \{ \psi \in \Pi_{U(1)} \ \Big\vert \ 
\psi (\tau, u) \ = \ \left ( \frac{a \tau + b}{c \tau + d} , \ 
\frac{u}{c \tau + d} \ \right ) \ {\rm where} \ 
\left (
\begin{array}{cc}
a & b \\
c & d \\
\end{array}
\right ) \in {\rm SL}(2, \mathbb{Z})
\ \right \} \ \simeq \ {\rm SL}(2, \mathbb{Z})$
\item $ \m{T}_{U(1)} \ = \ 
\left \{ \psi \in \Gamma_{U(1)}  \ \Big\vert \ 
\psi ( \tau, u) = ( \tau, u + a_1 + a_2 \tau) \ \ {\rm where} \ \ 
(a_1, a_2) \in \Zee \oplus \Zee   
\right \} \ \simeq \ \mathbb{Z} \oplus \mathbb{Z} .$
\end{itemize}
The group $ \Pi_{U(1)} $ can then be generated as a semi-direct product $ \m{T}_{U(1)} \rtimes \m{S}_{U(1)}$.
\par Consider then the odd spin structure (trivial double cover of the frame bundle)
on $E$. In the presence of a metric on $E$, there is then an induced Dirac operator. Moreover, this operator 
can be coupled to the flat connections creating therefore a family of elliptic differential operators. 
Interpreting the flat line bundles as holomorphic line bundles, one can see these operators as 
$ \overline{\partial}_L $ operators. The family determines then a determinant $ \mathbb{C} $-fibration 
(see \cite{f5} for more details):
\begin{equation}
\label{flatdet}
{\rm Det}(\m{D}_L) \rightarrow \m{M}_{E,  U(1)} 
\end{equation}
carrying a natural determinant section $ {\rm det} \colon \m{M}_{E,  U(1)} 
\rightarrow {\rm Det}(\m{D}_L) $. One can use then the model $(\ref{modelnou})$ to describe a set of 
automorphy factors for this fibration.
\begin{lem}
\label{flatdetlem}
The determinant fibration $ ( \ref{flatdet}) $ can be obtained by factoring the trivial holomorphic line bundle 
over $ \mathbb{H} \times \mathbb{C} $ by a set of automorphy factors $\varphi_g$, $g \in \Pi_{U(1)}$ with:
\begin{enumerate}
\item $ \varphi_g(\tau, u) = (-1)^{q_1+q_2} \ 
e^{\pi i(-2uq_2-\tau q_2^2)} $ for 
$ g \in \m{T}_{U(1)} $ determined by $(q_1, q_2) \in \Zee \oplus \Zee $ 
\item $ \varphi_g(\tau, u) = \chi(m)^2 \ e^{\pi i \ \frac{cu^2}{c \tau + d}} $ 
for $g \in \m{S}_{U(1)}$ associated to $ m = \abcd \in {\rm SL}(2, \Zee) $.
\end{enumerate}
Moreover, the determinant section $ {\rm det} $ in obtained in this setting as the factorization of the holomorphic function:
$$ f \colon \m{H} \times \mathbb{C} \rightarrow \mathbb{C}, \ \ 
f( \tau, u) = \frac{\vartheta_1(\tau, u)}{\eta(\tau)}. $$
Here $ \vartheta_1(\tau, u) $ is Siegel's twisted theta function:
$$ \vartheta_1(\tau, u) \ = \ \sum_{n \in \Zee} e^{2 \pi i \left ( u- \frac{1}{2} \right ) 
\left ( n + \frac{1}{2} \right ) +  
\pi i \tau \left (n+ \frac{1}{2} \right )^2 }.$$   
\end{lem} 
\begin{proof}
We follow arguments from \cite{f5}. The determinant line $ (\ref{flatdet}) $ can be
endowed with Quillen metric and connection compatible with the 
holomorphic structure. The Quillen norm of the determinant section 
$ {\rm det} $ 
is then computed as regularized determinant. According to \cite{ray1} \cite{f5}, one has:
\begin{equation}
\label{eq150}
\vert \vert \ {\rm det}( \tau, u) \ \vert \vert \ = 
\  \left | \
e^{ \pi i u u_2 } \  
\frac { \vartheta _1 ( \tau , u) }{ \eta ( \tau ) } \ \right | 
\end{equation}
with $ u_2 $ coming from the decomposition $ u = u_1 + \tau u_2 $. 
\par Let us consider the pull-back $ \pi ^* {\rm Det}(\m{D}_L)  $ through the 
covering map $ \pi \colon \m{H} \times \mathbb{C}  \rightarrow 
\mathcal{M}_{E,  U(1)} $. The determinant section $ {\rm det} $ pulls back 
to a new section $ \widetilde{{\rm det}} $ in $ \pi ^* {\rm Det}(\m{D}_L) $. 
The Quillen metric and connection pull-back to $ \pi ^* {\rm Det}(\m{D}_L) $ 
as well. One computes then the curvature of the pull-back connection as: 
\begin{equation}
\Omega \ = \ \overline{\partial} \partial \ {\rm log} \| \widetilde{{\rm det}} \| ^2
\end{equation}
Combining with $ (\ref{eq150}) $ we get:
\begin{equation}
\Omega \ = \ \overline{\partial} \partial \ (2 \pi i u u_2) \ 
= \ 2 \pi i \ \overline{\partial} \partial \left ( 
\frac{u^2 - u \overline{u} }{ \tau - \overline{ \tau } } \right ) 
\end{equation}
We choose now a holomorphic trivialization for $ \pi ^* {\rm Det}(\m{D}_L) $. 
It is
not possible to find a covariant trivialization (with respect to the
pull-back Quillen connection $ \nabla_Q $) since the curvature does 
not vanish. However, one can deform the connection to:
$$ \nabla ' = \nabla_Q - \mu \ , \ \mu \in 
\Omega^1(\m{H} \times \mathbb{C}). $$
and remove the curvature. In performing this deformation, we keep in mind that we need 
a holomorphic trivialization so, we wish to keep the connection 
compatible with the holomorphic structure. In order to achieve this, the deformation 1-form must satisfy:
\begin{equation}
\label{eq154}
d \mu \ = \ \Omega \ \ \ \text{and} \ \ \mu^{0,1} \ = 0 \ . 
\end{equation}
We make the choice:
\begin{equation}
\label{eq155}
\mu \ = \ 2 \pi i \partial \ \left ( 
\frac{u^2-u \overline{u}}{\tau - \overline{\tau}}
\right ) 
\end{equation}
which satisfies conditions $ (\ref{eq154}) $. The new connection $\nabla '$ is flat and therefore 
it is possible now to choose a holomorphic covariant section denoted $ \sigma $. Clearly, 
$ \widetilde{{\rm det}}  =  f \sigma $ with f holomorphic function. In connection with 
$ (\ref{eq150})$ one gets then:
$$
\vert \ f \ \vert \ = \ \left \vert 
\frac{\vartheta_1(\tau, \ z)}{\eta(\tau) }
\right \vert .
$$
Since f is holomorphic, we see that, up to multiplication by an unitary 
 complex number: 
\begin{equation}
\label{eq200}
f \ = \ \frac{\vartheta_1(\tau, \ z)}{\eta(\tau) } \ . 
\end{equation}
The automorphy factors result then from the well-known transformation
rules for $ (\ref{eq200}) $. They are (as described in \cite{sig}) :
\begin{equation}
f( \tau , u + q_1 +\tau q_2 ) = (-1)^{q_1+q_2} \ 
e^{  \pi i ( - 2uq_2- \tau q_2^2) }
 \  f( \tau , u) \ \ {\rm for} \ (q_1, q_2) \in \Zee \oplus \Zee 
\end{equation}
\begin{equation}
f \left (
\frac{ a \tau + b }{ c \tau +d } \ , \ \frac{u}{c \tau + d}
\right )
=\chi(m)^2 \  e^{  \frac{ \pi i c z^2 }{c \tau +d} } 
\ f( \tau , z) \ \ {\rm for} \ m = \abcd \in {\rm SL}(2, \Zee) . 
\end{equation}
\end{proof} 
\noindent There is a particular case of this construction that we need to also analyze. 
Suppose that we give up the coupling with flat connections and just consider the Dirac 
operators induced by the flat metrics of the elliptic curves. As showed in $ (\ref{part2}) $, 
these operators can be interpreted as skew-symmetric operators. In this setting, there exists 
a square root for the determinant, the pfaffian. This $ \mathbb{C} $-fibration:
\begin{equation}
\label{pfafff}
{\rm Pfaff}(\overline{\partial}_{\tau}) \rightarrow 
\m{M}_E = \psl \backslash \m{H} .
\end{equation}
appeared earlier in the fermionic combinations $ (\ref{phopfaf}) $ and 
$ (\ref{newpfaff}) $. Removing the flat connections from lemma 
$ \ref{flatdetlem} $ one obtains:
\begin{coro}
\label{coro1}
The holomorphic pfaffian fibration $ (\ref{pfafff}) $ can be obtained by factoring out of the
trivial line over $ \m{H} $ by a set of automorphy factors:
$$ \varphi^o_m (\tau) \ = \ \chi(m) \ \ {\rm for} \ \ m = \abcd \in 
{\rm SL}(2, \Zee). $$ 
\end{coro} 
\subsection{Automorphy Factors for Fermionic Pfaffians  }
We employ now the data of previous section to compute automorphy factors for the 
line bundles $ (\ref{phopfaf}) $ and $ (\ref{newpfaff}) $. Each of the two is a tensor 
product involving pfaffian bundles. Automorphy factors for the simplest one are given by 
Corollary $ \ref{coro1} $. In this section we compute the automorphy factors for the 
other two.
\par The first one to look at is the adjoint pfaffian:
\begin{equation}
\label{adjpfaff}
{\rm Pfaff}(\m{D}_{\rm ad}) \rightarrow \m{M}_{E, \ G}.
\end{equation}
This holomorphic line bundle is obtained as the pfaffian of the 
family of Dirac operators coupled
to flat connections on the complex vector bundle associated to the
adjoint representation of $G$. This complex vector bundle has a natural
real structure. The pfaffian carries a canonical holomorphic section
$ {\rm pfaff}_{{\rm ad}} \colon \m{M}_{E, \ G} \rightarrow 
 {\rm Pfaff}(\m{D}_{{\rm ad}}) $. Working under coordinate model:
$$ \m{M}_{E, \ G} \ = \   \Pi_G \backslash V_G $$
we can describe the bundle as follows:
\begin{lem}
\label{lemad}
The pfaffian line bundle $ (\ref{adjpfaff}) $ can be obtained by factoring 
out the trivial line over the universal cover $ V_G $ by means of automorphy 
factors $\varphi_g^{{\rm ad}}$, $g \in \Pi_G$ defined on generators as:
\begin{enumerate}
\item $ \varphi^{{\rm ad}}_g(\tau, z) = 
e^{c_G \pi i(-2<z, q_2> -\tau <q_2, q_2>)} $ for 
$g \in \m{T}_G$ associated to $ (q_1, q_2) \in \Lambda_G \oplus \Lambda_G $
\item $ \varphi^{{\rm ad}}_g(\tau, z) = \chi(m)^{N} \ 
e^{c_G \pi i \ \frac{c<z,z>}{c \tau + d}} $ 
for $g\in \m{S}_G$ associated to $ m = \abcd \in {\rm SL}(2, \Zee) $
\item $ \varphi^{{\rm ad}}_g(\tau, z) = 1 $ for $ g \in \m{W}_G $.
\end{enumerate}
Here $c_G=30 $ is the Coxeter number of the group $G$ (it takes the
same value for both choices $ E_8 \times E_8 $ and 
$ {\rm Spin}(32)/ \Zee_2 $) and $ N={\rm dim}(G) = 496 $.
\end{lem} 
\begin{proof}
Let $ r $ be a root of $G$. One can also look at $r$ as a particular weight 
$ r \colon T \rightarrow U(1) $ where $T$ is a choice of maximal torus $T \subset G$. 
For each flat $T$-connection $A$ one can use the weight to associate a flat hermitian 
line bundle, or a holomorphic line bundle. This association leads to an analytic 
projection:
$$ \sigma_r \colon \m{M}_{E, \ T} \rightarrow \m{M}_{E, U(1)}. $$
Employing coordinate models one can see this map can be seen as being factored 
from:
$$ \w{\sigma}_r \colon V_G \rightarrow \m{H} \oplus \mathbb{C}, \ \ \ 
\w{\sigma}_r ( \tau, z) = ( \tau, \ <r^{\vee}, z>) .$$ 
In this setting, the adjoint determinant appears as a tensor product over all roots of $G$ 
(including the zero ones):
\begin{equation}
\label{detnu}
\w{{\rm Det}}(\m{D}_{{\rm ad}}) = \bigotimes_r \ \sigma_r ^* 
{\rm Det}(\m{D}_L). 
\end{equation}    
After factoring out the action of the Weyl group, the line bundle $ (\ref{detnu}) $ 
descends to the determinant fibration $ {\rm Det}(\m{D}_{{\rm ad}}) \rightarrow \m{M}_{E, \ G} $. 
Lemma $ \ref{flatdetlem} $ provides a set of automorphy factors $\varphi_g^{r}$, $g \in \Pi_G$ 
for each $ \sigma_r ^* {\rm Det}(\m{D}_L) $. These are defined on generators of $\Pi_G$ as:
\begin{enumerate}
\item $ \varphi^r_g(\tau, z) = (-1)^{\left ( <r^{\vee}, q_1>+<r^{\vee}, q_2> \right ) }  \ 
 e^{\pi i(- 2<r^{\vee}, z><r^{\vee}, q_2> - 
\tau <r^{\vee}, q_2>^2)} $ for $g \in \m{T}_G$ associated to $ (q_1, q_2) \in \Lambda_G \oplus \Lambda_G $ 
\item $ \varphi^r_g(\tau, z) = \chi(m)^2 \ e^{\pi i \ 
\frac{c<r^{\vee}, z>^2}{c \tau + d}} $ 
for $g \in \m{S}_G$ associated to $ m = \abcd \in {\rm SL}(2, \Zee) $.
\end{enumerate} 
One can then construct a set of automorphy factors $\w{\varphi}_g^{{\rm ad}}$, $g \in \Pi_G$ 
for the adjoint determinant by just taking the product of above quantities over all roots of $G$. 
An useful shortcut here 
is provided by the formula:
\begin{equation}
\sum_{r} \ <r^{\vee}, z_1><r^{\vee}, z_2> = 2c_G  <z_1, z_2>.
\end{equation}   
One obtains therefore that $\w{\varphi}_g^{{\rm ad}}$ can be given on generators of $ \Pi_G$ as:
\begin{enumerate}
\item $ \w{\varphi}^{{\rm ad}}_g(\tau, z) 
= e^{2 c_G \pi i(- 2<z, q_2> - 
\tau <q_2, q_2>)} $ for 
$g \in \m{T}_G$ associated to $ (q_1, q_2) \in \Lambda_G \oplus \Lambda_G $ 
\item $ \w{\varphi}^{{\rm ad}}_g(\tau, z) = \chi(m)^{2N} \ e^{2 c_G \pi i \ 
\frac{c<z, z>^2}{c \tau + d}} $ 
for $ g \in \m{S}_G $ associated to $ m = \abcd \in {\rm SL}(2, \Zee) $.
\end{enumerate}
It is straightforward then that a set of automorphy factors $\varphi_g^{{\rm ad}}$ for 
the adjoint pfaffian can be obtained by factoring the exponents of $\w{\varphi}_g^{{\rm ad}}$ 
by two. The Weyl action factors are constant $1$ due to the Weyl invariance of the
pfaffian section.
\end{proof} 
\vspace{.1in}
\noindent As opposed to the adjoint pfaffian, the spin pfaffian line bundle appearing in fermionic combination
$ (\ref{phopfaf}) $ corresponds to a lower rank representation. However,
such representations do not exist on $G$ and therefore, in order to lower the rank one has to reduce the 
bundle group. For the sake of clarity, let us assume that we analyze the case $G= \eeight $.
Set then $ H = {\rm Spin}(16) \times{\rm Spin}(16) $. There exists then a Lie group homomorphism:
$ \alpha \colon {\rm Spin}(16) \rightarrow E_8 $ which doubles to 
$ \alpha \times \alpha \colon H \rightarrow G $. Any flat $H$-connection
induces then an associated flat $G$-connection and this correspondence
produces an analytic morphism between the respective moduli spaces:
\begin{equation}
\label{mapmod}
\sigma \colon \m{M}_{E,  H} \rightarrow \m{M}_{E,  G}. 
\end{equation}
This map can be explicitly realized in coordinate models. Let $ x_1, \ x_2,  \cdots \ x_8 $ be the fundamental (co)weights 
of $ {\rm SO}(16) $. 
The (co)roots of $ {\rm Spin}(16) $ can then be seen as:
$$ \frac{1}{2}  \left ( \pm x_i \pm x_j \right ) , 
\ \ 1 \leq i < j \leq 8 $$
and they generate the (co)root lattice $\Lambda_{{\rm Spin}(16)}$ of ${\rm Spin}(16)$. The (co)root 
lattice of $ E_8 $ makes a sub-lattice of $\Lambda_{{\rm Spin}(16)}$. 
The (co)roots are obtained in two series:
$$  \pm x_i \pm x_j, \ \ {\rm these \ are \ from} \ 
{\rm Spin}(16)/ \Zee_2 \ {\rm and \ there \ are} \ 112 \ {\rm of \ them} $$   
$$ \frac{1}{2} \left ( \pm x_1 \pm x_2 \cdots \pm x_8 \right ), \ \ 
{\rm where \ there \ are \ an \ even \ number \ of \ minus \ signs}. $$
The Lie group homomorphism $ \alpha \colon {\rm Spin}(16) \rightarrow E_8 $ 
induces then a morphism of coroot lattices:
$$  \alpha_{\Lambda} \colon \Lambda_{{\rm Spin}(16)} \rightarrow 
\Lambda_{E_8}, \ \ 
\alpha_{\Lambda}  \left ( \frac{1}{2} \left ( \pm x_i \pm x_j \right ) \right ) = 
\pm x_i \pm x_j.  
$$
This homomorphism, in turn, canonically extends to a homomorphism $ \alpha \colon \Lambda_H \rightarrow 
\Lambda_G $ commuting with the Weyl projection $ p \colon W_H \rightarrow W_G $. One gets therefore 
the map:
\begin{equation}
{\rm id} \oplus \left ( {\rm id} \otimes \alpha \right ) \colon
\Pi_H \backslash \m{H} \times \left ( \mathbb{C} \otimes_{\mathbb{Z}} 
\Lambda_H \right ) \rightarrow 
\Pi_G \backslash \m{H} \times \left ( \mathbb{C} \otimes_{\mathbb{Z}} 
\Lambda_G \right ) 
\end{equation} 
which models exactly the complex variety morphism $ ( \ref{mapmod}) $.
\par As mentioned earlier, unlike $ E_8 $, $ {\rm Spin}(16) $ admits a
representation in rank 16, the spin representation. Taking the direct sum of two such representations,
one obtains a representation $ \rho $ for $H$ in dimension $32$. 
Coupling the Dirac operators on elliptic curves with 
flat connections induced by representation $ \rho $, one obtains a 
$ \rho $-pfaffian $ \mathbb{C} $-fibration:
\begin{equation}
\label{rhopf}
{\rm Pfaff}(\m{D}_{\rho}) \rightarrow \m{M}_{E, \ H}.
\end{equation} 
It turns out that in the coordinate model:
$$ \m{M}_{E, \ H} = \Pi_H \backslash \m{H} \times \left ( \mathbb{C} \otimes_{\mathbb{Z}} 
\Lambda_H \right ) 
$$
we have:
\vspace{.1in}
\begin{lem}
\label{lemrho}
The holomorphic pfaffian fibration $ (\ref{rhopf}) $ can be by factoring out the 
trivial line bundle over the universal cover $ \m{H} \oplus 
\left ( \mathbb{C} \otimes_{\mathbb{Z}} \Lambda_H \right ) $ 
by means of a set of automorphy factors $ \varphi_g^{\rho}$, $g \in \Pi_H$ which can be given 
on generators of $\Pi_H$ as:
\begin{enumerate}
\item $ \varphi^{\rho}_g(\tau, z) = 
e^{ \pi i(- 2<z, q_2> -\tau <q_2, q_2>)} $ for 
$g \in \m{T}_H$ associated to $ (q_1, q_2) \in \Lambda_H \oplus \Lambda_H $
\item $ \varphi^{\rho}_g(\tau, z) = \chi(m)^{16} \ 
e^{ \pi i \ \frac{c<z,z>}{c \tau + d}} $ 
for $g \in \m{S}_H$ associated to $ m = \abcd \in {\rm SL}(2, \Zee) $
\item $ \varphi^{\rho}_g(\tau, z) = 1 $ for $ g \in \m{W}_H $.
\end{enumerate}
\end{lem} 
\vspace{.1in}
\begin{proof}
The representation $\rho$ is determined by thirty two weights on the maximal torus of $H$ 
represented by two copies of $ x_1, \ x_2, \cdots
x_8, $ $ -x_1, \ -x_2, \cdots - x_8 $. Each weight $x$ determines an analytic map:
$$ 
\sigma_x \colon \m{M}_{E, \ T} \rightarrow \m{M}_{E, \ U(1)} 
$$ 
which is can be described at the level of universal covers as:
$$
\m{H} \oplus 
\left ( \mathbb{C} \otimes_{\mathbb{Z}} \Lambda_H \right ) \rightarrow
\m{H} \oplus \mathbb{C}, \ \ \  (\tau, u \otimes \lambda) 
\rightsquigarrow ( \tau, u \cdot x(\lambda) ).
$$
The pull-back of the determinant bundle on $\m{M}_{E, \ U(1)}$ through
the morphism $ \sigma_x $ determines a holomorphic line bundle
$ \sigma^*_x {\rm Det}(\m{D}_L) \rightarrow \m{M}_{E, \ T} $. Tensoring over
all weights we obtain:
$$ \w{{\rm Det}}(\m{D}_{\rho}) = \bigotimes_w \ \sigma^*_x 
{\rm Det}(\m{D}_L) . $$ 
This bundle and its holomorphic section are invariant under the
action of the Weyl group $ W_H $. Factoring out the Weyl action, one obtains 
therefore the $\rho$-determinant fibration $ {\rm Det}(\m{D}_{\rho}) \rightarrow \m{M}_{E, \ H} $ which 
is the square of $(\ref{rhopf})$. Once again, lemma 
$ \ref{flatdetlem} $ can be used to determine a set of automorphy factors $\varphi^x_g$, $g \in \Pi_H$ for 
each $ \sigma^*_x {\rm Det}(\m{D}_L) $. These factors can be read on generators of $ \Pi_H$ as:
\begin{enumerate}
\item $ \varphi^x_g(\tau, z) = (-1)^{x(q_1) + x(q_2)} \ 
e^{\pi i(- 2x(z)x(q_2)-\tau x(q_2)^2)} $ for $g \in \m{T}_H$ associated to
$ (q_1, q_2) \in \Lambda_H \oplus \Lambda_H $ 
\item $ \varphi^x_g(\tau, z) = \chi(m)^2 \ 
e^{\pi i \ \frac{c \cdot x(z)^2}{c \tau + d}} $ for $ g \in \m{S} $ associated to  
for $ m = \abcd \in \psl $.
\end{enumerate}
Taking the product over all weights and using the fact that:
$$ 2 \cdot \sum_{i=1}^{8} x_i(a)\cdot x_i(b) \ = \ <a,b> $$ 
for all $ a, b \in \Lambda_H $, we obtain a set of automorphy factors 
$ \w{\varphi}^{\rho}_g$, $g \in \Pi_H$ for the $\rho $-determinant fibration. They are: 
\begin{enumerate}
\item $ \w{\varphi}^{\rho}_g(\tau, z) = 
e^{ 2 \pi i(- 2<z, q_2> -\tau <q_2, q_2>)} $ for $ g \in \m{T}_H$ associated to
$ (q_1, q_2) \in \Gamma_H \oplus \Gamma_H $
\item $ \w{\varphi}^{\rho}_g(\tau, z) = \chi(m)^{32} \ 
e^{ 2 \pi i \ \frac{cu^2}{c \tau + d}} $ 
for $g \in \m{S}_H$ associated to $ m = \abcd \in {\rm SL}(2, \Zee) $.
\end{enumerate} 
The factor corresponding to the Weyl action is identically 1. The pfaffian fibration is obtained
by grouping together the weights of opposite sign. Therefore, one can obtain a set $ \varphi^{\rho}_g$ of 
automorphy factors for the $\rho$-pfaffian fibration $(\ref{rhopf})$ by just dividing the 
exponents in $ \w{\varphi}^{\rho}_g$ by two. 
\end{proof} 
\noindent There is a similar $ {\rm Spin}(32)/ \Zee_2 $ analog for Lemma $ \ref{lemrho} $. In this case 
one takes $ H $ to be $ {\rm Spin}(32) $ and $ \rho$ to be the spin representation of $ {\rm Spin}(32)$. 
\subsection{Heterotic Fibration = Character Fibration}
We are now in position to finish the proof of theorem $ \ref{main} $. Recall
the basic facts. The moduli space of $G$-heterotic data makes an 18-dimensional complex
variety $ \mhet $ which is the total space of the heterotic
Chern-Simons $ \mathbb{C}^* $-fibration:
\begin{equation}
\label{hetf}
\m{M}_{\rm het} \rightarrow \m{M}_{E, \ G}.
\end{equation}
Under the coordinate description:
$$ \m{M}_{E,  G} \ = \ \Pi_G \backslash \m{H} \times \left ( \mathbb{C} \otimes _{\mathbb{Z}} \Lambda_G \right )  , $$
the holomorphic type of fibration $ (\ref{hetf}) $ is given, equivariantly, 
by the class of a set of automorphy factors $ \varphi^{{\rm het}}_g$, $ g \in \Pi_G $, which describe a way 
to factor $(\ref{hetf}) $ from the trivial line bundle over the universal cover 
$ \m{H} \times \left ( \mathbb{C} \otimes _{\mathbb{Z}} \Lambda_G \right ) $ of the base. 
\par The fermionic tensor products constructed as pfaffian phase combinations in section $\ref{sec4}$
were holomorphic $\mathbb{C}$-fibrations:
$$ \m{L}_{{\rm ad}} \rightarrow \m{M}_{E, \ G} \ \ \ {\rm and} \ \ \
\m{L}_{\rho} \rightarrow \m{M}_{E, \ H} $$ 
obtained as in $ (\ref{newpfaff}) \ (\ref{phopfaf}) $:
\begin{equation}
\label{ano11}
\m{L}_{{\rm ad}} \ = \ {\rm Pfaff}(\m{D}_{{\rm ad}}) \otimes \left ( 
{\rm Pfaff}(\m{D}_g) \right )^{\otimes 224} 
\end{equation}
and 
\begin{equation}
\label{ano22}
\m{L}_{\rho} \ = \ {\rm Pfaff}(\m{D}_{\rho}) \otimes \left ( 
{\rm Pfaff}(\m{D}_g) \right )^{\otimes -8}.  
\end{equation}
\begin{prop} \
\begin{itemize}
\item The adjoint anomaly fibration: 
$$  \m{L}_{{\rm ad}} \rightarrow \m{M}_{E, \ G} = \Pi_G \backslash \m{H} \times \left ( \mathbb{C} \otimes _{\mathbb{Z}} \Lambda_G \right ) $$ 
can be described by a set of automorphy factors $ \phi^{{\rm ad}}_g$, $g \in \Pi_G$ given on generators of 
$ \Pi_G$ as:
\begin{enumerate}
\item $ \phi^{{\rm ad}}_g(\tau, z) = 
e^{c_G \pi i(-2<z, q_2> -\tau <q_2, q_2>)} $ for $ g \in \m{T}_G$ associated to 
$ (q_1, q_2) \in \Lambda_G \oplus \Lambda_G $
\item $ \phi^{{\rm ad}}_g(\tau, z) =  
e^{c_G \pi i \ \frac{c<z,z>}{c \tau + d}} $ 
for $g \in \m{S}_G$ associated to $ m = \abcd \in {\rm SL}(2, \Zee) $.
\item $ \phi^{{\rm ad}}_g(\tau, z) = 1 $ for $ g \in \m{W}_G $. 
\end{enumerate}
For both cases $ G = \eeight $ and $ G = {\rm Spin}(32)/ \Zee_2 $ the Coxeter
number $c_G =30$. 
\item The $ \rho$-anomaly fibration:
$$  \m{L}_{\rho} \rightarrow \m{M}_{E, \ H} = \Pi_H \backslash \m{H} \times \left ( \mathbb{C} \otimes _{\mathbb{Z}} \Lambda_H \right )   $$ 
can be described by a set of automorphy factors $ \phi^{\rho}_g$, $g \in \Pi_H$ given on generators of 
$ \Pi_H$ as:
\begin{enumerate}
\item $ \phi^{\rho}_g(\tau, z) = 
e^{ \pi i(-2<z, q_2> -\tau <q_2, q_2>)} $ for $ g \in \m{T}_H$ associated to
$ (q_1, q_2) \in \Lambda_H \oplus \Lambda_H $
\item $ \phi^{\rho}_g(\tau, z) =  
e^{ \pi i \ \frac{c<z,z>}{c \tau + d}} $ 
for $g \in \m{S}_H$ associated to $ m = \abcd \in {\rm SL}(2, \Zee) $.
\item $ \phi^{\rho}_g(\tau, z) = 1 $ for $ g \in \m{W}_H $.
\end{enumerate}
\end{itemize} 
\end{prop}
\begin{proof}
This is a straightforward computation based on Lemmas $ (\ref{lemad}) $, $ (\ref{lemrho}) $ and
Corollary $ (\ref{coro1}) $.
\end{proof}
\noindent In this setting, after checking the above formulas upon the character bundle automorphy factors of
Proposition $ \ref{kac} $ one concludes:
\begin{rem}
\label{rema}
Under the analytic map $ \sigma \colon \m{M}_{E,  H} \rightarrow \m{M}_{E,  G} $ 
the two fermionic fibrations $ (\ref{ano11})$ and $ (\ref{ano22})$ are related 
to the $ \Lambda_G$-character fibration $ \m{Z} \rightarrow \m{M}_{E,  G} $ of $ (\ref{thetaline})$ 
as follows:
\begin{equation}
\m{L}_{\rho} \ \simeq \ \sigma^* \m{Z}, \ \ \ 
\m{L}_{{\rm ad}} \ \simeq \m{Z}^{\otimes c_G}. 
\end{equation}
\end{rem}    
\noindent This remark we exploit in order to extract information about the automorphy factors 
$\varphi^{{\rm het}}_g$ of $ (\ref{hetf})$. The considerations in section $ \ref{sec4} $ allow us to 
affirm the following two holomorphic fibration isomorphisms:
\begin{enumerate}
\item $ \m{L}_{\rho} \ \simeq \ \sigma^* \m{L}_{{\rm het}}  $  
\item $ \m{L}_{{\rm ad}} \ \simeq \m{L}_{{\rm het}}^{\otimes c_G}. $
\end{enumerate}
These statements can subsequently be translated into the following equivalences of automorphy factors:
\begin{enumerate}
\item  $ \left [ \left  ( \varphi^{{\rm het}}_a \right )_{a \in 
\Pi_G} \right ]^{c_G} \ \sim 
\left ( \phi^{{\rm ad}}_a \right )_{a \in \Pi_G}   $ 
\item  $ \left  ( \varphi^{{\rm het}}_{p(b)} \right )_{b \in \Pi_H} \ 
\sim \left ( \phi^{\rho}_b \right )_{b \in \Pi_H} $ 
\end{enumerate}
The equivalence relation ``$ \sim $'' asserts the factors on the two sides realize the
same group cohomology class in 
$$ H^1(\Pi_G, H^0(\m{O}^*_{ \m{H} \times \left ( \mathbb{C} \otimes  \Lambda_G \right )   } )) $$
in the first assertion or
$$ H^1(\Pi_H, H^0(\m{O}^*_{\m{H} \times \left ( \mathbb{C} \otimes  \Lambda_H \right )}) $$
in the second assertion or, in geometrical terms, the two sets of automorphy factors determine
isomorphic equivariant line bundles on the universal cover of the base space. As defined earlier, 
$ c_G $ represents the Coxeter number of $G$. For both $ G = \eeight $ and 
$ G = {\rm Spin}(32)/ \Zee_2 $ we have $ c_G =30 $.  
\par The two equivalence relations above and Remark $ \ref{rema} $ allows us to conclude that 
that the two sets of automorphy factors $ \varphi^{{\rm het}}_a $ and 
$ \varphi^{{\rm ch }}_a $ are related:
\begin{equation}
\label{finalcount} 
\varphi^{{\rm het}}_a \ \sim f_a \cdot \varphi^{{\rm ch }}_a 
\end{equation} 
where $ f_a \in H^1(\Pi_G, H^0(\m{O}^*_{\m{H} \times \left ( \mathbb{C} \otimes  \Lambda_G \right )} )) $ and 
satisfies:
\begin{enumerate}
\item $ ( f_a )^{30} = 1 \ \ {\rm for \ any} \ a \in \Pi_G. $ 
\item $ f_{p(b)} = 1 \ \ {\rm for \ any} \ b \in \Pi_H. $ 
\end{enumerate}
Since all $ f_a $ are holomorphic functions, the first condition above
implies that they are all constant functions. Therefore, the set $f_a$ can be reduced to a 
group homomorphism $ f \colon \Pi_G \rightarrow \Zee_{30} $ with $ p(\Pi_H) \subset {\rm Ker}(f) $. 
\par We note immediately that $ f_a =0 $ for $ a \in \m{S}_G $ as $ \m{S}_G = p(\m{S}_H)$. Let us 
next analyze the restriction of $f$ on $ \m{W}_G$. In the case $ G = {\rm Spin}(32)/\mathbb{Z}_2$, 
$H= {\rm Spin}(32)$ one has $p(\m{W}_H) = \m{W}_G$ and therefore $f$ vanishes on $\m{W}_G$. 
For the other case, $ G = \eeight$, $H= {\rm Spin}(16) \times {\rm Spin}(16)$, and the image 
$ p(\m{W}_H) $ is a subgroup of index $ 135=3^3 \cdot 5 $ in  $\m{W}_G$. This can be deduced form 
the orders of the Weyl groups of $ {\rm Spin}(16)$ and $ E_8$, whose explicit computations one can find, 
for example, in \cite{bourbaki}. Therefore, the image of the character $f$ restricted to $\m{W}_G$
is either $\{0\}$ or a subgroup of $ \mathbb{Z}_{30} $ of odd order. But $\m{W}_G$ is generated 
by elements of order two. Hence $ {\rm Im}(f)= \{0\}$. We conclude then that, in both cases,
the character $f$ vanishes on the Weyl transformations in $ \m{W}_G$. Let us then continue with 
an analysis of the behavior of $f$ on $ \m{T}_G$. This group is acted upon by $ \m{W}_G$ and, under 
the identification, $ \m{T}_G \simeq \Lambda_G \times \Lambda_G$, this action is the Weyl action. 
Since the Weyl action is transitive on the roots, we conclude that the 
homomorphism $f$ takes the same value on all roots in $ \Lambda_G \times \Lambda_G$. However,
some of the roots are coming from the roots of $H$ and therefore all roots must belong to 
$ {\rm Ker}(f)$. This shows that $ f_a = 0$ for any $ a \in \m{T}_G $ and since we have already 
proved that $f$ vanishes on $ \m{W}_G$ and $ \m{S}_G$, conclude that $f$ is identically zero on $ \Pi_G$.     
\par By $ (\ref{finalcount}) $ the automorphy factors $ \varphi^{{\rm het}}_a $ 
and  $ \varphi^{{\rm ch }}_a $ are equivalent and hence the heterotic $ \mathbb{C}^* $-fibration 
$ (\ref{hetf})$ is holomorphically isomorphic to the $ \Lambda_G$-character fibration 
$ \m{Z} \rightarrow \m{M}_{E,  G} $ of $ (\ref{thetaline})$. This completes the proof of theorem 
$\ref{main}$.
%
%
%
%
\section{Appendix: The Push-Forward Map.}
\label{appendix} 
In this appendix we clarify the push-forward mechanism
for differential cocycles. This feature has been used in section $ \ref{sec4}$ in 
order to construct the Chern-Simons bundle $(\ref{a11})$ associated to a family 
of heterotic parameters. 
\par As mentioned in \cite{f1}, differential cohomology classes can be integrated. 
In other words, there exists a push-forward homomorphism of the following type: 
\begin{equation}
\int_{E} \colon  \check{H}^n(X \times E) \rightarrow  \check{H}
^{n-d}(X) 
\end{equation} 
where X and E are smooth closed manifolds and $ \text{dim}_R E = d$.
Moreover, there are (non-canonical) extensions of this cohomology map 
to push-forward morphisms for non-flat differential cochains: 
\begin{equation}
\label{intmappp}
\int_{E} \colon N \check{C}^n(X \times E) \rightarrow N \check{C}
^{n-d}(X) 
\end{equation} 
Such kind of map integrates non-flat $n$-cochains from $ X \times E $ 
down to non-flat $(n-d)$-cochains on X. And carries certain
properties as diffeomorphism invariance, sensitivity to orientation,
Stokes' theorem and gluing law. The goal of this appendix is to describe in 
detail the construction of $ (\ref{intmappp}) $. 
\par Recall that there is an integration map for differential forms: 
\begin{equation}
\int_{E} \colon \Lambda ^n(X \times E) \rightarrow \Lambda
^{n-d}(X) 
\end{equation} 
which gives rise to a morphism of the de Rham complexes: 
\begin{equation}
\int_{E} \colon \Lambda ^{\bullet}(X \times E) \rightarrow \Lambda
^{\bullet-d}(X) 
\end{equation} 
The induced homomorphism 
\begin{equation}
\int_{E} \colon H_{dR}^n(X \times E) \rightarrow H_{dR}
^{n-d}(X) 
\end{equation} 
represents the cap product with the fundamental homology class
$ [E] \in H^d(E) $. We would like to reproduce the same pattern
for differential cohomology. 
\par The push-forward mechanism, in the framework we are going to use, is 
in some sense a generalization of the higher gerbe connection holonomy, which 
a non-flat n-cocycle defines along a closed embedded $n$-manifold. This has been 
discussed briefly in section $\ref{sec222}$. namely, if one assumes that 
$ X $ reduces to just a point and $ d=n $ then the particular map : 
\begin{equation}
\int_{E} \colon N \check{Z}^n(E) \rightarrow N \check{Z}
^0(\{ {\rm point} \} ) \simeq {\bf R}  
\end{equation} 
returns a real number whose exponential gives precisely the holonomy of the 
n-cocycle along E.  Accordingly, we construct the push-forward mechanism along 
the lines holonomy was defined in section $ \ref{sec222}$.
\par  Let X and E be closed smooth manifolds and the real dimension of E is $d$. Assume that 
X and E are endowed with (contractible) open coverings: 
$$ X = \ \bigcup_{a \in \mathcal{A}} \ U_a ,  \ 
E = \ \bigcup_{b \in \mathcal{B}} \ U_b \  $$
and that the covering of $E$ is admissible, in the sense that it admits 
a subordinated dual cell decomposition. 
This creates a product covering of $ X \times E$:
$$ X \times E =  \bigcup_{(a b) \in \mathcal{A} \times \mathcal{B} } 
\ U_{(a b)} \ \text{with} \ U_{(a b)}= U_a \times U_b. $$ 
As with holonomy, we make a choice of dual cell decomposition for E, 
$ (\Delta_i)_{ i \in I } $, $ \Delta_i $ representing the top d-cells. We
keep the same orientation convention. $ \Delta_{(i_1 > i_2 \cdots
> i_k) } $ is the $ (d+1-k)$-face obtained by intersecting $ \Delta_{i_1} 
, \ \Delta_{i_2} , \ \cdots , \ \Delta_{i_k} $ if such a thing does
exist. The orientation we consider on 
$ \Delta_{(i_1 > i_2 > \cdots > i_k) } $ is the one obtained when
regarding $ \Delta_{(i_1 > i_2 > \cdots > i_k) } $ as a boundary
component in $ \Delta_{(i_1 > i_2 > \cdots > i_{k-1}) } $. As before, 
suppose the decomposition is subordinated to the covering 
$ (U_b)_{b \in \mathcal{B}} $ through a subordination map 
\begin{equation}
\rho \colon I \rightarrow \mathcal{B} , \ \Delta_i \subset U_{\rho(i)}
\end{equation} 
We now describe the push-forward map. Let $ \omega $ be a non-flat
n-cochain on $ X \times E $. We assume that the product covering 
$ U_{(a b)}= U_a \times U_b $ is chosen to be sufficiently small such 
that $ \omega $ can be represented in this covering as a multiplet: 
\begin{equation}
\omega = \ ( \ H, \ \omega^n_{(a_1 b_1)}, \ 
\omega^{n-1}_{(a_1 b_1)(a_2 b_2)}, \ 
\omega^{n-2}_{(a_1 b_1)(a_2 b_2)(a_3 b_3)},  \ \cdots , \ 
\omega^{-1}_{(a_1 b_1)(a_2 b_2) \cdots (a_{n+2} b_{n+2})} \ )
\end{equation}  
The integration map associates then to each $ \omega $ a
non-flat (n-d)-cochain: 
$$ A \ = \ \int_{E} \omega \ \in N \check{C}^{n-d} (X). $$ 
and such an object can be represented in the covering 
$ (U_a)_{a \in \mathcal{A}} $ as a multiplet: 
\begin{equation}
A \ = \ ( \ T, \ A^{n-d}_{a_1}, \ A^{n-d-1}_{a_1 a_2}, \ 
\cdots , \ A^{-1}_{a_1 a_2 \cdots a_{n-d+2}} \ ).
\end{equation} 
\begin{dfn}
\label{integmech} 
The components T and $ A_{(a)} $ are defined as:  
\begin{equation}
T= \int_{E} H 
\end{equation} 
and for higher indices $ (a) = (a_1 a_2 \cdots a_r) $ , 
\begin{equation}
\label{mess}
A^{n-d+1-r}_{(a)} \ = \sum_{k=1}^{d+1} \ \ 
(-1)^{(n-d+1)(k+1)}  
\sum_{(i)=(i_1 > i_2 > \cdots > i_k)} \ 
\int_{\Delta_{(i)}} \ T^{(a)}_{\rho((i))} \omega
\end{equation} 
if $ 1 \leq r \leq n+1-d $ and  
\begin{equation}
\label{expr2}
A^{-1}_{(a)} = (-1)^{(n-d+1)d}  
\sum_{(i)=(i_1 > i_2 > \cdots > i_{d+1})} \ 
\int_{\Delta_{(i)}} \ T^{(a)}_{\rho((i))} \omega 
\end{equation} 
for $ r = n+2-d $.  
For a multi-index $ (i)=(i_1 > i_2 > \cdots > i_k)  $ in the second sum, 
$ \rho((i)) $ is defined as $ ( \rho(i_1) \rho(i_2) \cdots \rho(i_k) ) $. 
The symbols: 
\begin{equation}
\label{symbt}
T^{(a_1 a_2 \cdots a_r)}_{(b_1 b_2 \cdots b_k)} \omega
\end{equation} 
represent local $(n+2-r-k)$-forms (local constant functions with values 
in $ 2 \pi \Zee $ if $ n+2-r-k=-1 $) living on 
$$ (U_{a_1} \cap U_{a_2} \cap \cdots \cap U_{a_r}) \times  
(U_{b_1} \cap U_{b_2} \cap \cdots \cap U_{b_k} ).  $$  
The exact formulation of symbols $ (\ref{symbt})$ is given by the next definition. 
By convention, in expressions $ (\ref{mess}) $ and $ (\ref{expr2}) $, 
integration over points 
($ \Delta_{(i)} $ with $ \vert (i) \vert = d+1 $) means just evaluating 
the respective forms on $ \Delta_{(i)} $. 
\end{dfn} 
\begin{dfn}
\label{d1}
Let the two multi-indices be:
$$ (a) = ( a_1, a_2, \cdots , a_r) \ {\rm and} \ 
(b) =  (b_1, b_2, \cdots , b_k) . $$
We set 
\begin{equation}
T ^{(a)} _{(b)} \omega = \
\sum_{\gamma \in \mathcal{D}} \ (-1)^{A(d)} \ \omega^{n+2-k-r} _d 
\end{equation} 
where $ \mathcal{D} $ stands for the set of paths in the rectangular
network $ (a_p , b_q) $ generated by 
$ (a)=( a_1, a_2, \cdots , a_r) $ and $ (b) = (b_1, b_2,\cdots ,b_k) $ 
joining $ (a_1, \ b_1) $ to $ (a_r, b_k) $ and moving only to
the right or upward. $ A(\gamma) $ stands for the area of the domain bounded 
by the path $ \gamma $  and the b-axis. For 
$$ \gamma = \left ( \ (a_{p_1} \ b_{q_1}),  (a_{p_2} \ b_{q_2}), \cdots
(a_{p_{r+k-1}} b_{q_{t+r-1}}) \ \right )   
$$
we define  
\begin{equation}
\label{symbols}
\omega^{n+2-k-r}_{\gamma} = 
\omega^{n+2-r-k}_{ (a_{p_1} \ b_{q_1}) (a_{p_2} \ b_{q_2}) \cdots
(a_{p_{r+k-1}}  b_{q_{k+r-1}}) }.
\end{equation} 
\end{dfn} 
\noindent We make the remark here that in $ (\ref{symbols}) $, if the number of nodes in 
the path $ \gamma $ is $ n+2 $ then 
$ \omega^{-1}_{\gamma} \in 2 \pi \Zee $. Relation $ (\ref{expr2}) $ let us
then conclude that $ A^{-1}_{(a)} \in 2 \pi \Zee $ for $ \vert (a) \vert =
n-d+2 $, and therefore, the multiplet $ A $ described in definition $ \ref{integmech} $ 
is indeed a differential cochain. One obtains therefore a push-forward 
map: 
\begin{equation}
\int_{E} \colon N \check{C}^n(X \times E) \rightarrow N \check{C}
^{n-d}(X) 
\end{equation} 
It satisfies the following propriety: 
\begin{teo}
\label{t1}
For every non-flat differential 
cocycle $ \omega \in N \check{Z}^N (X \times E) $, 
the push-forward: 
$$ \int_{E} \omega $$ 
is a (n-d)-dimensional non-flat differential cocycle on X. 
\end{teo} 
\noindent Moreover, theorem $ \ref{t1} $  is just a particular case of a 
Stokes-like integration argument.  
\begin{teo}
\label{t2}
The push forward map defined above for differential cochains commutes with
the differentiation operator $ \check{d} $. Namely, 
\begin{equation}
\int_{E} \check{d} \omega \ = \ \check{d} \left ( \ \int_{E} 
\omega \ \right )  
\end{equation} 
\end{teo}  
\noindent Based on above statement, one concludes that the integration mechanism 
constructed so far sends differential cocycles to differential cocycles
and coboundaries of flat differential cochains to coboundaries of flat
cochains. This induces therefore a push-forward cohomology map: 
\begin{equation} 
\int_{E} \colon \check{H}^n (X \times E) \rightarrow 
\check{H}^{n-d}(X).  
\end{equation}  
\vspace{.1in}
\begin{proof} (of Theorem $ \ref{t2}$) Theorem $ \ref{t1} $  will 
follow as a corollary. Assume that the input $ \omega $ is a non-flat differential 
n-cochain on a product open covering of $ X \times E $ 
and $ \eta = \check{d} \omega \in N \check{C}^{n+1} ( X \times E) $. 
We are going to show that the push-forward cochains:  
\begin{equation}
A \ = \ \int_{E} \omega \ = \ 
( \ T, \ A^{n-d}_{a_1}, \ A^{n-d-1}_{a_1 a_2}, \ 
\cdots , \ A^{-1}_{a_1 a_2 \cdots a_{n-d+2}} \ )
\end{equation} 
and 
\begin{equation}
B \ = \ \int_{E} \eta \ = \ 
( \ Q, \ B^{n-d}_{a_1}, \ B^{n-d-1}_{a_1 a_2}, \ 
\cdots , \ B^{-1}_{a_1 a_2 \cdots a_{n-d+2}} \ )
\end{equation} 
satisfy $ \check{d} A = B $. That means: 
$$ dT=Q, \ \ T_{\vert U_{a_1}}-dA_{a_1}=B_{a_1}  \ \ \mbox{and} \ 
(\delta A)_{(a)}+(-1)^{\vert (a) \vert} \ dA_{(a)} = \ B_{(a)} \ . $$ 
The first identity can be quickly proved: 
$$
dT \ = \ d \left ( \int_{E} H \right ) \ = \ \int_{E} d_x H \ = \ 
\int_{E} dH \ - \ \int_{E} d_y H \ = \ \int_{E} dH \ = \ Q.  
$$ 
Here $ d_x $ and $ d_y $ are derivatives along $ X$, respectively $ Y $ 
direction. 
\par Let us prove that $ T_{\vert U_{a_1}}-dA_{a_1}=B_{a_1} $. We have: 
$$ T-dA_{a_1}= \int_{E} H - dA_{a_1} = $$ 
$$ 
= \int_{E} H - 
\left ( 
\sum_{k=1}^{d+1} \ \ (-1)^{(n-d+1)(k+1)} 
\sum_{(i)=(i_1 > i_2 > \cdots > i_k)}  
\int_{\Delta_{(i)}} \ d_x T^{a_1}
_{(\rho(i_1)\rho(i_2) \cdots \rho(i_k))} \omega \right )= $$
$$ = \ \sum_{i_1} \int_{\Delta_{i_1}} H \ - $$
$$
- \ \sum_{k=1}^{d+1} (-1)^{(n-d+1)(k+1)} 
\sum_{(i)=(i_1 > i_2 > \cdots > i_k)}  
\int_{\Delta_{(i)}}  ( d  T^{a_1}
_{(\rho(i_1)\rho(i_2) \cdots \rho(i_k))} \omega
-  d_yT^{a_1}
_{(\rho(i_1)\rho(i_2) \cdots \rho(i_k))} \omega ). $$ 
But Stockes' Theorem provides: 
$$ \sum_{(i)=(i_1 > i_2 > \cdots > i_{k-1})}  
\int_{\Delta_{(i)}} d_y  T^{a_1}
_{(\rho(i_1)\rho(i_2) \cdots \rho(i_{k-1}))} \omega= $$  
$$
= \sum_{(i)=(i_1 > i_2 > \cdots > i_k)}  
(-1)^{n-d+1+k+1} \ \int_{\Delta_{(i)}} [ \delta_b T^{a_1}_{(.)} \omega]
_{(\rho(i_1)\rho(i_2)) \cdots \rho(i_k))}.  $$ 
Here $ \delta_b $ means Cech differentiation with respect to the
second indices. Continuing along these lines one obtains that
$ T-d A_{a_1} $ equals: 
$$ \sum_{i_1} \int_{\Delta_{i_1}} (H -d T^{a_1}_{\rho(i_1)} \omega ) \ \ 
- \ \sum_{k=2}^{d+1} (-1)^{(n-d+1)(k+1)} $$
$$ \sum_{(i)=(i_1 > i_2 > \cdots > i_k)}   
\left ( \int_{\Delta_{(i)}} ( d T^{a_1}
_{(\rho(i_1)\rho(i_2) \cdots \rho(i_k))} 
+  (-1)^k \ 
[ \delta_b T^{a_1}_{(.)} \omega]
_{(\rho(i_1)\rho(i_2)) \cdots \rho(i_k))} ) 
\right ). $$ 
However, $ \check{d} \omega = \eta $ and therefore: 
$$ H -d T^{a_1}_{\rho(i_1)} \omega \ = \ 
H - d \omega^{n}_{(a_1 \rho(i_1))} \ = \ \eta^{n+1}_{(a_1 \rho(i_1))} \
 = \ T^{a_1}_{\rho(i_1)} \eta $$ 
and 
$$  d T^{a_1}
_{(\rho(i_1)\rho(i_2) \cdots \rho(i_k))} \omega \ 
+  (-1)^k \ 
[ \delta_b T^{a_1}_{(.)} \omega]
_{(\rho(i_1)\rho(i_2)) \cdots \rho(i_k))} \ = $$ 
$$ \ 
d \omega^{n+1-k}_{(a \rho(i_1))(a \rho(i_2)) \cdots 
(a \rho(i_k))} \ +  (-1)^k \  
[\delta \omega^{n+1-k}]_{(a \rho(i_1))(a \rho(i_2)) \cdots 
(a \rho(i_k))} \ = 
$$ 
$$
\ = \ (-1)^{k} \ \eta^{n+2-k}_{(a \rho(i_1))(a \rho(i_2)) \cdots 
(a \rho(i_k))} \ = \ (-1)^{k} T^{a_1}
_{(\rho(i_1)\rho(i_2) \cdots \rho(i_k))} \eta.  
$$  
\par Hence: 
$$ 
T-d A_{a_1} \ =  \ 
$$
$$ 
= \sum_{i_1} \int_{\Delta_{i_1}} T^{a_1}_{\rho(i_1)} \eta + 
\sum_{k=2}^{d+1} (-1)^{(n-d+1)(k+1)} 
\sum_{(i)=(i_1 > i_2 > \cdots > i_k)}   
\int_{\Delta_{(i)}}  (-1)^{k} T^{a_1}
_{(\rho(i_1)\rho(i_2) \cdots \rho(i_k))} \eta \ = \  
$$ 
$$ 
= \ \sum_{k=1}^{d+1} (-1)^{(n-d+2)(k+1)} 
\sum_{(i)=(i_1 > i_2 > \cdots > i_k)}   
\int_{\Delta_{(i)}}  T^{a_1}
_{(\rho(i_1)\rho(i_2) \cdots \rho(i_k))} \eta \ \ =  \ B_{a_1}. 
$$ 
\vspace{.1in}
We are now in position to prove the cocycle condition for higher indices. 
Namely: 
$$ [\delta A]_{(a)}+(-1)^{\vert (a) \vert} \ dA_{(a)} \ = \  B_{(a)} \ \ 
\mbox{for} \ (a)=(a_1 a_2 \cdots a_r), \ \vert (a) \vert = r \leq n-d+2. $$ 
We have: 
$$
dA_{(a)} \ = \sum_{k=1}^{d+1} \ \ 
(-1)^{(n-d+1)(k+1)}  
\sum_{(i)=(i_1 > i_2 > \cdots > i_k)} \ 
\int_{\Delta_{(i)}} \ d_xT^{(a)}_{\rho((i))} \omega 
$$ 
and 
$$  
(\delta A)_{(a)} \ = \sum_{k=1}^{d+1} \ \ 
(-1)^{(n-d+1)(k+1)}  
\sum_{(i)=(i_1 > i_2 > \cdots > i_k)} \ 
\int_{\Delta_{(i)}} \ [ \delta_a T^{(.)}_{\rho((i))} \omega]^{(a)}.
$$ 
Taking into account 
$$ d_xT^{(a)}_{\rho((i))} \omega \ = \  d T^{(a)}_{\rho((i))} \omega \ -
\   d_y T^{(a)}_{\rho((i))} \omega $$ 
we write: 
$$   
[\delta A]_{(a)}+(-1)^{r} \ dA_{(a)} \ = \ $$ 
$$ 
= \sum_{k=1}^{d+1} \ \ 
(-1)^{(n-d+1)(k+1)}  
\sum_{(i)=(i_1 > i_2 > \cdots > i_k)} \ 
\int_{\Delta_{(i)}} 
\left ( 
[ \delta_a T^{(.)}_{\rho((i))} \omega]^{(a)} \ + \
(-1)^{r} d_x T^{(a)}_{\rho((i))} \omega 
\right ) \ = $$
$$ 
= \  \sum_{k=1}^{d+1} \ \ 
(-1)^{(n-d+1)(k+1)}  
\sum_{(i)=(i_1 > i_2 > \cdots > i_k)} \ $$
$$ 
\int_{\Delta_{(i)}} \
\left ( 
[ \delta_a T^{(.)}_{\rho((i))} \omega]^{(a)} \ + \
(-1)^{r} d T^{(a)}_{\rho((i))} \omega \ - \ 
(-1)^{r} d_y T^{(a)}_{\rho((i))} \omega 
\right ). \ 
$$ 
But 
$$ 
\sum_{(i)=(i_1 > i_2 > \cdots > i_{k-1})} \
\int_{\Delta_{(i)}} \ d_y T^{(a)}_{\rho((i))} \omega \ = \ 
(-1)^{n-d+r+k+1} \    
\sum_{(i)=(i_1 > i_2 > \cdots > i_k)} \
\int_{\Delta_{(i)}} \ [\delta_b T^{(a)}_{(.)} \omega]_{\rho((i))}. 
$$ 
Therefore we continue: 
\begin{equation}
\label{eq2027}
[\delta A]_{(a)}+(-1)^{r} \ dA_{(a)} \ = \ 
\end{equation} 
$$ 
= \  \sum_{i_1} \ \int_{\Delta_{i_1}}  
\left ( [ \delta_a T^{(.)}_{\rho(i_1)} \omega]^{(a)} \ + \ 
(-1)^{r} d T^{(a)}_{\rho(i_1)} \omega \right ) \ + \
\sum_{k=2}^{d+1} \ \ 
(-1)^{(n-d+1)(k+1)}  
\sum_{(i)=(i_1 > i_2 > \cdots > i_k)} \ $$
$$ 
\int_{\Delta_{(i)}} \
\left ( 
[ \delta_a T^{(.)}_{\rho((i))} \omega]^{(a)} \ + \
(-1)^{r} d T^{(a)}_{\rho((i))} \omega \ + \ 
(-1)^{k+1} [\delta_b T^{(a)}_{(.)} \omega ] _{\rho((i))} 
\right ). \ 
$$  
\noindent We claim now that the right part of identity $ (\ref{eq2027}) $ is 
exactly $ B_{(a)} $. 
To clarify this assertion we need the following lemma: 
\begin{lem}
\label{l1}
Let $ \omega \in N \check{C}^n (X \times E) $ be a non-flat differential
n-cocycle and 
$ \check{d} \omega = \eta \in N \check{C}^{n+1} (X \times E) $. We consider
the symbols $ T^{(a)}_{(b)} \omega $ and $ T^{(a)}_{(b)} \eta $ introduced 
by Definition $ \ref{d1} $ . They satisfy the following relation: 
\begin{equation}
\label{eq2026}
T^{(a)}_{(b)} \eta  \ = \ 
(-1)^{\vert (a) \vert + \vert (b) \vert +1 } \ d T^{(a)}_{(b)} \omega \ + \
(-1)^{\vert (b) \vert +1 } \ [\delta_a T^{(.)}_{(b)} \omega]^{(a)} \ + \
[\delta_b T^{(a)}_{(.)} \omega ] _{(b)}.
\end{equation}  
If $ \vert (a) \vert = 1 $ or $ \vert (b) \vert = 1 $  the second or
the third term in $ (\ref{eq2026}) $ disappears.
\end{lem}
\noindent Before giving a proof for Lemma $ \ref{l1} $ , let us notice that this 
this ends the proof of theorem $ \ref{t2} $ . Indeed, one can see 
that the terms 
appearing in summation $ (\ref{eq2027}) $ can be immediately rewritten as: \\
\begin{equation}
[ \delta_a T^{(.)}_{\rho(i_1)} \omega]^{(a)} \ + \ 
(-1)^{r} d T^{(a)}_{\rho(i_1)} \omega \ = \ T^{(a)}_{\rho(i_1)} \eta
\end{equation} 
and
\begin{equation}
[ \delta_a T^{(.)}_{\rho((i))} \omega]^{(a)} \ + \
(-1)^{r} d T^{(a)}_{\rho((i))} \omega \ + \ 
(-1)^{k+1} [\delta_b T^{(a)}_{(.)} \omega ] _{\rho((i))} \ = \
(-1)^{k+1} T^{(a)}_{\rho((i))} \eta .
\end{equation} 
Therefore, $ [\delta A]_{(a)}+(-1)^{r} \ dA_{(a)} \ = $ 
$$
= \ \sum_{i_1} \ \int_{\Delta_{i_1}} T^{(a)}_{\rho(i_1)} \eta \ + \ 
\sum_{k=2}^{d+1} \ \ 
(-1)^{(n-d+1)(k+1)}  
\sum_{(i)=(i_1 > i_2 > \cdots > i_k)} \ \int_{\Delta_{(i)}} \
(-1)^{k+1} T^{(a)}_{\rho((i))} \eta \ = \ 
$$
$$ 
= \ \sum_{k=1}^{d+1} \ \ 
(-1)^{(n-d+2)(k+1)}  
\sum_{(i)=(i_1 > i_2 > \cdots > i_k)} \ \int_{\Delta_{(i)}} \
T^{(a)}_{\rho((i))} \eta \ = \ B_{(a)}. 
$$ 
This finishes the proof. 
\end{proof}
\noindent We still have to provide a proof for lemma $ \ref{l1} $ that was used above.
\begin{proof} (of Lemma $ \ref{l1}$.) Let us recall the 
definition of the symbols $ T^{(a)}_{(b)} $. 
For $ (a) = ( a_1, a_2, \cdots , a_r) $ and
$ (b) = (b_1, b_2, \cdots , b_k) $, these objects are expressed as: 
\begin{equation}
T ^{(a)} _{(b)} \omega = \
\sum_{\gamma \in \mathcal{D}} \ (-1)^{A(\gamma)} \ \omega^{n+2-k-r} _{\gamma} 
\end{equation} 
where $ \mathcal{D} $ is the set of paths  
$$ \gamma = \left ( \ (a_{p_1} \ b_{q_1}),  (a_{p_2} \ b_{q_2}), \cdots
(a_{p_{r+k-1}} b_{q_{r+k-1}}) \ \right )  
$$ 
in the rectangular network 
$$  
\{ (a_i \  b_j) \ \vert \ i \in \{ 1 \cdots r \} , \ j \in \{1 \cdots k \} \}
$$ 
starting at $ (a_1 \ b_1) $ , ending up at $ (a_r \ b_k) $ and moving only
upward or to the right. $ A(\gamma) $ represents the area of the domain 
bounded by the 
path $ \gamma $ and the b-axis. With this in mind: 
$$ 
T^{(a)}_{(b)} \eta \ = \ \sum_{\gamma \in \mathcal{D}} \ (-1)^{A(d)} \ 
\eta^{n+2-k-r}_{\gamma} \ = \ 
$$ 
$$
= \ \sum_{\gamma \in \mathcal{D}} \ (-1)^{A(\gamma)} \ 
\left ( \ 
(\delta \omega^{n+2-k-r})_{\gamma} \ + \ (-1)^{r+k-1} \ 
d \omega^{n+1-k-r}_{\gamma} \  
\right )  \ =   
$$  
$$
= \ \sum_{\gamma \in \mathcal{D}} \ (-1)^{A(\gamma)} \ 
(\delta \omega^{n+2-k-r})_{\gamma}
\ + \ (-1)^{r+k-1} \ 
\left ( 
\sum_{\gamma \in \mathcal{D}} \ (-1)^{A(\gamma)} \ 
\omega^{n+1-k-r}_{\gamma} \  
\right ).
$$ 
\noindent The second term on the left hand side above is exactly 
$ (-1)^{\vert (a) \vert + \vert (b) \vert +1 } \ d T^{(a)}_{(b)} \omega $ 
,so, in order to complete the proof of the lemma we just have to show:  
\begin{equation} 
\label{eq2028}
\sum_{\gamma \in \mathcal{D}} \ (-1)^{A(\gamma)} \ 
(\delta \omega^{n+2-k-r})_{\gamma} \ =
 \ (-1)^{\vert (b) \vert +1 } \ [\delta_a T^{(.)}_{(b)} \omega]^{(a)} \ + \
[\delta_b T^{(a)}_{(.)} \omega ] _{(b)} 
\end{equation} 
\noindent One can prove this by observing that both sides of $ (\ref{eq2028}) $ 
are made of terms of type: 
$$   \pm  \ 
\omega ^{n+3-r-k}_{(a_{p_1} b_{q_1})(a_{p_2} b_{q_2}) \cdots 
(a_{p_{k+r-2}} q_{i_{k+r-2}})} \ .
$$ 
We have to show that, after cancellation, the same terms appear in both
sides, with the same signatures. Let 
$$ \gamma =[ \ (a_{p_1} b_{q_1})(a_{p_2} b_{q_2}) \cdots (a_{p_{k+r-1}} 
b_{q_{k+r-1}}) \ ]  $$ 
 be a path in $ \mathcal{D} $. We have: 
$ (a_{p_1} b_{p_1}) = (a_1 b_1) $ and $ (a_{p_{k+r-1}} b_{q_{k+r-1}})
=(a_r b_k) $. Let $ t \in {1, \cdots \ , k+r-1} $. We denote by
$ \gamma_t $ the broken path obtained after removing the point 
$ (a_{p_t} b_{q_t}) $ from $ \gamma $. In other words: 
$$ 
\gamma_t \ = \ [ \ (a_{p_1} b_{q_1}) \ (a_{p_2} b_{q_2}) \ \cdots 
\ (a_{p_{t-1}} b_{q_{t-1}}) \ (a_{p_{t+1}} b_{q_{t+1}}) \ 
\cdots \ (a_{p_{k+r-1}} b_{q_{k+r-1}}) \ ] \ .  
$$ 
Let 
$$
\omega^{n+3-k-r}_{\gamma_t} \ = 
\ \omega^{n+3-k-r}_{(a_{p_1} b_{q_1}) (a_{p_2} q_{i_2}) \cdots 
(a_{p_{t-1}} b_{q_{t-1}}) (a_{p_{t+1}} b_{q_{t+1}})  
\cdots (a_{p_{k+r-1}} b_{q_{k+r-1}}) } \ .
$$ 
$ \omega^{n+3-k-r}_{\gamma_t} $ appears in the right hand side of 
equation $ (\ref{eq2028}) $ with signature: 
$$
(-1)^{ A(\gamma) +t+1}
$$ 
We analyze what happens on the other side. There are four cases to be 
considered.  
\begin{enumerate}
\item $ b_{q_{t-1}} = b_{q_t} = b_{q_{t+1}} $. This happens if the 
path $ \gamma $ keeps going horizontally in the vicinity of 
the node $ (a_{p_{t}} b_{q_{t}}) $. 
\item $ a_{p_{t-1}} = a_{p_t} = a_{p_{t+1}} $. This happens when the path
keeps climbing vertically in the vicinity of the node 
$ (a_{p_{t}} b_{q_{t}}) $.
\item $ b_{q_{t-1}} = b_{q_t} = b_{q_{t+1}}-1 $. The path $ \gamma $ 
makes a lower L with center in $ (a_{p_{t}} b_{q_{t}}) $. 
\item $ a_{p_{t-1}} = a_{p_t} = a_{p_{t+1}}-1 $. The path makes a
upper L with center in  $ (a_{p_{t}} b_{q_{t}}) $. 
\end{enumerate}
\vspace{.1in}
For each case we get a term 
$$  (-1)^{ A(\gamma)+t+1} \ \omega^{n+3-k-r}_{\gamma_t} \ = \ 
(-1)^{ A(\gamma)+t+1} \
 \omega^{n+3-k-r}_{(a_{p_1} b_{q_1}) (a_{p_2} b_{q_2}) \cdots 
(a_{p_{t-1}} b_{q_{t-1}}) (a_{p_{t+1}} b_{q_{t+1}})  
\cdots (a_{p_{k+r-1}} b_{q_{k+r-1}}) } $$ 
on the right side of the equation. Let us notice that the last two
cases, 3 and 4, cancel each other in some way. For example, if case 3 
happens then the corresponding term 
$$ (-1)^{ A(\gamma) +t+1} \ \omega^{n+3-k-r}_{\gamma_t} $$ 
gets canceled by the term 
$$ (-1)^{ A(\widetilde{\gamma}) +t+1} \ 
\omega^{n+3-k-r}_{\widetilde{\gamma}_t} $$ 
where $ \widetilde{\gamma} $ is the path : 
$$ \widetilde{\gamma} \ = \ 
[ \ (a_{p_1} b_{q_1})(a_{p_2} b_{q_2}) \cdots (a_{p_{t-1}}, b_{q_{t-1}})
(a_{p_{t}-1}, b_{q_{t}+1})(a_{p_{t+1}}, b_{q_{t+1}}) \cdots 
(a_{p_{k+r-1}} b_{q_{k+r-1}}) \ ] $$  
($ \widetilde{\gamma} $ is obtained from $ \gamma $ by downgrading the 
upper L at step t to a lower L). Clearly: 
$$ \omega^{n+3-k-r}_{\widetilde{\gamma}_t} \ = \ 
\omega^{n+3-k-r}_{\gamma_t} $$ 
(as we are removing exactly the step where the two paths differ) and 
since $ A(\widetilde{\gamma}) = A(\gamma) -1 $ we get: 
$$ (-1)^{ A(\gamma) +t+1} \ \omega^{n+3-k-r}_{\gamma_t} \ + \ 
(-1)^{ A(\widetilde{\gamma}) +t+1} \ 
\omega^{n+3-k-r}_{\widetilde{\gamma}_t} = \ 0. $$ 
It is also clear that such terms do not appear on the left side of the
equation $ (\ref{eq2028}) $. 
\par It remains to look at the first two cases. Let us analyze the first case. 
In this situation, the term $ \omega^{n+3-k-r}_{\gamma_t} $ appears in 
$$ ( \delta_b T^{(a_1 a_2 \cdots a_r)}_{(.)} 
\omega )_{(b_1 b_2 \cdots b_k)} $$ 
but not in 
$$ ( \delta_a T^{(.)}_{(b_1 b_2 \cdots b_k)} 
\omega )_{(a_1 a_2 \cdots a_r)} \ . 
$$ 
We check the sign which $ \omega^{n+3-k-r}_{\gamma_t} $ carries on the 
left side of equation $ (\ref{eq2028}) $. 
Assume $ q_t=f $ and $ p_t= g $. Then $ \gamma_t $ can be interpreted
as a path in the rectangular network $ (a_p b_q) $ 
generated by $ (a_1 , \ a_2 , \ \cdots \ a_{r} ) $ and 
$ (b_1 , \ b_2 , \ \cdots \ b_{f-1} , \ b_{f+1} , \ 
\cdots \ b_{f_k} ) $. The area bounded by $ \gamma_t $ and b-axis is: 
$$ A(\gamma_t) \ = A(\gamma) - (\gamma-1) . $$ 
Therefore $ \omega^{n+3-k-r}_{\gamma_t} $ appears in 
$$ T^{(a_1 a_2 \cdots a_{r} )}
_{(b_1 b_2 \cdots b_{f-1} b_{f+1} \cdots b_k )} 
 \ \omega $$ 
with sign 
$$ (-1)^{A(\gamma_t)} \ = (-1)^{A(\gamma)-g+1} \ . $$    
But $ T^{(a_1 a_2 \cdots a_{r} )}
_{(b_1 b_2 \cdots b_{f-1} b_{f+1} \cdots b_k )} 
 \ \omega $ is a term in  $ ( \delta_b T^{(a_1 a_2 \cdots a_r)}_{(.)} 
\omega )_{(b_1 b_2 \cdots b_k)} $ with sign $ (-1)^{f+1} $ and 
$ ( \delta_b T^{(a_1 a_2 \cdots a_r)}_{(.)} 
\omega )_{(b_1 b_2 \cdots b_k)} $ appears on the left hand side of 
expression $ (\ref{eq2028}) $ with positive sign. 
Therefore, the following 
term appears once in the left hand side of $ (\ref{eq2028}) $: 
$$  (-1)^{ A(\gamma)-g +1+f +1} \ \omega^{n+3-k-r}_{\gamma_t} . $$   
Since $ f+g=t+1 $,  $ A(\gamma) -g +1+f +1 \equiv A(\gamma) +t+1 \ 
\text{(mod \ 2)} \ $ and hence, the sign for $ \omega^{n+3-k-r}_{\gamma_t} $ 
in the left hand side of $ (\ref{eq2028}) $ coincides with the 
sign in the right hand side. 
\par One clears the second case in a similar manner. Then, drawing a conclusion, 
the entire right hand side of expression $ (\ref{eq2028}) $ appears exactly the same in the left 
hand side. It
can immediately be checked that the reverse holds too. Each term
in the left hand side of equation $ (\ref{eq2028}) $  appears in the 
mirror on the right hand side. Therefore, the equality holds and the proof of 
the lemma is complete. 
\end{proof} 
\vspace{.1in}
\par The integration mechanism: 
\begin{equation}
\int_{E} \ \omega = A \ = ( \ T, \ A^{n-d}_{a_1}, \ A^{n-d-1}_{a_1 a_2}, \ 
\cdots , \ A^{-1}_{a_1 a_2 \cdots a_{n-d+2}} \ )
\end{equation}
with: 
\begin{equation}
T= \int_{E} H,  
\end{equation} 
\begin{equation}
A_{(a)} \ = \sum_{k=1}^{d+1} \ \ (-1)^{(n-d+1)(k+1)} \    
\sum_{(i)=(i_1 > i_2 > \cdots > i_k)} \ 
\int_{\Delta_{(i)}} \ T^{(a)}_{\rho((i))} \omega \ , \ \ 
\vert (a) \vert \leq n-d+1
\end{equation} 
and 
\begin{equation}
A^{-1}_{a_1 a_2 \cdots a_{n-d+2}} \ = \ (-1)^{(n-d+1)d} \    
\sum_{(i)=(i_1 > i_2 > \cdots > i_k)} \ 
\int_{\Delta_{(i)}} \ 
T^{a_1 a_2 \cdots a_{n-d+2}}_{\rho((i))} \omega \ , 
\end{equation} 
is now complete. However, 
there is one more issue to be discussed. The push-forward map obtained: 
\begin{equation} 
\int_{E} \colon N \check{C}^n (X \times E) \rightarrow 
N \check{C}^{n-d}(X) 
\end{equation} 
commutes with differentiation operator $ \check{d} $, realizing therefore
a morphism of cochain complexes. 
$$
\begin{array}{ccccccccc}
\cdots & \stackrel{\check{d}}{\rightarrow} & N \check{C}^{n-1}(X \times E )  & 
\stackrel{\check{d}}{\rightarrow} & 
N \check{C}^{n}(X \times E )  & 
\stackrel{\check{d}}{\rightarrow} & 
N \check{C}^{n+1}(X \times E )  &
\stackrel{\check{d}}{\rightarrow} & \cdots \\
& &  
\Big\downarrow\vcenter{\rlap{$ \int_E $ }} & &
\Big\downarrow\vcenter{\rlap{$ \int_E $ }} & & 
\Big\downarrow\vcenter{\rlap{$ \int_E $ }} & & \\
\cdots & \stackrel{\check{d}}{\rightarrow} & N \check{C}^{n-d-1}(X)  & 
\stackrel{\check{d}}{\rightarrow} & 
N \check{C}^{n-d}(X)  & 
\stackrel{\check{d}}{\rightarrow} & 
N \check{C}^{n-d+1}(X)  &
\stackrel{\check{d}}{\rightarrow} & \cdots \\
\end{array}
$$
This morphism is far from being canonical. 
It depends on choices of dual cell decomposition $ \Delta_i $ 
of the compact manifold E 
and subordination map $ \rho $. However, one can measure its variation when 
using different sets of choices.
\begin{prop} \
\label{pp1} 
Assume that $ ( \Delta_i, \ \rho ) $ 
and $ ( \Delta'_i , \ \rho' ) $ are two distinct pairs 
of dual cell decompositions and corresponding subordination maps. Let:
$$ \int'_{E} \ \ \text{and} \ \ \int_{E} $$ 
be the associated push-forward morphisms. There exists a homotopy
operator: 
\begin{equation}
k^n \colon N \check{C}^{n}(X \times E ) \rightarrow 
N \check{C}^{n-d-1}(X) 
\end{equation}
such that for any differential non-flat $n$-cocycle 
$ \omega \in N \check{Z}^n (X \times E) $, 
\begin{equation}
\label{ttttt}
\int_{E} \omega \ - \ \int'_{E} \omega \ = \ 
\check{d} k^n (\omega) + k^{n+1} ( \check{d} \omega ).  
\end{equation}
\end{prop} 
\begin{proof}The proof goes in 
two steps.
Initially, we show that the above statement holds if the two cell 
decompositions are the same and only the subordination maps are different. 
Secondly, we prove that refining a decomposition in the same subordination 
map does not change the integration process. 
\par Let us proceed with the first step and construct the homotopy 
operator.  
$$
\begin{array}{ccccccccc}
\cdots & \stackrel{\check{d}}{\rightarrow} & N \check{C}^{n-1}(X \times E )  & 
\stackrel{\check{d}}{\rightarrow} & 
N \check{C}^{n}(X \times E )  & 
\stackrel{\check{d}}{\rightarrow} & 
N \check{C}^{n+1}(X \times E )  &
\stackrel{\check{d}}{\rightarrow} & \cdots \\
& \stackrel{k^{n-1}}{\swarrow  } & 
\Big\downarrow\vcenter{\llap{$ \int'_E \ $ }} 
\Big\downarrow\vcenter{\rlap{$ \int_E $ }} & 
\stackrel{k^{n}}{\swarrow  } &
\Big\downarrow\vcenter{\llap{$ \int'_E \ $ }} 
\Big\downarrow\vcenter{\rlap{$ \int_E $ }} &
\stackrel{k^{n+1}}{\swarrow  } &
\Big\downarrow\vcenter{\llap{$ \int'_E \ $ }} 
\Big\downarrow\vcenter{\rlap{$ \int_E $ }} & & \\
\cdots & \stackrel{\check{d}}{\rightarrow} & N \check{C}^{n-d-1}(X)  & 
\stackrel{\check{d}}{\rightarrow} & 
N \check{C}^{n-d}(X)  & 
\stackrel{\check{d}}{\rightarrow} & 
N \check{C}^{n-d+1}(X)  &
\stackrel{\check{d}}{\rightarrow} & \cdots \\
\end{array}
$$
We keep the notations used earlier.
The subordination maps are defined as: 
$$ \rho, \rho' \colon I \rightarrow \mathcal{B}. $$ 
Let $ \omega $ be a n-cocycle in $ N \check{Z}^n(X \times E) $ described
in the product open covering as: 
\begin{equation}
\omega = ( H, \ \omega^n_{a}, \ \omega^{n-1}_{a_1 a_2}, \ 
\omega^{n-2}_{a_1 a_2 a_3}, \ .... \ \omega^{0}_{a_1 a_2 ....a_{n+1}}, \ 
\omega^{-1}_{a_1 a_2 .... a_{n+2}} ).
\end{equation} 
We define $ \Theta = k^{n}(\omega) \in N \check{C}^{n-d-1}(X) $ 
represented in the open covering $ (U_a)_{a \in \m{A}} $ of $ X $ as:
\begin{equation}
\label{homot}
\Theta \ = \ \left ( \ 0, \Theta^{n-d-1}_{a_1} , \ 
\Theta^{n-d-2}_{a_1 a_2} , \ \cdots \ 
\Theta^{-1}_{a_1 a_2 \cdots a_{n-d+1}} \ \right ). 
\end{equation}
The components of $ (\ref{homot}) $ are obtained as follows. Let 
$ (a) = (a_1 a_2 \cdots a_r) $, $ 1 \leq k \leq n-d+1 $. Then:
\begin{equation}
\Theta^{n-d-r}_{(a)} \ = \ \sum_{k=1}^{d+1} \ (-1)^{(n-d)(k+1)} \ 
\sum_{(i) \ = (i_1 > i_2 > \cdots > i_k)} \ 
\int_{\Delta_{(i)}} \left ( \ \sum_{t=1}^{k} (-1)^{t+1} 
T^{(a)}_{(\rho(i_1) \rho(i_2) \cdots \rho(i_t) \rho'(i_t) \cdots \rho'(i_k))}
\omega \ \right ) .
\end{equation}  
The integrals over points are considered to be restrictions over 
the corresponding $X$-slices. If $ r = n-d+1 $, the above summation
is made for $ k= d+1 $ only. 
\par Denote by $ A_{(a)} $ respectively $ A'_{(a)} $, the
components of the $ (n-d) $-cochains on X obtained from integrating
$ \omega $ using the subordination $ \rho $ respectively 
$ \rho' $. We compute:
$$    
(\delta \Theta)_{(a)} \ + (-1)^{ \vert (a) \vert } \ d \Theta_{(a)} \ = \ 
\sum_{k=1}^{d+1} \ (-1)^{(n-d)(k+1)} \ 
$$
$$ 
\sum_{t=1}^{k} (-1)^{t+1} \sum_{(i)=(i_1 > \cdots > i_k)} \ \int_{\Delta_{(i)}} 
\left ( \ 
[ \delta_a T^{(.)}
_{(\rho(i_1) \cdots \rho(i_t) \rho'(i_t) \cdots \rho'(i_k))} \omega 
]^{(a)} +  (-1)^{\vert (a) \vert } 
d_x T^{(a)}
_{(\rho(i_1) \cdots \rho(i_t) \rho'(i_t) \cdots \rho'(i_k))} \omega \ 
\right ) \ = 
$$  
$$ 
= \ \sum_{k=1}^{d+1} \sum_{t=1}^{k} \ (-1)^{(n-d)(k+1)+t+1} \
\sum_{(i)= (i_1 > \cdots > i_k)} \ 
( \ [ \delta_a T^{(.)}
_{(\rho(i_1) \cdots \rho(i_t) \rho'(i_t) \cdots \rho'(i_k))} \omega 
]^{(a)} \ + 
$$ 
$$ 
+ \ (-1)^{\vert (a) \vert } \ 
d T^{(a)}
_{(\rho(i_1) \cdots \rho(i_t) \rho'(i_t) \cdots \rho'(i_k))} \omega \ 
 + \ \ (-1)^{\vert (a) \vert +1} \ 
d_y T^{(a)}
_{(\rho(i_1) \cdots \rho(i_t) \rho'(i_t) \cdots \rho'(i_k))} \omega \ ).  
$$ 
But, for a fixed $ k $:
$$ \sum_{t=1}^{k} \ \sum_{(j)=(j_1 > \cdots >j_{k-1})} \ (-1)^{t+1} 
\int_{\Delta_{(j)}} \   
d_y T^{(a)}
_{(\rho(j_1) \cdots \rho(j_t) \rho'(i_t) \cdots \rho'(i_{k-1}))} \omega 
\ = 
$$
$$ = \ (-1)^{(n-d+ \vert (a) \vert +k +1)} \ \sum_{t=1}^{k+1} \ 
\sum_{(i)=(i_1 > \cdots > i_k)} 
$$
$$ 
(-1)^{t+1} 
\int_{\Delta_{(i)}} \left (
[ \delta_b T^{(a)}_{(.)} \omega ]
_{(\rho(i_1) \cdots \rho(i_t) \rho'(i_t) \cdots \rho'(i_k))} \ 
+ \ ( T^{(a)}
_{(\rho(i_1) \rho(i_2) \cdots \rho(i_k))} \omega \ - \
T^{(a)}
_{(\rho'(i_1) \rho'(i_2) \cdots \rho'(i_k))} \omega \ )  
\right ). 
$$ 
On can see the above equality by just applying Stokes' Theorem and
carefully removing the canceling terms. Using that in the earlier expression
we get:
$$    
(\delta \Theta)_{(a)} \ + (-1)^{ \vert (a) \vert } \ d \Theta_{(a)} \ = \ 
\ \sum_{k=1}^{d+1} \sum_{t=1}^{k} \ (-1)^{(n-d)(k+1)+t+1}  
$$ 
$$ 
\sum_{(i)= (i_1 > \cdots >i_k)} \ 
( \ [ \delta_a T^{(.)}
_{(\rho(i_1) \cdots \rho(i_t) \rho'(i_t) \cdots \rho'(i_k))} \omega 
]^{(a)} \ +  \ (-1)^{\vert (a) \vert } \ 
d T^{(a)}
_{(\rho(i_1) \cdots \rho(i_t) \rho'(i_t) \cdots \rho'(i_k))} \omega \ + 
$$
$$ 
+ \ \ (-1)^{k} \ 
[ \delta_b T^{(a)}_{(.)} \omega ]
_{(\rho(i_1) \cdots \rho(i_t) \rho'(i_t) \cdots \rho'(i_k))} 
+ \ (-1)^{k+1} \ 
( T^{(a)}
_{(\rho(i_1) \rho(i_2) \cdots \rho(i_k))} \omega \ - \
T^{(a)}
_{(\rho'(i_1) \rho'(i_2) \cdots \rho'(i_k))} \omega \ )
 \ ). 
$$
According to Lemma $ \ref{l1} $, the first three terms inside the above 
summation make:
$$ 
(-1)^{k} \ 
T^{(a)}_{(\rho(i_1) \cdots \rho(i_t) \rho'(i_t) \cdots \rho(i_k))}
\check{d} \omega.
$$ 
Therefore:
$$    
(\delta \Theta)_{(a)} \ + (-1)^{ \vert (a) \vert } \ d \Theta_{(a)} \ = \ 
$$
$$
= - \ \sum_{k=1}^{d+1} \sum_{t=1}^{k} \ (-1)^{(n-d+1)(k+1)+t+1} \ 
\sum_{(i) \ = (i_1 > i_2 > \cdots > i_k)} \ 
\int_{\Delta_{(i)}} \left ( \ 
T^{(a)}_{(\rho(i_1) \rho(i_2) \cdots \rho(i_t) \rho'(i_t) \cdots \rho'(i_k))}
\check{d} \omega \ \right ) +  
$$ 
$$
+ \ \sum_{k=1}^{d+1}  (-1)^{(n-d+1)(k+1)} \
\sum_{(i)=(i_1 > \cdots > i_k)} \ \int_{\Delta_{(i)}} \  
\left ( \ T^{(a)}
_{(\rho(i_1) \rho(i_2) \cdots \rho(i_k))} \omega \ - \
T^{(a)}
_{(\rho'(i_1) \rho'(i_2) \cdots \rho'(i_k))} \omega \ \right )  .
$$
The first sum on the right hand-side of above equality represents 
exactly the $(a)$-component of $ k^{n+1}( \check{d} \omega ) $. 
Therefore we obtain:
$$
(\delta \Theta)_{(a)} \ + (-1)^{ \vert (a) \vert } \ d \Theta_{(a)} \ = 
\ - \left [ k^{n+1}( \check{d} \omega ) \right ] _{(a)} + A_{(a)} - A'_{(a)}, 
$$
In other words:
$$
\left [ \check{d} k^n(\omega) \right ]_{(a)} + 
\left [ k^{n+1}( \check{d} \omega ) \right ] _{(a)} \ = \
A_{(a)} - A'_{(a)}.  
$$
Equality $ (\ref{ttttt}) $ follows. 
The first part of the proof is complete.
\par We continue with the second step. Let us
assume that we are dealing with two different cell decompositions 
$ ( \Delta_i )_{i \in I} $ and $ ( \Delta_j )_{j \in J} $ and the latter
one is a refinement of the former. Say, there is a map: 
$$ \varphi \colon J \rightarrow I  \ \text{such that} \ 
\Delta_j \subset \Delta_{\varphi(j)}. $$ 
Also, assume that both decompositions are subordinated to 
covering $ (U_b)_{b \in \m{B}} $ and the subordination map for $ \Delta_i $ is 
$$ \rho \colon I \rightarrow \mathcal{B}. $$ 
The subordination on 
$ \Delta_j $ can be taken as $ \rho \circ \varphi $. 
\par Say the integration of $ \omega $ using the first pair 
$ ( \Delta_i, \ \rho) $ is : 
\begin{equation}
\int_{E} \ \omega = A \ = ( \ T, \ A^{n-d}_{a_1}, \ A^{n-d-1}_{a_1 a_2}, \ 
\cdots , \ A^{-1}_{a_1 a_2 \cdots a_{n-d+2}} \ )
\end{equation}  
whereas the integration of $ \omega $ using 
$ ( \Delta_j, \ \rho \circ \varphi ) $ is: 
\begin{equation}
\widetilde{\int_{E}} \ \omega = \widetilde{A} \ = ( \ T, \ 
\widetilde{A}^{n-d}_{a_1}, \ \widetilde{A}^{n-d-1}_{a_1 a_2}, \ 
\cdots , \ \widetilde{A}^{-1}_{a_1 a_2 \cdots a_{n-d+2}} \ ).
\end{equation}   
Our goal is to prove that these two multiplets are actually the same. 
We have:
\begin{equation}
\label{eq2030}
A_{(a)} \ = \sum_{k=1}^{d} \ \ (-1)^{(n-d+1)(k+1)} \    
\sum_{(i)=(i_1 > i_2 > \cdots > i_k)} \ 
\int_{\Delta_{(i)}} \ T^{(a)}_{\rho((i))} \omega \ .
\end{equation} 
and 
\begin{equation}
\label{eq2031}
\widetilde{A}_{(a)} \ = \sum_{k=1}^{d} \ \ (-1)^{(n-d+1)(k+1)} \    
\sum_{(j)=(j_1 > j_2 > \cdots > j_k)} \ 
\int_{\Delta_{(j)}} \ T^{(a)}_{\rho(\varphi(j))} \omega \ .
\end{equation} 
(the order relations on J and I are such chosen to
make $ \varphi $ an increasing map). 
\par But, $ \omega_{(c)} $ vanishes as soon as the multi-index $ (c) $ 
includes two identical sub-indices (by convention). Therefore, in
sum $ (\ref{eq2031}) $  all terms: 
$$  \int_{\Delta_{(j)}} \ T^{(a)}_{\rho(\varphi(j))} \omega \ $$ 
for which $ (j) $ contains two sub-indices $ j_u $ and $ j_v $ 
with $ \varphi (j_u) = \varphi (j_v) $, vanish. After removing
those, the terms in expression $ (\ref{eq2031}) $ can be summed up to make
the ones in $ (\ref{eq2030}) $. The second step is complete.
\end{proof}
Let us finish up the appendix by making the following observations related 
to Proposition $ \ref{pp1} $. First, although the integration mechanism we devised 
for differential cochains depends on ingredient choices (dual cell 
decomposition and subordination assignment), the induced cohomology
push-forward map does not depend on those. One can extend the arguments in the 
proof and take into account a variation of the product open covering too. 
Under a change of admissible open covering, the integration mechanism
stays again in the same homotopy class. One can therefore 
conclude that the induced cohomology push-forward map:
$$ \int_{E} \colon  \check{H}^n(X \times E) \rightarrow  \check{H}
^{n-d}(X)  $$
does not depend at all on the choice of product open covering, 
dual cell decomposition and subordination relation. 
\par Second, the construction can be extended to work for 
differential cochains having as derivatives global forms. Namely, let
$ \m{B}^n(X) $ be the set of non-flat differential n-cochains $ \omega $ 
such that $ \check{d} \omega \in \Omega^n(X) $. We know from section 
$ \ref{fir} $
that the flat coboundary equivalence relation extends to $ \m{B}^n(X) $. 
One can
show then that, for $ \omega_1 , \omega_2 \in \m{B}^n(X \times E) $, 
defined on admissible product coverings, with $ \omega_1 \sim \omega_2 $  
$$ \int_E \omega_1 \ \sim \ \int_E \omega _2 \ {\rm in} \ 
\m{B}^{n-d}(X) $$
regardless of the choice of dual cell decomposition subordinated to the
two distinct open coverings. One obtains, therefore, a canonical push-forward
homomorphism:
$$ \int_{E} \colon  \m{B}^n(X \times E) / \sim \ \rightarrow  \m{B}
^{n-d}(X) / \sim.  $$

%
%
%
%
\bibliographystyle{plain}

\end{document}